%%%%%%AMS-TeX file for the paper
%
%     ASKEY-WILSON POLYNOMIALS AND THE QUANTUM $SU(2)$ GROUP:
%     SURVEY AND APPLICATIONS
%
%by H.T. Koelink
%To appear in Acta Applicandae Mathematicae.
%%%%%%%%%%%%%%%%%AMS TeXfile%%%%%%%%%%%%%%%%%%%%%%%%%%%%%%%%%%%%%%%%%%%%
\input amstex.tex
\documentstyle{amsppt}
\magnification=1200
\baselineskip=13pt
\hsize=6.5truein
\vsize=8.9truein
%%%%%%%%%%%%%%%%% M a c r o s %%%%%%%%%%%%%%%%%%%%%%%%%%%%%%%%%%%%%%%%%%%
%Section, Theorem and Formula Numbering%%%%%%%%%%%%%%%%%%%%%%%%%%%%%%%%%%
%%%%%%%%%%%%%%%%%%%%%%%%%%%%%%%%%%%%%%%%%%%%%%%%%%%%%%%%%%%%%%%%%%%%%%%%%
\countdef\sectionno=1
\countdef\eqnumber=10
\countdef\theoremno=11
\countdef\countrefno=12
\countdef\cntsubsecno=13
\sectionno=0
\def\newsection{\global\advance\sectionno by 1
                \global\eqnumber=1
                \global\theoremno=1
                \global\cntsubsecno=0
                \the\sectionno}

\def\newsubsection#1{\global\advance\cntsubsecno by 1
                     \xdef#1{{\S\the\sectionno.\the\cntsubsecno}}
                     \ \the\sectionno.\the\cntsubsecno.}

\def\theoremname#1{\the\sectionno.\the\theoremno
                   \xdef#1{{\the\sectionno.\the\theoremno}}
                   \global\advance\theoremno by 1}

\def\eqname#1{\the\sectionno.\the\eqnumber
              \xdef#1{{\the\sectionno.\the\eqnumber}}
              \global\advance\eqnumber by 1}

\global\countrefno=1

\def\refno#1{\xdef#1{{\the\countrefno}}\global\advance\countrefno by 1}

%               %%%%*     MSSYMB.TeX    %%%%*                  8 Jul 87
%
%       This file contains the definitions for the symbols in the two
%       "extra symbols" fonts created at the American Math. Society.

\catcode`\@=11
%  The following 13 lines establish the use of the Euler Fraktur font.
%  To use this font, remove % from beginning of these lines.
\font\teneuf=eufm10
\font\seveneuf=eufm7
\font\fiveeuf=eufm5
\newfam\euffam
\textfont\euffam=\teneuf
\scriptfont\euffam=\seveneuf
\scriptscriptfont\euffam=\fiveeuf
% Use the next 6 lines with AMS-TeX:
\def\frak{\relaxnext@\ifmmode\let\next\frak@\else
 \def\next{\Err@{Use \string\frak\space only in math mode}}\fi\next}
\def\goth{\relaxnext@\ifmmode\let\next\frak@\else
 \def\next{\Err@{Use \string\goth\space only in math mode}}\fi\next}
\def\frak@#1{{\frak@@{#1}}}
\def\frak@@#1{\noaccents@\fam\euffam#1}
%  End definition of Euler Fraktur font.
\catcode`\@=12

%%%%%%%%%%%%%%%%%%%%%%%%%%%%Abbreviations%%%%%%%%%%%%%%%%%%%%%%%%%%%%%%%%
\def\R{{\Bbb R}}
\def\N{{\Bbb N}}
\def\C{{\Bbb C}}
\def\Z{{\Bbb Z}}
\def\Zp{{\Bbb Z}_+}
\def\hZp{{1\over 2}\Zp}
\def\A{{\Cal A}_q(SU(2))}
\def\Al{{\Cal A}_q^l(SU(2))}
\def\AS{{\Cal A}_{q,\s}(S^2)}
\def\U{{\Cal U}_q({\frak{su}}(2))}
\def\a{\alpha}
\def\b{\beta}
\def\g{\gamma}
\def\d{\delta}
\def\s{\sigma}
\def\t{\tau}
\def\l{\lambda}
\def\tl{\tau_\lambda}
\def\r{\rho_{\t ,\s}}
\def\p{\pi_\theta}
\def\hp{\pi_{\theta/2}}
\def\vp{\varphi}
\def\L{{\Cal L}_\s^l}

%%%%%%%%%%%%%%%%%%%%%%%%%Reference Numbering%%%%%%%%%%%%%%%%%%%%%%%%%

\refno{\Abe}
\refno{\Andr}
\refno{\AndrA}
\refno{\Aske}
\refno{\AskeI}
\refno{\AskeWtwee}
\refno{\AskeW}
\refno{\BadeK}
\refno{\BroctD}
\refno{\BurbK}
\refno{\DijkK}
\refno{\Domb}
\refno{\Drin}
\refno{\ErdeHTF}
\refno{\FlorV}
\refno{\Flor}
\refno{\GaspRPK}
\refno{\GaspR}
\refno{\IsmaW}
\refno{\Jimb}
\refno{\KalnMM}
\refno{\KalnMiMu}
\refno{\Klim}
\refno{\Koel}
\refno{\KoelAF}
\refno{\KoelUP}
\refno{\KoelBeen}
\refno{\KoelBtwee}
\refno{\KoelK}
%\refno{\Koor}
\refno{\Koordrie}
\refno{\KoorOPTA}
\refno{\KoorAF}
\refno{\Koortwee}
\refno{\KoorHG}
\refno{\MasuMNNSUeen}
\refno{\MasuMNNSUtwee}
\refno{\MasuMNNU}
\refno{\NikiU}
\refno{\Noum}
\refno{\NoumM}
\refno{\NoumMCMP}
\refno{\NoumMDMJ}
\refno{\NoumMtwee}
\refno{\NoumMCompM}
\refno{\Podl}
%\refno{\Rahm}
\refno{\RahmV}
\refno{\RideW}
\refno{\Stan}
\refno{\Swee}
\refno{\Vain}
\refno{\VaksK}
\refno{\VaksS}
\refno{\VAsscK}
\refno{\VDael}
\refno{\Vara}
\refno{\Vile}
\refno{\VileK}
\refno{\Woro}
\refno{\Worotwee}
\refno{\Zhed}

%%%%%%%%%%%%%%%%%%%%%%%%%%%%%%%%%%%%%%%%%%%%%%%%%%%%%%%%%%%%%%%%%%%%%**
%%%%%%%%%%%%%%%%%%%%%%%Beginning of Text%%%%%%%%%%%%%%%%%%%%%%%%%%%%%%%
%%%%%%%%%%%%%%%%%%%%%%%%%%%%%%%%%%%%%%%%%%%%%%%%%%%%%%%%%%%%%%%%%%%%%**
\topmatter
\title Askey-Wilson polynomials and the quantum $SU(2)$ group: \\
survey and applications\endtitle
\rightheadtext{Askey-Wilson polynomials and the quantum $SU(2)$ group}
\author H.T. Koelink\endauthor
\affil  Katholieke Universiteit Leuven\endaffil
\address Departement Wiskunde, Katholieke Universiteit Leuven,
Celestijnenlaan 200 B, B-3001 Leuven (Heverlee), Belgium\endaddress
\email erik\%twi\%wis\@cc3.KULeuven.ac.be\endemail
\date to appear in Acta Applicandae Mathematicae\enddate
\thanks Supported by a NATO-Science Fellowship of the Netherlands
Organization for Scientific Research (NWO). \endthanks
\abstract
Generalised matrix elements of the irreducible representations of the
quantum $SU(2)$ group are defined using certain orthonormal bases of the
representation space. The generalised matrix elements are relatively
infinitesimal invariant with respect to Lie algebra like elements of the
quantised universal enveloping algebra of $sl(2)$. A full proof of the
theorem announced by Noumi and Mimachi [Proc. Japan Acad. Sci. {\bf 66},
Ser.~A (1990), pp.~146--149] describing the generalised matrix elements in
terms of the full four-parameter family of Askey-Wilson polynomials is
given. Various known and new applications of this interpretation are
presented.
\endabstract
%Contents
\toc
\widestnumber\head{10}
\head 1. Introduction\endhead
\head 2. Representation theory of $SU(2)$\endhead
\head 3. Preliminaries on basic hypergeometric orthogonal polynomials\endhead
\head 4. Preliminaries on the quantum $SU(2)$ group\endhead
\head 5. Koornwinder's $(\s,\t)$-spherical elements\endhead
\head 6. Generalised matrix elements\endhead
\head 7. Generalised matrix elements and Askey-Wilson polynomials\endhead
\head 8. Some applications\endhead
\head 9. Discrete orthogonality relations\endhead
\head 10. Spherical and associated spherical elements\endhead
\head 11. Characters\endhead
\endtoc
\endtopmatter
\document

%&&&&&&&&&&&&&&&&&&&&&&&&&&&&&&&&&&&&
%%N E W   S E C T I O N%%%%%%%%%%%%**
%%%%%%%%%%%%%%%%%%%%%%%%%%%%%%%%%%%%%

\head\newsection . Introduction
\endhead

The theory of special functions is an old but still lively branch of
mathematics, which has its origins in the works of the great
mathematicians of the 18th and 19th century trying to solve the differential
equations of mathematical physics. One of the highlights in the theory
is the introduction of the hypergeometric series
$$
1+{{ab}\over{1\cdot c}}z + {{a(a+1)b(b+1)}\over{1\cdot 2\cdot
c(c+1)}}z^2 +{{a(a+1)(a+2)b(b+1)(b+2)}\over{1\cdot 2\cdot 3\cdot
c(c+1)(c+2)}}z^3 + \ldots
$$
by Gauss in 1812. It turned out that many special functions can be
written in terms of Gauss's hypergeometric series or in terms of one of
the generalised hypergeometric series.

An example of such a hypergeometric special function is the Jacobi
polynomial which forms a set of orthogonal polynomials, cf. \S 2. Jacobi
polynomials, as many other special functions of hypergeometric type,
satisfy various interesting properties. For the Jacobi polynomials it
was understood in the 1950's by the work of Gelfand and \v Sapiro, cf. the
references in \cite{\Vile}, that some of the properties
of the Jacobi polynomials were reflections of underlying structures
given by the complex group $SL(2,\C)$, or its real form $SU(2)$, to
which we get back in \S 2.
This is just one of the many examples of a very fruitful relation
between group theory, and in particular representation theory of Lie
groups, and special functions of hypergeometric type, which is still a
topic of current research, cf. the survey by Klimyk \cite{\Klim}.

In 1846, 34 years after the introduction of the hypergeometric series by
Gauss, Heine introduced the basic (or $q$-)hypergeometric series
$$
1 + {{(1-q^a)(1-q^b)}\over{(1-q)(1-q^c)}}z +
{{(1-q^a)(1-q^{a+1})(1-q^b)(1-q^{b+1})}\over
{(1-q)(1-q^2)(1-q^c)(1-q^{c+1})}}z^2 + \ldots,
$$
where the base $q$ is usually looked at as a fixed real number between $0$
and $1$. The limit $q\uparrow 1$ of the basic hypergeometric series
yields the hypergeometric series.

Although, since its introduction, work was done on $q$-hypergeometric
series and applications were known in several fields, such as the famous
Rogers-Ramanujan identities in number theory, cf. Andrews \cite{\Andr},
progress was much slower than for special functions of hypergeometric
type. The collaboration of Andrews and Askey starting in the mid 1970's
was the beginning of an outburst of results on $q$-hypergeometric series
and in particular on various $q$-hypergeometric orthogonal polynomials.
This culminated in the introduction in their 1985 memoir by Askey and
Wilson  of a very general four-parameter set of orthogonal polynomials
nowadays known as the Askey-Wilson polynomials.

Until recently there was not much knowledge on possible natural structures
on which the $q$-hypergeometric series could `live' in a similar fashion
as the special functions of hypergeometric type `living' on certain Lie
groups. This problem seems to be solved with the introduction of quantum
groups by Drinfeld \cite{\Drin}, Jimbo \cite{\Jimb} and Woronowicz
\cite{\Woro}, \cite{\Worotwee} around 1986. Quantum groups are
no longer groups, but we consider them as deformations of the algebra of
functions on a group so that the deformed algebra still carries
properties that resemble group actions. More information can be found in
the survey papers by Koornwinder \cite{\KoorOPTA} and Noumi
\cite{\Noum}, see also Klimyk \cite{\Klim} and Vilenkin and Klimyk
\cite{\VileK , Vol.~3, Ch.~14}.

One of the oldest examples of a quantum group is the quantum $SU(2)$
group which is a quantum group analogue of the compact group $SU(2)$.
A first indication of the relation between $q$-hypergeometric functions
and quantum groups is the interpretation of the so-called little
$q$-Jacobi polynomials as matrix elements of irreducible unitary
representations of the quantum $SU(2)$ group by Vaksman and Soibelman
\cite{\VaksS}, Masuda et al. \cite{\MasuMNNU}, Koornwinder
\cite{\Koordrie}. Since the Askey-Wilson polynomials can be viewed as
$q$-analogues of the Jacobi polynomials, see the title of the
Askey-Wilson memoir \cite{\AskeW}, it is natural to try to interpret the
Askey-Wilson polynomials on the quantum $SU(2)$ group.

A decisive step in this direction is taken by Koornwinder
\cite{\Koortwee}. In that paper he shows how to create sufficiently many
`quantum subgroups' of the quantum $SU(2)$ group to get a
two-parameter family of Askey-Wilson polynomials as zonal spherical
functions. We get back to this paper in more detail in \S 5. As
follow-ups to Koornwinder's paper there is the paper \cite{\Koel}
by the author in
which a quantum group theoretic derivation of the Rahman-Verma
\cite{\RahmV} addition formula for the continuous $q$-Legendre
polynomials is given and the announcement by Noumi and Mimachi
\cite{\NoumM}, \cite{\NoumMtwee} in which they claim the interpretation
of the full four-parameter family of Askey-Wilson polynomials as matrix
elements of irreducible unitary representations of the quantum $SU(2)$
group.

This paper grew out of an attempt to provide the announcements by Noumi
and Mimachi \cite{\NoumM}, \cite{\NoumMtwee} with full proofs of the
sort of proofs given in \cite{\Koel} in which the associated
spherical elements on the quantum $SU(2)$ group are calculated
explicitly. The purpose of
this paper is to give a detailed proof of the relation between
Askey-Wilson polynomials and matrix elements of irreducible unitary
representations of the quantum $SU(2)$ group and to present a
survey of some known as well as of some seemingly new applications of
this interpretation of the Askey-Wilson polynomials. In this light the
current paper may be viewed as a sequel to Koornwinder's survey
\cite{\KoorOPTA}.

An alternative to the approach in this paper to $q$-special functions
living on quantum groups is the approach in which representations of the
so-called quantised universal enveloping algebra are $q$-exponentiated
to obtain matrix coefficients expressible as $q$-special functions. This
approach is the analogue of exponentiating a Lie algebra representation
to find a representation of the corresponding Lie group. However, since
the alternative approach uses a $q$-analogue of the exponential function
the matrix coefficients neither end up on some group nor on a quantum
group. Details of this `local' approach can be found in papers by
Floreanini and Vinet \cite{\FlorV}, Kalnins, Manocha and Miller
\cite{\KalnMM}, Zhedanov \cite{\Zhed}, see also the references in
\cite{\FlorV}.

There are two other closely related non-compact quantum groups for which
the relation with $q$-special functions is worked out in some detail. The
quantum $SU(1,1)$ group is obtained from the quantum $SU(2)$ group by
redefining a $\ast$-operation, which means that another real form of the
quantum $SL(2,\C)$ group is chosen. For more information on the
representation theory of the quantum $SU(1,1)$ group and the corresponding
special functions, which are $q$-analogues of Jacobi functions,
the papers by Masuda et al. \cite{\MasuMNNSUeen},
\cite{\MasuMNNSUtwee} and by Burban and Klimyk \cite{\BurbK} may be consulted.
The quantum group of plane motions is obtained from the quantum $SU(2)$ group
by a suitable contraction procedure similar to the classical case. The
representation theory of this quantum group and its connection with basic
Bessel functions is worked out in papers by Vaksman and Korogodski\u\i
\cite{\VaksK} and the author \cite{\KoelBeen}, \cite{\KoelBtwee}.
The alternative `local' approach for these quantum groups
has also been developed, cf. Floreanini and Vinet \cite{\FlorV} and references
therein and Kalnins, Miller and Mukherjee \cite{\KalnMiMu}.

The organisation of the paper is as follows. In \S 2 we consider briefly
the representation theory of the compact Lie group $SU(2)$ and its
relation with special functions. We also discuss the approach to
the group case adapted for deformation. In order to make the paper more
self-contained sections 3 and 4 contain some information on
$q$-hypergeometric orthogonal polynomials and on the quantum $SU(2)$
group. These sections are also used to fix the notation. Sections 3 and
4 are by no means complete, but they present sufficient information
concerning these subjects to be able to read the rest of the paper.
References to the literature are given as well.

In \S 5 we elaborate on Koornwinder's \cite{\Koordrie} approach to
introduce sufficiently many `quantum subgroups' and on the
corresponding spherical elements. In \S 6 we introduce and investigate
generalised matrix elements of the quantum $SU(2)$ group, which are
linked to Askey-Wilson polynomials in \S 7. In \S 8 we discuss some
known applications of the interpretation of Askey-Wilson polynomials on
the quantum $SU(2)$ group, such as addition formulas for $q$-Legendre
polynomials, quantum spheres and non-negative linearisation coefficients
for $q$-Legendre polynomials. The last three sections are
concerned with seemingly new applications. In \S 9 we consider some
discrete orthogonality relations. The relation between spherical
elements and associated spherical elements is studied in \S 10.
Finally, in \S 11 we consider the characters of the irreducible unitary
representations.

\demo{Acknowledgement} I thank George Gasper and especially Mizan Rahman for
their help in checking a large formula in \S 9. Thanks are also due to
Masatoshi Noumi for explaining the results and ideas of \cite{\NoumMCompM},
which led to an improvement of the presentation.
\enddemo

%&&&&&&&&&&&&&&&&&&&&&&&&&&&&&&&&&&&&
%%N E W   S E C T I O N%%%%%%%%%%%%**
%%%%%%%%%%%%%%%%%%%%%%%%%%%%%%%%%%%%%

\head\newsection . Representation theory of $SU(2)$
\endhead

The representation theory of the compact Lie group $SU(2)$ is the
subject of this section. Its relation with special functions is shortly
discussed. The point of view to this group and its representation
theory, which is favourable for the quantum group approach, is given.
This section is of an introductory nature. Complete results on the
relation between $SU(2)$ and special functions can be found in Vilenkin
\cite{\Vile , Ch.~4}, Vilenkin and Klimyk \cite{\VileK , vol.~1, Ch.~6}.
For more information on Lie groups and their representation theory the
reader may consult Br\"ocker and tom Dieck \cite{\BroctD} and
Varadarajan \cite{\Vara}.

The group
$$
SU(2) = \left\lbrace \pmatrix \a &\b\\ \g&\d\endpmatrix \mid
\a\d-\b\g = 1,\ \bar\a=\d,\ \bar\b = - \g\right\rbrace
\tag\eqname{\vgldefinitieSUtwee}
$$
is a compact real three-dimensional Lie group of $2\times 2$ unitary
matrices with determinant $1$. The group $SU(2)$ is a compact real form
of the complex Lie group $SL(2,\C)$, the group of $2\times 2$ complex
matrices with determinant $1$. By Weyl's unitary trick, cf. \cite{\Vara
, \S~4.11} the finite dimensional holomorphic representations of
$SL(2,\C)$ correspond to the finite dimensional representations of
$SU(2)$.

By $K$ we denote the diagonal subgroup
$$
K=S(U(1)\times U(1)) = \left\lbrace \pmatrix \a &0\\ 0&\d\endpmatrix \mid
\a\d = 1,\ \bar\a=\d\right\rbrace \subset SU(2).
$$
So $K$ is equal to the one-dimensional torus, i.e. the unit circle
${\Bbb T}$
in the complex plane viewed as a commutative one-dimensional group with
multiplication given by $e^{i\theta}e^{i\phi}=e^{i(\theta+\phi)}$. The
characters of this commutative group are given by
$$
\psi_n\biggl(\pmatrix \a &0\\ 0&\d\endpmatrix \biggr) = \a^n, \qquad
n\in\Z .
$$

For each $l\in\hZp$ there exists a $(2l+1)$-dimensional irreducible
unitary representation $t^l$, i.e. a homomorphism $t^l\colon SU(2)\to
GL(V^l)$ for a $(2l+1)$-dimensional vector space $V^l$ such that all
operators $t^l(g)$ are unitary and that there is no invariant subspace
of $V^l$ for all operators $t^l(g)$ other than $V^l$ and $\{ 0\}$. The
restriction to $K$ of $t^l$ splits multiplicity-free;
$$
t^l \vert_K = \bigoplus_{n=-l}^l \psi_{-2n}.
$$
Let $\{ e^l_n\mid n=-l,-l+1,\ldots,l\}$ be an orthonormal basis for
$V^l$ such that $t^l(k)e^l_n=\psi_{-2n}(k) e^l_n$, $k\in K$. The matrix
elements $SU(2)\ni g\mapsto t^l_{n,m}(g)=\langle t^l(g)e^l_m, e^l_n
\rangle$ can be calculated in terms of Jacobi polynomials
$$
P^{(\a,\b)}_n(x)= {{(\a+1)_n}\over{n!}} \, {}_2F_1\left(
{{-n,n+\a+\b+1}\atop{\a+1}}; {{1-x}\over 2}\right)
\tag\eqname{\vgldefJacobipol}
$$
for $\a,\b\in\Zp$. Here we use the standard notation
$$
{}_2F_1\left( {{a,b}\atop{c}}; z\right) = \sum_{k=0}^\infty
{{(a)_k(b)_k}\over{(c)_k k!}} z^k, \qquad (a)_k = \prod_{i=0}^{k-1}
(a+i),
$$
for Gauss's hypergeometric series ${}_2F_1$ and for Pochhammer's
symbol $(a)_k$.

The Schur orthogonality relations for the matrix elements are
$$
\int_{SU(2)} t^l_{n,m}(g) \overline{t^k_{i,j}(g)} \, dg =
\d_{l,k}\d_{n,i}\d_{m,j} (2l+1)^{-1},
$$
where $dg$ denotes the normalised Haar measure on $SU(2)$. For $n=i$,
$m=j$, these relations yield the orthogonality relations
for the Jacobi polynomials
$$
\int_{-1}^1 P^{(\a,\b)}_l(x)P^{(\a,\b)}_k(x) (1-x)^\a (1+x)^\b \, dx = 0,
\qquad l\not= k,
$$
for $\a,\b\in\Zp$.

Since $t^l(g)$ is a unitary matrix for $g\in SU(2)$ we have
$\sum_k t^l_{m,k}(g) \overline{t^l_{n,k}(g)} = \d_{n,m}$. This
implies the finite discrete orthogonality relation
$$
\sum_{x=0}^N K_m(x;p,N)K_n(x;p,N) {N\choose x} p^x (1-p)^{N-x} =
\d_{n,m} (1-p)^np^{-n} {N\choose n}^{-1}
$$
for the Krawtchouk polynomials
$$
K_n(x;p,N) = {}_2F_1\left( {{-n,-x}\atop{-N}};{1\over p}\right),
\tag\eqname{\vgldefKrawtchouk}
$$
where $N\in\Zp$, $n\in\{ 0,1,\ldots,N\}$, and $0<p<1$.

The commutative $\ast$-algebra of matrix elements of finite dimensional
unitary representations of $SU(2)$ is a dense $\ast$-subalgebra of the
$C^\ast$-algebra of continuous functions on $SU(2)$. Moreover, this
algebra equals the $\ast$-algebra $P(SU(2))$ of polynomials in the
coordinate functions $\a$, $\b$, $\g$, $\d$, cf.
\thetag{\vgldefinitieSUtwee}, subject to the relations
$\a \d-\b \g=1$, $\a^\ast=\d$, $\b^\ast=\g$, where the $\ast$-operator is
complex conjugation. This algebra has much more structure, which stems
from the group structure of $SU(2)$. So multiplication, unit and
inverse give rise to the mappings
$$
\gather
\bigl(\Delta p\bigr) (g,h) =p(gh), \qquad \Delta \colon P(SU(2)) \to
P(SU(2)\times SU(2)) \simeq P(SU(2)) \otimes P(SU(2)), \\
\varepsilon(p) = p(\pmatrix 1&0\\ 0&1\endpmatrix ),\qquad
\varepsilon \colon P(SU(2)) \to \C, \\
\bigl( Sp\bigr) (g) = p(g^{-1}), \qquad S\colon P(SU(2)) \to P(SU(2)).
\endgather
$$
and with these mappings $P(SU(2))$ becomes a Hopf $\ast$-algebra, cf. \S
4.1 for the definition. %forward reference
(The reader is invited to check the axioms of a
Hopf $\ast$-algebra in this particular case.)

If we forget about the $\ast$-structure on the Hopf $\ast$-algebra
$P(SU(2))$ we can consider the Hopf algebra $P(SU(2))$ as the algebra of
polynomials in coordinate functions on the complex group $SL(2,\C)$. The
$\ast$-operator determines a real form of $SL(2,\C)$ and in this case we
can recover the real group $SU(2)$ from $P(SU(2))$ by
$$
SU(2) = \Bigl\lbrace g\in SL(2,\C) \mid p^\ast(g) = \overline{p(g)},\
\forall p\in P(SU(2)) \Bigr\rbrace .
$$

For the matrix elements
$t^l_{n,m}$ of the irreducible unitary representation $t^l$ we get from
$t^l(gh)=t^l(g)t^l(h)$, $t^l(\pmatrix 1&0\\ 0&1\endpmatrix ) = id$,
$t^l(g^{-1})=t^l(g)^\ast$ the identities
$$
\Delta(t^l_{n,m})=\sum_{k=-l}^l t^l_{n,k}\otimes t^l_{k,m}, \quad
\varepsilon (t^l_{n,m})=\delta_{n,m}, \quad
S(t^l_{n,m})=(t^l_{m,n})^\ast ,
$$
cf. (4.18).

The Lie algebra ${\frak{su}}(2) = \{ X\in{\frak{gl}}(2,\C)\mid trX=0,
X^\ast+X=0\}$ is equal to $ihH+bB-\bar b C$ for $h\in\R$, $b\in\C$. Here
$$
H = \pmatrix 1&0\\ 0&-1\endpmatrix, \qquad
B = \pmatrix 0&1\\ 0&0 \endpmatrix, \qquad
C = \pmatrix 0&0\\ 1&0 \endpmatrix,
$$
satisfying $[H,B]=2B$, $[H,C]=-2C$, $[B,C]=H$, span the complexification
${\frak{sl}}(2,\C)$ of ${\frak{su}}(2)$. The representations $t^l$ of
$SU(2)$ can be differentiated yielding representations of the Lie
algebra ${\frak{su}}(2)$ which in turn can be extended to the whole
universal enveloping algebra of ${\frak{sl}}(2,\C)$.
The universal enveloping
algebra ${\frak{U}}({\frak{sl}}(2,\C))$ is another example of a Hopf
algebra if we define $\Delta(1) = 1\otimes 1$, $\varepsilon (1) = 1$,
$S(1)=1$, and
$$
\Delta(X)=1\otimes X+X\otimes 1, \quad \varepsilon(X)=0,\quad
S(X)=-X,\qquad X\in{\frak{sl}}(2,\C),
\tag\eqname{\vglLiealgelminU}
$$
and by extending $\Delta$, $\varepsilon$ as homomorphisms and $S$ as an
antihomomorphism. We can introduce a $\ast$-operator making
${\frak{U}}({\frak{sl}}(2,\C))$ into a Hopf $\ast$-algebra and such that
the $\ast$-operator can be used to recover the real Lie algebra
${\frak{su}}(2)$ from it.
Here the $\ast$-operator is defined by $H^\ast=H$,
$B^\ast=C$, $C^\ast=B$, satisfying $[X^\ast,Y^\ast]=[Y,X]^\ast$. The
$\ast$-operator is chosen in such a way that the Lie algebra
${\frak{su}}(2)$ consists of the $-1$-eigenspace of the $\ast$-operator
when restricted to the Lie algebra ${\frak{sl}}(2,\C)$.

For $X\in{\frak{sl}}(2,\C)$ and $p\in P(SU(2))$ we have the pairing
$$
\langle X, p\rangle = {{d}\over{dt}}\big\vert_{t=0} p(\exp tX),
\tag\eqname{\vgldefklassiekparing}
$$
by using the exponential function from the Lie algebra ${\frak{sl}}(2,\C)$
to the Lie group $SL(2,\C)$. This pairing extends to a pairing of
${\frak{U}}({\frak{sl}}(2,\C))$ and $P(SU(2))$ and this pairing makes the
Hopf $\ast$-algebras ${\frak{U}}({\frak{sl}}(2,\C))$ and $P(SU(2))$ in
duality as Hopf $\ast$-algebras, cf. \S 4.4. %forward reference
The relation
$$
t^l_{n,m}(k_\theta g k_\psi) = \langle t^l(g) t^l(k_\psi) e^l_m ,
t^l(k_{-\theta})e^l_n\rangle = e^{-2im\psi-2in\theta} t^l_{n,m}(g)
$$
for $k_\theta=diag(e^{i\theta},e^{-i\theta})\in K$, $g\in SU(2)$, can be
rephrased as
$$
\langle id\otimes iH,\Delta(t^l_{n,m})\rangle = -2im t^l_{n,m},\qquad
\langle iH\otimes id,\Delta(t^l_{n,m})\rangle = -2in t^l_{n,m},
$$
where we use the duality \thetag{\vgldefklassiekparing} in the second,
respectively first factor of the tensor product. If, instead of $K$, we
start with two, possibly different, one-parameter subgroups of $SU(2)$,
such as $\exp tX$ and $\exp tY$ for $X,Y\in {\frak{su}}(2)$, then we may
try to find matrix elements $s^l_{n,m}$ such that
$$
\langle id\otimes X,\Delta(s^l_{n,m})\rangle = -2im s^l_{n,m},\qquad
\langle Y\otimes id,\Delta(s^l_{n,m})\rangle = -2in s^l_{n,m},
$$
i.e. matrix elements which are relatively right (left) invariant with
respect to the one-parameter subgroup $\exp tX$ ($\exp tY$). However,
this does not lead to a different interpretations of Jacobi polynomials,
since it is equivalent to an affine transformation of the argument of
the Jacobi polynomials. The underlying reason for this is that
all one-parameter subgroups in $SU(2)$ are conjugated.

In the quantum $SU(2)$ group case this approach leads to interpretations
of different $q$-analogues of the Jacobi polynomials

%&&&&&&&&&&&&&&&&&&&&&&&&&&&&&&&&&&&&
%%N E W   S E C T I O N%%%%%%%%%%%%**
%%%%%%%%%%%%%%%%%%%%%%%%%%%%%%%%%%%%%

\head\newsection . Preliminaries on basic hypergeometric orthogonal
polynomials
\endhead

In this section we collect various results on basic hypergeometric
orthogonal polynomials that appear in this paper. Per subsection
references to the literature are given. The reader acquainted with
$q$-hypergeometric series is advised to browse through this section. The
notation follows the book \cite{\GaspR} by Gasper and Rahman. See
Andrews \cite{\Andr} for more
information concerning applications of $q$-series.

%&&&&&  NEW SUBSECTION &&&&&&&&&&&&&&&&&&&&&&&&&&&&&&&&&&&&&&&&&&&&&&&&&
\subhead\newsubsection{\ssectqhypergeometricseries}
Basic hypergeometric series
\endsubhead (\cite{\GaspR , Ch.~1})
We consider $q\in (0,1)$ fixed.
The $q$-shifted factorials are defined by
$(a;q)_k = \prod_{i=0}^{k-1}(1-aq^i)$ for $a\in\C$, $k\in\Zp$. The empty
product corresponding to $k=0$ equals $1$ by definition. Since the
absolute value of $q$ is less than $1$, the limit $k\to\infty$ of any
$q$-shifted factorial $(a;q)_k$ exists. This infinite product is denoted
by $(a;q)_\infty$. The following abbreviation for $q$-shifted factorials
$$
(a_1,a_2,\ldots, a_r;q)_k = \prod_{p=1}^r (a_p;q)_k, \qquad
k\in\Zp\cup\{\infty\}
$$
is standard. Occasionally we use $\left[ {n\atop k}\right]_q =
(q^n;q^{-1})_k/(q;q)_k$, $n,k\in\Zp$, $k\leq n$, as a $q$-analogue of
the binomial coefficient. Note that $(x;q)_n$ is a polynomial of degree
$n$ in $x$ and that $(ae^{i\theta},a^{-i\theta};q)_n = \prod_{i=0}^n
(1-2aq^i\cos\theta+a^2q^{2i})$ is a polynomial of degree $n$ in
$\cos\theta$ for $a\not= 0$.

Define basic hypergeometric series, or $q$-hypergeometric
series, with upper parameters $a_1,\ldots,a_{r+1}$ and lower parameters
$b_1,\ldots,b_r$ ($r\in\Zp$) of argument $z$ and base $q$ by
$$
{}_{r+1}\vp_r \left(
{{a_1,\ldots,a_{r+1}}\atop{b_1,\ldots,b_r}};q,z\right) =
\sum_{k=0}^\infty {{(a_1,\ldots,a_{r+1};q)_k}\over{
(b_1,\ldots,b_r,q;q)_k}} z^k .
\tag\eqname{\vgldefqhypergeoseries}
$$
For generic values of the parameters the radius of convergence of this
series is $1$. In this paper the $q$-hypergeometric series
\thetag{\vgldefqhypergeoseries} are terminating series. This is the case
if the upper parameter $a_1=q^{-n}$, $n\in\Zp$, since $(q^{-n};q)_k=0$
for $k>n$. If one of the lower parameters is of the form $q^{-N}$ for
$N\in\Zp$, then we require that there is an upper parameter of the form
$q^{-n}$ for $n\in\{ 0,\ldots,N\}$ in order to have a well-defined
$q$-hypergeometric series. Here we follow the convention that a
${}_{r+1}\vp_r$-series with $q^{-N}$ as both an upper and lower parameter
is considered as a terminating series of degree $N$.

%&&&&&  NEW SUBSECTION &&&&&&&&&&&&&&&&&&&&&&&&&&&&&&&&&&&&&&&&&&&&&&&&&
\subhead\newsubsection{\ssectAWpols}
Askey-Wilson polynomials
\endsubhead (\cite{\AskeW}, \cite{\GaspR , \S 7.5})
The Askey-Wilson polynomials are defined by
$$
p_n(\cos\theta;a,b,c,d\mid q) = a^{-n} (ab,ac,ad;q)_n
\, {}_4\vp_3 \left( {{q^{-n}, abcdq^{n-1}, ae^{i\theta},a^{-i\theta}}
\atop{ab,\ ac,\ ad}};q,q\right),
$$
cf. Askey and Wilson \cite{\AskeW, (1.15)}. The Askey-Wilson polynomial is
symmetric in the parameters $a$, $b$, $c$ and $d$, cf. \cite{\AskeW, p.6}.
To stress the fact that Askey-Wilson polynomials generalise Jacobi
polynomials (see the title of Askey-Wilson memoir \cite{\AskeW}),
we use the notation, cf. Noumi and Mimachi \cite{\NoumM , (4.1)},
$$
p_n^{(\a,\b)}(x;s,t\mid q) = p_n(x;q^{1/2}t/s,q^{1/2+\a}s/t,
-q^{1/2}/(st), -stq^{1/2+\b}\mid q),
\tag\eqname{\vgldefAWalsqJacobi}
$$
which generalises Rahman's notation ($s=t=1$) for continuous
$q$-Jacobi polynomials, cf. \cite{\AskeW , (4.17)}.
%?? Nog te checken: artikel Rahman
With this notation we have, cf. \thetag{\vgldefJacobipol},
$$
\align
\lim_{q\uparrow 1} {{p_n^{(\a,\b)}(x;s,t\mid q)}\over{(q;q)_n}} = &
(s+s^{-1})^n(t+t^{-1})^n {{(\a+1)_n}\over{n!}}\\
&\times {}_2F_1\left( {{-n,n+\a+\b+1}\atop{\a+1}};
{{s/t-2\cos\theta+t/s}\over{(s+s^{-1})(t+t^{-1})}}\right).
\endalign
$$

%&&&&&  NEW SUBSECTION &&&&&&&&&&&&&&&&&&&&&&&&&&&&&&&&&&&&&&&&&&&&&&&&&
\subhead\newsubsection{\ssectAWorthog}
Orthogonality relations for Askey-Wilson polynomials
\endsubhead (\cite{\AskeW}, \cite{\GaspR , \S 7.5})
The orthogonality relations for the Askey-Wilson polynomials depend on the
values of the parameters $a$, $b$, $c$ and $d$. First we introduce some
notation;
$$
\aligned
w\bigl( {1\over 2}(z+z^{-1})\bigr) &=
{{(z^2,z^{-2};q)_\infty}\over{(az,a/z,bz,b/z,cz,c/z,dz,d/z;q)_\infty}}, \\
h_n &= {{(1-q^{n-1}abcd) (q,ab,ac,ad,bc,bd,cd;q)_n}\over
{(1-q^{2n-1}abcd) (abcd;q)_n}} h_0, \\
h_0 &= {{(abcd;q)_\infty}\over{(q,ab,ac,ad,bc,bd,cd;q)_\infty}},
\endaligned
$$
where we suppressed the dependence on $a$, $b$, $c$ and $d$ in the notation
for $w$ and $h$. For $z=e^{i\theta}$ we use $w(\cos\theta)$.

\proclaim{Proposition \theoremname{\propAWorthogonality}}
Let $a$, $b$, $c$ and $d$ be real and let all the pairwise products of
$a$, $b$, $c$ and $d$ be less than $1$. Then the Askey-Wilson polynomials
$p_n(x)= p_n(x;a,b,c,d\mid q)$ satisfy the orthogonality relations
$$
{1\over {2\pi}} \int_0^\pi p_n(\cos\theta)p_m(\cos\theta)w(\cos\theta)d\theta
+\sum_k p_n(x_k)p_m(x_k)w_k = \delta_{n,m} h_n.
$$
The points $x_k$ are of the form ${1\over 2}(eq^k+e^{-1}q^{-k})$ for $e$ any
of the parameters $a$, $b$, $c$ or $d$ with absolute value greater than $1$;
the sum is over $k\in\Zp$ such that $\vert eq^k\vert >1$ and $w_k$ is the
residue of $z\mapsto w\bigl({1\over 2}(z+z^{-1})\bigr)$ at
$z=eq^k$ minus the residue at $z=e^{-1}q^{-k}$.
\endproclaim

The orthogonality relations remain valid for complex parameters $a$, $b$, $c$
and $d$, if they occur in conjugate pairs. If all parameters have
absolute value less than $1$, the Askey-Wilson orthogonality measure is
absolutely continuous.

We use the notation $dm(x)=dm(x;a,b,c,d\mid q)$ for the normalised
orthogonality measure. So for any polynomial $p$
$$
\int_\R p(x) \, dm(x) = {1\over{h_0}} \Bigl( {1\over{2\pi}}
\int_{-1}^1 p(x)w(x) {{dx}\over{\sqrt{1-x^2}}} + \sum_k p(x_k)w_k
\Bigr).
\tag\eqname{\vglnormalisedAWmeasure}
$$

%&&&&&  NEW SUBSECTION &&&&&&&&&&&&&&&&&&&&&&&&&&&&&&&&&&&&&&&&&&&&&&&&&
\subhead\newsubsection{\ssectcontqultraspherical}
Continuous $q$-ultraspherical polynomials \endsubhead
(\cite{\AskeI}, \cite{\AskeW, \S 4}, \cite{\GaspR , \S 7.4}) The
continuous $q$-ultraspherical polynomials are special cases of the
Askey-Wilson polynomials, but were already introduced by Rogers in 1894.
They are defined by
$$
C_n(\cos\theta ; \b \mid q) = \sum_{k=0}^n {{(\b;q)_k
(\b;q)_{n-k}}\over{(q;q)_k (q;q)_{n-k}}} e^{i(n-2k)\theta}.
\tag\eqname{\vgldefcontqultraspherpols}
$$
The orthogonality relation for the continuous $q$-ultraspherical
polynomials is absolutely continuous for $-1<\b<1$. Up to a
normalisation factor the continuous $q$-ultraspherical polynomials
correspond to Askey-Wilson polynomials with parameters $a=\sqrt{\b}$,
$b=\sqrt{q\b}$, $c=-\sqrt{q\b}$, $d=-\sqrt{\b}$, cf. \cite{\AskeW ,
\S 4}. There is also a relation between continuous $q$-ultraspherical
polynomials of base $q^2$ and Askey-Wilson polynomials of base $q$ given
by
$$
C_n(\cos\theta;q^{1+2\a}\mid q^2) =
{{(q^{1+2\a};q)_n}\over{(q^{2+2\a},q^2;q^2)_n}}
p_n^{(\a,\a)}(\cos\theta;1,1\mid q),
\tag\eqname{\vglquadraticcqusAWp}
$$
cf. \cite{\AskeW , (4.20), (4.2)}, \cite{\GaspR , (7.5.34)}.

In case $\b=0$ we call $H_n(\cos\theta\mid q) = (q;q)_n C_n(\cos\theta;
0 \mid q)$ the continuous $q$-Hermite polynomials, cf. \cite{\AskeI , \S
6}. The three-term recurrence relation for the continuous $q$-Hermite
polynomials is
$$
2x H_n(x\mid q) = H_{n+1}(x\mid q) + (1-q^n) H_{n-1}(x\mid q).
\tag\eqname{\vgldrietermcontqHermite}
$$
The corresponding orthogonality relation can be obtained from
proposition \propAWorthogonality\ and we get
$$
\int_0^\pi H_n(\cos\theta\mid q) H_m(\cos\theta\mid q)
(e^{2i\theta},e^{-2i\theta};q)_\infty \, d\theta = \d_{n,m} {{2\pi
(q;q)_n}\over{(q;q)_\infty}}.
\tag\eqname{\vglorthrelcontqHermite}
$$

%&&&&&  NEW SUBSECTION &&&&&&&&&&&&&&&&&&&&&&&&&&&&&&&&&&&&&&&&&&&&&&&&&
\subhead\newsubsection{\ssectlittleqJacobi}
Little $q$-Jacobi polynomials
\endsubhead (\cite{\AndrA}, \cite{\GaspR , \S 7.3})
The little $q$-Jacobi polynomials
$$
p_n(x;a,b;q) = {}_2\vp_1 \left( {{q^{-n},abq^{n-1}}\atop{aq}};q,qx
\right)
\tag\eqname{\vgldeflittleqJacobi}
$$
are $q$-analogues of the Jacobi polynomials. The orthogonality relations
are
$$
\sum_{x=0}^\infty {{(bq;q)_x}\over{(q;q)_x}} (aq)^x p_n(q^x;a,b;q)
p_m(q^x;a,b;q) = 0, \qquad n\not= m,
\tag\eqname{\vglorthorellittleqJac}
$$
cf. Andrews and Askey \cite{\AndrA , thm.~9}. The little $q$-Jacobi
polynomials can be obtained as a limit case of the Askey-Wilson
polynomials of \ssectAWpols , cf. Koornwinder \cite{\Koortwee ,
prop.~6.3} and the orthogonality relations, cf. proposition
\propAWorthogonality , go over into the orthogonality relations
\thetag{\vglorthorellittleqJac} for the little $q$-Jacobi polynomials.

%&&&&&  NEW SUBSECTION &&&&&&&&&&&&&&&&&&&&&&&&&&&&&&&&&&&&&&&&&&&&&&&&&
\subhead\newsubsection{\ssectqHahn}
$q$-Hahn polynomials
\endsubhead (\cite{\AskeWtwee}, \cite{\GaspR , \S 7.2})
The $q$-Hahn polynomials are defined by, cf. \cite{\GaspR , (7.2.22)},
$$
Q_n(x) = Q_n(x;a,b,N;q) = {}_3\vp_2 \left( {{q^{-n},abq^{n+1},x}\atop
{aq,q^{-N}}};q,q\right)
\tag\eqname{\vgldefqHahnpol}
$$
for $N\in\Zp$ and $n\in\{ 0,1,\ldots,N\}$. These polynomials are
orthogonal with respect to a finite discrete measure;
$$
\sum_{x=0}^N Q_n(q^{-x})Q_m(q^{-x}) {{(aq;q)_x (bq;q)_{N-x}}\over
{(q;q)_x (q;q)_{N-x}}} (aq)^{-x} = 0, \qquad n\not= m.
\tag\eqname{\vglorthorelqHahn}
$$
The $q$-Hahn polynomials are special cases of the more general
orthogonal polynomials, the so-called $q$-Racah polynomials introduced
by Askey and Wilson \cite{\AskeWtwee}.

%&&&&&  NEW SUBSECTION &&&&&&&&&&&&&&&&&&&&&&&&&&&&&&&&&&&&&&&&&&&&&&&&&
\subhead\newsubsection{\ssectqKrawtchouk}
(Dual) $q$-Krawtchouk polynomials
\endsubhead (\cite{\Stan}, \cite{\AskeWtwee})
The $q$-Krawtchouk are defined by, cf. \cite{\Stan},
$$
K_n(x;q^\s,N;q) = {}_3\vp_2 \left( {{q^{-n},x,-q^{n-N-\s}}\atop
{q^{-N},0}};q,q\right)
$$
and these polynomials are orthogonal with respect to a finite discrete
measure;
$$
{{q^{N+\s}}\over{(-q^\s;q)_{N+1}}} \sum_{x=0}^N
\bigl( K_nK_m\bigr) (q^{-x};q^\s,N;q) w_x(q^\s,N)
= \d_{n,m} \bigl( h_n(q^\s,N)\bigr)^{-1}
\tag\eqname{\vglorthorelqKrawtchouk}
$$
with
$$
\gather
w_x(q^\s,N) = (-q^{N+\s})^x {{(q^{-N};q)_x}\over{(q;q)_x}}, \\
h_n(q^\s,N) = {{(1+q^{2n-N-\s})}\over{(-q^{n-2N-\s})^n}}
{{(-q^{-N-\s},q^{-N};q)_n}\over{(q,-q^{1-\s};q)_n}}.
\endgather
$$
The dual $q$-Krawtchouk polynomials are obtained from the $q$-Krawtchouk
polynomials by interchanging the roles of the degree $n$ and argument
$q^{-x}$,
$$
R_n(q^{-x}-q^{x-N-\s};q^\s,N;q) =
{}_3\vp_2 \left( {{q^{-n},q^{-x},-q^{x-N-\s}}\atop
{q^{-N},0}};q,q\right) = K_x(q^{-n};q^\s,N;q).
\tag\eqname{\vgldefdualqKrawtchoukpol}
$$
The orthogonality relations for the dual $q$-Krawtchouk polynomials
follow from \thetag{\vglorthorelqKrawtchouk};
$$
{{q^{N+\s}}\over{(-q^\s;q)_{N+1}}} \sum_{x=0}^N
\bigl( R_nR_m\bigr) (q^{-x}-q^{x-N-\s};q^\s,N;q)
h_x(q^\s,N) = \d_{n,m} \bigl( w_n(q^\s,N)\bigr)^{-1}
\tag\eqname{\vglorthoreldualqKrawtchouk}
$$
with $w_n(q^\s,N)$ and $h_x(q^\s,N)$ defined as above. Again both the
$q$-Krawtchouk and dual $q$-Krawtchouk polynomials are special cases of
the $q$-Racah polynomials \cite{\AskeWtwee}.

%&&&&&&&&&&&&&&&&&&&&&&&&&&&&&&&&&&&&
%%N E W   S E C T I O N%%%%%%%%%%%%**
%%%%%%%%%%%%%%%%%%%%%%%%%%%%%%%%%%%%%

\head\newsection . Preliminaries on the quantum $SU(2)$
group\endhead

The necessary tools and definitions concerning the quantum $SU(2)$ group
are given in this section. There is no notational standard on this
subject. References to the literature are given per subsection.

%&&&&&  NEW SUBSECTION &&&&&&&&&&&&&&&&&&&&&&&&&&&&&&&&&&&&&&&&&&&&&&&&&
\subhead\newsubsection{\ssectHopfstaralgs}
Hopf $\ast$-algebras\endsubhead (\cite{\Abe}, \cite{\Swee},
\cite{\VDael})
A complex associative algebra ${\Cal A}$ with unit $1$ and
multiplication $m\colon {\Cal A}\otimes{\Cal A}\to {\Cal A}$, $m\colon
a\otimes b\mapsto ab$ and unit $e\colon \C\to{\Cal A}$, $e\colon z\mapsto
z1$ is a bialgebra if there exist algebra homomorphisms $\Delta \colon
{\Cal A}\to{\Cal A}\otimes{\Cal A}$, the comultiplication, and
$\varepsilon \colon {\Cal A}\to\C$, the counit, satisfying the
coassociativity axiom $(\Delta\otimes id)\circ\Delta =
(id\otimes \Delta)\circ\Delta$ and the counit axiom
$(\varepsilon\otimes id)\circ\Delta = id =
(id\otimes \varepsilon)\circ\Delta$.
A bialgebra ${\Cal A}$ is a Hopf algebra if there exists a linear map
$S\colon {\Cal A}\to{\Cal A}$, the antipode, satisfying the antipode
axiom $m\circ(S\otimes id)\circ\Delta = e\circ\varepsilon
= m\circ(id\otimes S)\circ\Delta$.
The antipode $S$ in the Hopf algebra ${\Cal A}$ is unique. Moreover, it
satisfies $S(1)=1$, $\varepsilon\circ S=\varepsilon$, $S(ab)=S(b)S(a)$
and $\Delta\circ S = \omega\circ (S\otimes S)\circ \Delta$, where
$\omega$ denotes the flip automorphism of ${\Cal A}\otimes {\Cal A}$,
$\omega(a\otimes b) =b\otimes a$.

A $\ast$-operator on an algebra ${\Cal A}$ is an antilinear,
antimultiplicative involution, i.e. $(\lambda a+ \mu b)^\ast =
\bar\lambda a^\ast+\bar\mu b^\ast$, $(ab)^\ast=b^\ast a^\ast$,
$(a^\ast)^\ast = a$ for $\lambda,\mu\in\C$, $a,b\in{\Cal A}$.
A Hopf algebra ${\Cal A}$ with a $\ast$-operator is a Hopf
$\ast$-algebra if $(\ast\otimes\ast)\circ\Delta = \Delta\circ\ast$,
$\varepsilon(a^\ast)=\overline{\varepsilon(a)}$ for all $a\in{\Cal
A}$, and then it is possible to prove that
$S\circ\ast\circ S\circ\ast =id$.

%&&&&&  NEW SUBSECTION &&&&&&&&&&&&&&&&&&&&&&&&&&&&&&&&&&&&&&&&&&&&&&&&&
\subhead\newsubsection{\ssectquanfunalg}
Quantised algebra of polynomials on $SU(2)$
\endsubhead (\cite{\Woro}, \cite{\Worotwee}, \cite{\KoorOPTA})
Fix $q\in(0,1)$. $\A$ is the complex unital associative algebra generated
by $\a$, $\b$, $\g$, $\d$ subject to the relations
$$
\eqalign{
&\a\b =q\b\a ,\quad \a\g = q\g\a ,\quad \b\d = q\d\b ,\quad \g\d = q\d\g ,\cr
&\b\g =\g\b ,\quad \a\d -q\b\g = \d\a - q^{-1}\b\g =1 .\cr}
\tag\eqname{\vglcommrelAq}
$$
The algebra $\A$ is an example of a Hopf $\ast$-algebra.
The comultiplication $\Delta$, the counit $\varepsilon$, the antipode $S$
and the $\ast$-operator are given on the generators by
$$
\gather
\eqalign{
&\Delta(\a)=\a\otimes\a + \b\otimes\g ,\quad
\Delta(\b )=\a\otimes\b + \b\otimes\d ,\cr
&\Delta(\g )=\g\otimes\a + \d\otimes\g ,\quad
\Delta(\d )=\g\otimes\b + \d\otimes\d , \cr}
\tag\eqname{\vgldefDeltaonAq}\\
\varepsilon\pmatrix \a &\b\\ \g&\d\endpmatrix = \pmatrix
1&0\\0&1\endpmatrix,\qquad
S\pmatrix \a &\b\\ \g&\d\endpmatrix = \pmatrix \d&-q^{-1}\b \\
-q\g&\a\endpmatrix,
\tag\eqname{\vglSenepsilonopAq} \\
\a^\ast = \d , \quad \b^\ast = -q\g ,\quad \g^\ast = -q^{-1}\b ,\quad
\d^\ast = a .
\tag\eqname{\vglsteropAq}
\endgather
$$
It is possible to identify ${\Cal A}_1(SU(2))$ with the polynomial
algebra $P(SU(2))$, cf. \S 2, on the Lie group $SU(2)$.

%&&&&&  NEW SUBSECTION &&&&&&&&&&&&&&&&&&&&&&&&&&&&&&&&&&&&&&&&&&&&&&&&&
\subhead\newsubsection{\ssectquantunivenvalg}
Quantised universal enveloping algebra \endsubhead (\cite{\Jimb},
\cite{\KoorOPTA})
The quantised universal enveloping algebra $\U$ is the complex unital
associative algebra generated by $A$, $B$, $C$, $D$ subject to the relations
$$
AD=1=DA, \quad AB=qBA,\quad AC=q^{-1}CA,\quad
BC-CB = {{A^2-D^2}\over{q-q^{-1}}}.
\tag\eqname{\vglcommrelUq}
$$
The algebra $\U$ is also a Hopf $\ast$-algebra. The comultiplication, counit,
antipode and $\ast$-operator are defined on the generators by
$$
\gather
\eqalign{
\Delta (A)=A\otimes A,\quad \Delta(B)=A\otimes B+B\otimes D,\cr
\Delta(C) = A\otimes C+C\otimes D, \quad \Delta(D)=D\otimes D,\cr}
\tag\eqname{\vgldefDeltaonUq} \\
\eqalign{
&\qquad\qquad \varepsilon(A)=\varepsilon(D)=1,\qquad
\varepsilon(C)=\varepsilon(B)=0,\cr
&S(A)=D,\quad  S(B)=-q^{-1}B,\quad S(C)=-qC,\quad S(D)=A,\cr}
\tag\eqname{\vgldefSenepsilononUq} \\
A^\ast=A,\quad B^\ast=C, \quad C^\ast=B,\quad D^\ast=D.
\tag\eqname{\vgldefsteronUq}
\endgather
$$
The element
$$
\Omega = {{q^{-1}A^2+qD^2-2}\over{(q^{-1}-q)^2}} + BC
 = {{qA^2+q^{-1}D^2-2}\over{(q^{-1}-q)^2}} + CB
\tag\eqname{\vgldefCasimir}
$$
is the Casimir element of the quantised universal enveloping algebra
$\U$. The Casimir element is self-adjoint, $\Omega=\Omega^\ast$, and
$\Omega$ belongs to the centre of $\U$.

%&&&&&  NEW SUBSECTION &&&&&&&&&&&&&&&&&&&&&&&&&&&&&&&&&&&&&&&&&&&&&&&&&
\subhead\newsubsection{\ssectdualHopfalg}
Duality for Hopf $\ast$-algebras \endsubhead (\cite{\VDael},
\cite{\MasuMNNSUeen})
The Hopf $\ast$-algebras are in duality as Hopf $\ast$-algebras if there
exists a doubly non-degenerate pairing $\langle\cdot,\cdot\rangle \colon
\U\times\A\to\C$, i.e. $\langle X,\xi\rangle=0$ for all $X\in\U$,
respectively for all $\xi\in\A$, implies $\xi=0$, respectively $X=0$,
such that
$$
\gather
\langle XY,\xi\rangle = \langle X\otimes Y, \Delta(\xi)\rangle,\quad
\langle X,\xi\eta\rangle = \langle \Delta(X), \xi\otimes\eta\rangle ,
\tag\eqname{\vgldualcomult}\\
\langle 1,\xi\rangle = \varepsilon(\xi),\quad \langle X, 1\rangle =
\varepsilon(X),
\tag\eqname{\vgldualunit}  \\
\langle S(X),\xi\rangle = \langle X,S(\xi)\rangle, \quad
\langle X^\ast, \xi\rangle = \overline{\langle X, S(\xi)^\ast\rangle} .
\tag\eqname{\vgldualSster}
\endgather
$$

The duality can be used to define a left and right action of $\U$ on $\A$.
For $X\in\U$ and $\xi\in\A$ we define elements $X.\xi$ and $\xi.X$ of
$\A$ by
$$
X.\xi = (id\otimes X)\circ \Delta (\xi),\qquad
\xi.X = (X\otimes id)\circ \Delta (\xi),
\tag\eqname{\vgldefactionUonA}
$$
where the pairing between $\A$ and $\U$ is used in the second, respectively
the first, part of the tensor product. Using the duality we can rewrite
\thetag{\vgldefactionUonA} as
$$
\langle Y,X.\xi\rangle = \langle YX,\xi\rangle,\qquad
\langle Y,\xi.X\rangle = \langle XY,\xi\rangle .
\tag\eqname{\vglweakactionUonA}
$$
If $\Delta(X)=\sum_{(X)} X_{(1)}\otimes X_{(2)}$, then
\thetag{\vgldualcomult} implies
$$
X.(\xi\eta) = \sum_{(X)} \bigl( X_{(1)}.\xi\bigr)
\bigl(X_{(2)}.\eta\bigr), \qquad
(\xi\eta).X = \sum_{(X)} \bigl( \xi.X_{(1)}\bigr)
\bigl(\eta.X_{(2)}\bigr).
\tag\eqname{\vglactionUonprodA}
$$
Using \thetag{\vgldefactionUonA} and \thetag{\vgldualSster} we prove
$$
X.\xi^\ast = \sum_{(\xi)} \xi_{(1)}^\ast \langle X,\xi^\ast_{(2)}\rangle
=\Bigl(\sum_{(\xi)} \xi_{(1)} \langle S(X)^\ast,\xi_{(2)}\rangle
\Bigr)^\ast = \bigl( S(X)^\ast.\xi\bigr)^\ast.
\tag\eqname{\vglactionUonAandstar}
$$

The explicit duality between $\A$ and $\U$ is given by
$$
\aligned
\Bigl\langle A, \pmatrix \a&\b\\ \g&\d\endpmatrix \Bigr\rangle =
\pmatrix q^{1/2}&0 \\ 0&q^{-1/2}\endpmatrix , \quad
\Bigl\langle B, \pmatrix \a&\b\\ \g&\d\endpmatrix \Bigr\rangle =
\pmatrix 0&1 \\ 0&0 \endpmatrix , \\
\Bigl\langle C, \pmatrix \a&\b\\ \g&\d\endpmatrix \Bigr\rangle =
\pmatrix 0&0 \\ 1&0 \endpmatrix , \quad
\Bigl\langle D, \pmatrix \a&\b\\ \g&\d\endpmatrix \Bigr\rangle =
\pmatrix q^{-1/2}&0 \\ 0&q^{1/2}\endpmatrix .
\endaligned
\tag\eqname{\vglexpldualityongenAU}
$$
The duality between the Hopf $\ast$-algebras $\A$ and $\U$ is completely
determined by \thetag{\vglexpldualityongenAU} and
\thetag{\vgldualcomult}.

%&&&&&  NEW SUBSECTION &&&&&&&&&&&&&&&&&&&&&&&&&&&&&&&&&&&&&&&&&&&&&&&&&
\subhead\newsubsection{\ssectreprthey}
Representation theory \endsubhead (\cite{\Woro}, \cite{\VaksS},
\cite{\MasuMNNU}, \cite{\Koordrie}, \cite{\KoorOPTA})
A square matrix $t=(t_{n,m})$ of elements of the Hopf $\ast$-algebra $\A$
is called a unitary matrix corepresentation if
$$
\Delta(t_{n,m})=\sum_k t_{n,k}\otimes t_{k,m}, \quad
\varepsilon (t_{n,m})=\delta_{n,m}, \quad
S(t_{n,m})=t_{m,n}^\ast .
\tag\eqname{\vgldefunitarycorep}
$$
Dropping the prefix unitary means dropping the last requirement of
\thetag{\vgldefunitarycorep}.
Using the duality a unitary matrix corepresentation of $\A$ gives a
$\ast$-representation of $\U$ by
$\bigl( t(X)\bigr)_{n,m} = \langle X, t_{n,m} \rangle$, for $X\in\U$.
There is precisely one irreducible unitary corepresentation of $\A$ in each
finite dimension (up to equivalence). The corresponding $\ast$-representation
of $\U$ can be realised on a $(2l+1)$-dimensional ($l\in\hZp$)
vector space with orthonormal basis $\{ e^l_n\}$, $n=-l,-l+1,\ldots,l$
as the matrix representation $t^l=(t^l_{n,m})_{n,m=-l,\ldots,l}$. The
action of the generators is given by
$$
\eqalign{
&t^l(A)e^l_n=q^{-n}e^l_n,\qquad t^l(D)e^l_n=q^ne^l_n, \cr
&t^l(B)e^l_n= {{\sqrt{(q^{-l+n-1}-q^{l-n+1})(q^{-l-n}-q^{l+n})}}\over
{q^{-1}-q}} e^l_{n-1} \cr
&t^l(C)e^l_n= {{\sqrt{(q^{-l+n}-q^{l-n})(q^{-l-n-1}-q^{l+n+1})}}\over
{q^{-1}-q}} e^l_{n+1}, \cr}
\tag\eqname{\vgldefreprU}
$$
where $e^l_{l+1}=0=e^l_{-l-1}$. The action of the Casimir follows from
\thetag{\vgldefCasimir} and \thetag{\vgldefreprU}. Explicitly,
$t^l(\Omega) = q^{1-2l}(1-q^{2l+1})^2/(1-q^2)^2$.

The matrix elements $t^l_{n,m}\in\A$ are explicitly known in terms of the
generators $\a$, $\b$, $\g$ and $\d$, Vaksman and Soibelman \cite{\VaksS
, prop.~6.6}, Masuda et al. \cite{\MasuMNNU , thm.~2.8}, Koornwinder
\cite{\Koordrie , thm.~5.3}. These expressions involve little
$q$-Jacobi polynomials;
$$
\aligned
t^l_{n,m} & = c^l_{n,m} \d^{n+m} \g^{n-m}p_{l-k}(-q^{-1}\b\g), \qquad
(n\geq m\geq -n), \\
t^l_{n,m} & = c^l_{m,n} \d^{n+m} \b^{m-n} p_{l-k}(-q^{-1}\b\g), \qquad
(m\geq n\geq -m), \\
t^l_{n,m} & = c^l_{-n,-m} \b^{m-n} \a^{-m-n} p_{l-k}(-q^{2m+2n-1}\b\g),
\qquad (-n\geq m\geq n), \\
t^l_{n,m} & = c^l_{-m,-n} \g^{n-m} \a^{-m-n} p_{l-k}(-q^{2m+2n-1}\b\g),
\qquad (-m\geq n\geq m),
\endaligned
\tag\eqname{\vglmatrixeltalslittleqJacobi}
$$
with $p_{l-k}(x)=p_{l-k}(x;q^{2\vert n-m\vert}, q^{2\vert
n+m\vert};q^2)$, $k=\max(\vert n\vert, \vert m\vert )$, a little
$q$-Jacobi polynomial, cf. \thetag{\vgldeflittleqJacobi}, and
the constant is given by
$$
c^l_{n,m} = \left[ {{l-m}\atop{n-m}}\right]_{q^2}^{1/2}
\left[ {{l+n}\atop{n-m}}\right]_{q^2}^{1/2} q^{-(n-m)(l-n)}.
$$

%&&&&&  NEW SUBSECTION &&&&&&&&&&&&&&&&&&&&&&&&&&&&&&&&&&&&&&&&&&&&&&&&&
\subhead\newsubsection{\ssectreprsA}
Representations of the Hopf $\ast$-algebra $\A$
\endsubhead (\cite{\VaksS}, \cite{\Woro}, \cite{\Worotwee})
%Woronowicz checken??
The irreducible $\ast$-re\-pre\-sen\-ta\-tions
of the Hopf $\ast$-algebra $\A$ have
been completely classified, cf. Vaksman and Soibelman \cite{\VaksS ,
thm.~3.2}. The one-dimensional $\ast$-representations $\p$ is defined by
$$
\p(\a)=e^{i\theta},\quad \p(\b)=0=\p(\g),\quad \p(\d)=e^{-i\theta}.
\tag\eqname{\vgldefonedimreprA}
$$
Then $\p(t^l_{n,m})=\delta_{n,m}e^{-2ni\theta}$. For $\lambda\not= 0$
we also have one-dimensional representations of $\A$ given by
$$
\t_\lambda(\a)=\lambda,\quad \t_\lambda(\b)=0=\tau_\lambda(\g),\quad
\t_\lambda(\d)=\lambda^{-1}.
$$
Note that $\t_\lambda$ is a $\ast$-representation of $\A$ if and only if
$\lambda=e^{i\theta}$, or $\t_\lambda=\p$, for some $\theta\in [0,2\pi)$.
Now $\t_\lambda(t^l_{n,m})=\delta_{n,m}\lambda^{-2n}$. The counit
$\varepsilon$ coincides with the special case $\t_1=\pi_0$.

Infinite dimensional $\ast$-representations of the Hopf $\ast$-algebra
$\A$ act in the Hilbert space $\ell^2(\Zp)$. For an orthonormal
basis $\{ f_n\}_{n\in\Zp}$ the action of the generators is given by
$$
\pi^\infty_\theta(\a)f_n = \sqrt{1-q^{2n}}f_{n-1}, \quad
\pi^\infty_\theta(\g)f_n = e^{i\theta}q^nf_n,
$$
with the convention $f_{-1}=0$. The operators corresponding to
$\b$ and $\d$ follow by \thetag{\vglsteropAq}

%&&&&&  NEW SUBSECTION &&&&&&&&&&&&&&&&&&&&&&&&&&&&&&&&&&&&&&&&&&&&&&&&&
\subhead\newsubsection{\ssectPWthmCGseries}
Peter-Weyl theorem and Clebsch-Gordan series
\endsubhead (\cite{\Woro}, \cite{\Worotwee}, \cite{\KoelK},
\cite{\KoorOPTA})
The irreducible unitary matrix corepresentations $t^l$ of $\A$
exhaust the set of irreducible unitary corepresentations of $\A$ up to
equivalence. The Peter-Weyl theorem for the quantum $SU(2)$ group states that
the matrix elements $t^l_{n,m}$, $n,m=-l,-l+1,\ldots,l$, $l\in\hZp$, form
a basis for the underlying linear space of $\A$. If we denote by $\Al$ the
span of the matrix elements $t^l_{n,m}$, $n,m=-l,-l+1,\ldots,l$, then
we have the decomposition
$$
\A = \bigoplus_{l\in\hZp} \Al.
\tag\eqname{\vglPeterWeylforA}
$$
The tensor product of two irreducible matrix corepresentations
$t^{l_1}\otimes t^{l_2}$ is the matrix corepresentation with matrix elements
defined by $(t^{l_1}\otimes t^{l_2})_{i,n;j,m}=t^{l_1}_{i,j}t^{l_2}_{n,m}$
for $i,j=-l_1,\ldots,l_1$; $n,m=-l_2,\ldots,l_2$. Using the duality this
can be defined equivalently as the $\ast$-representation of $\U$ given
by $X\mapsto (t^{l_1}\otimes t^{l_2})\Delta(X)$. The Clebsch-Gordan series for
the quantum $SU(2)$ group states that
$$
t^{l_1}\otimes t^{l_2} = \sum_{l=\vert l_1-l_2\vert}^{l_1+l_2}
\hskip-0.3truecm{}^\oplus\ \, t^l,
$$
or in terms of the decomposition \thetag{\vglPeterWeylforA},
$$
{\Cal A}^{l_1}(SU(2))\cdot {\Cal A}^{l_2}(SU(2)) =
\sum_{l=\vert l_1-l_2\vert}^{l_1+l_2} {\Cal A}^l(SU(2)).
\tag\eqname{\vglCGseriesonA}
$$

%&&&&&  NEW SUBSECTION &&&&&&&&&&&&&&&&&&&&&&&&&&&&&&&&&&&&&&&&&&&&&&&&&
\subhead\newsubsection{\ssectHaarfSchuror}
Haar functional and Schur orthogonality relations
\endsubhead (\cite{\Woro}, \cite{\KoorOPTA})
On $\A$ there exists a functional which is an analogue of the invariant
integration on $SU(2)$. This functional, $h\colon\A\to\C$, called the Haar
functional, is uniquely determined by the properties

\noindent
(i) $h(1)=1$,

\noindent
(ii) $h(\xi^\ast\xi) \geq 0$ for all $\xi\in\A$,

\noindent
(iii) $\bigl( h\otimes id\bigr)\bigl(\Delta(a)\bigr) = h(a)1 =
\bigl( id\otimes h\bigr)\bigl(\Delta(a)\bigr)$.

\noindent
Condition (iii) is the left and right invariance of the Haar functional.

The Schur orthogonality for the matrix elements of the irreducible unitary
matrix corepresentations of $\A$ can be phrased as
$$
h\bigl( (t^l_{i,j})^\ast t^k_{n,m}\bigr) =
\delta_{l,k}\delta_{n,i}\delta_{m,j}q^{2(l-i)} {{1-q^2}\over{1-q^{4l+2}}}.
\tag\eqname{\vglstanmateltSchurorthrel}
$$

%&&&&&&&&&&&&&&&&&&&&&&&&&&&&&&&&&&&&
%%N E W   S E C T I O N%%%%%%%%%%%%**
%%%%%%%%%%%%%%%%%%%%%%%%%%%%%%%%%%%%%

\head\newsection . Koornwinder's $(\s,\t)$-spherical elements \endhead

In this section we treat the spherical functions on the quantum $SU(2)$
group introduced by Koornwinder \cite{\Koortwee} using the concept of
infinitesimal invariance. These results are the basis for the further
interpretation of the Askey-Wilson polynomials as generalised matrix
elements.

To discover which elements in the quantised universal enveloping algebra
$\U$ play the role of Lie algebra like elements we recall that in the
universal enveloping algebra ${\frak{U}}({\frak{sl}}(2,\C))$ the Lie
algebra elements satisfy, cf. \thetag{\vglLiealgelminU},
$\Delta(X)=1\otimes X + X\otimes S(1)$ for $X\in{\frak{sl}}(2,\C)$,
where $1$ satisfies $\Delta(1)=1\otimes 1$.
The elements in $\U$ satisfying $\Delta(X)=X\otimes X$, which are the
so-called group-like elements, are precisely $A^n$, $n\in\Z$, cf. Masuda
et al. \cite{\MasuMNNSUeen , lemma~1(i)}. The Lie algebra
${\frak{sl}}(2,\C)$ is three-dimensional over $\C$ and the elements in
$\U$ satisfying
$$
\Delta(X) = A^n\otimes X + X \otimes S(A^n)
$$
span a three-dimensional space if and only if $n=1$, cf.
\cite{\Koortwee , lemma~3.1}, \cite{\NoumMtwee , \S 1}. Moreover, the elements
satisfying $\Delta(X) = A\otimes X + X \otimes D$ are in the linear span
of $B$, $C$ and $A-D$.

For $\s\in\R$ we define
$$
X_\s = iq^{1\over 2}B-iq^{-{1\over 2}}C - {{q^\s-q^{-\s}}\over{q-q^{-1}}}
(A-D) \in\U .
\tag\eqname{\vgldsefXsigma}
$$
We also define
$$
X_\infty = D-A = \lim_{\s\to\infty} (q^{-1}-q)q^\s X_\s =
\lim_{\s\to -\infty} (q-q^{-1})q^{-\s} X_\s.
\tag\eqname{\vgldefXoneindig}
$$
Using \ssectquantunivenvalg\ it is easy to check that
$$
\Delta(X_\s) = A\otimes X_\s + X_\s\otimes D, \qquad S(X_\s)=-X_\s,
\qquad (X_\s A)^\ast = X_\s A.
\tag\eqname{\vglpropXsigma}
$$

An element $\xi\in\A$ is a $(\s ,\t )$-spherical element if
$$
X_\s .\xi = 0, \quad {\hbox{\rm and}}\quad \xi.X_\t = 0.
\tag\eqname{\vgldefspherical}
$$

\proclaim{Theorem \theoremname{\thmzonaleen}} {\rm (\cite{\Koortwee ,
prop.~4.7, thm.~5.3})}
The space of $(\s ,\t )$-spherical elements in the Hopf $\ast$-algebra
$\A$ is a commutative algebra generated by the element
$$
\eqalign{
\r = {1\over 2}\bigl(& \a^2+\d^2 +q\g^2 +q^{-1}\b^2 + i(q^{-\s}-q^\s )
(q\d\g +\b\a )\cr
& -i(q^{-\t}-q^\t)(\d\b +q\g\a) + (q^{-\s}-q^\s)(q^{-\t}-q^\t)\b\g \bigr)\cr}
\tag\eqname{\vgldefrhosigmatau}
$$
satisfying $\r^\ast=\r$, $\r\in{\Cal A}^1(SU(2))$.
The Haar functional on the subalgebra of $(\s ,\t )$-spherical elements is
given by
$$
h\bigl( p(\r )\bigr) = \int_\R p(x)\, dm(x;a,b,c,d\mid q^2)
\tag\eqname{\vgldefhaaropzonal}
$$
for any polynomial $p$, where $a=-q^{\s+\t+1}$, $b=-q^{-\s-\t+1}$,
$c=q^{\s-\t+1}$, $d=q^{-\s+\t+1}$ and $dm(x;a,b,c,d\mid q)$ denotes the
normalised Askey-Wilson measure, cf. \thetag{\vglnormalisedAWmeasure}.
\endproclaim

Since theorem \thmzonaleen\ is an important tool in our proof of theorem
7.5, %forward reference
we give the idea of the proof of theorem \thmzonaleen .
We start with the
observation that the $(\s,\t)$-spherical elements form a subalgebra of
$\A$, cf. proposition 6.4(i). %forward reference
Using the Peter-Weyl theorem it is
possible to prove that the $(\s,\t)$-spherical elements in $\Al$ are
completely
characterised by the kernels of the operators $t^l(X_\s)$ and
$t^l(X_\t^\ast)$, cf. the proof of proposition 6.4(ii). %forward reference
Koornwinder explicitly calculates the simplest non-trivial
$(\s,\t)$-spherical element $\r$ by considering the kernels of the
operators $t^1(X_\s)$
and $t^1(X_\t^\ast)$.

By an analysis of the spectrum of the operators $t^l(X_\s)$ and
$t^l(X_\t^\ast)$, which is essentially rephrased in proposition
\the\sectionno.2, %forward reference
it is shown that $\Al$, $l\in{1\over 2}+\Zp$, doesn't
contain $(\s,\t)$-spherical elements and that $\Al$, $l\in\Zp$,
contains a unique, up to
a scalar, $(\s,\t)$-spherical element, which has to be a polynomial
in $\r$. This
polynomial is identified with an Askey-Wilson polynomial with parameters as
in theorem \thmzonaleen\ by identifying the action of the Casimir operator
of $\U$ with the second order $q$-difference operator for the Askey-Wilson
polynomials.

This is done as follows. For this polynomial $p_l(\r)\in\Al$ we have
for $k\in\Z$ by \ssectreprthey\ $\langle A^k\Omega,p_l(\r)\rangle =
q^{1-2l}(1-q^{2l+1})^2(1-q^2)^{-2} p_l\bigl( {1\over2}(q^k+q^{-k})\bigr)$,
since testing against $A^k$ is an algebra homomorphism. By rewriting
$A^k\Omega$ as a linear combination of $A^{k+2}$, $A^k$ and $A^{k-2}$
modulo $\U X_\s + X_\t \U$ we see that this expression also equals a
linear combination of $p_l\bigl( {1\over2}(q^m+q^{-m})\bigr)$,
$m=k+2,k,k-2$. This gives precisely the second order $q$-difference
equation \cite{\AskeW , (5.7)} for the Askey-Wilson polynomials with the
parameters as in theorem \thmzonaleen , which has a unique polynomial
solution.

Equation \thetag{\vgldefhaaropzonal} follows by identifying
the Haar functional on the subalgebra of $(\s,\t)$-spherical elements
with the Askey-Wilson orthogonality measure via the Schur orthogonality
relations \thetag{\vglstanmateltSchurorthrel}.

\proclaim{Proposition \theoremname{\propspectrumdualqK}}
{\rm (\cite{\Koortwee , thm.~4.3, lemma~4.4})}
The self-adjoint operator $t^l(X_\s A)$ has an orthonormal basis of
eigenvectors $v^{l,j}(\s)=\sum_{n=-l}^l v^{l,j}_n(\s)e^l_n$
corresponding to the eigenvalue
$$
\l_j(\s ) = {{q^{-2j-\s}-q^{\s+2j}+q^\s - q^{-\s}}\over{q-q^{-1}}},
\quad j=-l,-l+1,\ldots,l.
$$
The coefficients $v^{l,j}_n(\s)$ are explicitly known by
$$
\align
v^{l,j}_n(\s) = &\, C^{l,j}(\s ) i^{n-l} q^{\s (l-n)}
q^{{1\over 2}(l-n)(l-n-1)}
\left( {{(q^{4l};q^{-2})_{l-n}}\over{(q^2;q^2)_{l-n}}}\right)^{1/2} \\
&\times R_{l-n}(q^{2j-2l}-q^{-2j-2l-2\s};q^{2\s},2l;q^2),
\endalign
$$
where $R_{l-n}$ is a dual $q$-Krawtchouk polynomial,
cf. \thetag{\vgldefdualqKrawtchoukpol}, and the constant
is given by
$$
C^{l,j}(\s) = q^{l+j} \left[ {{2l}\atop{l-j}}\right]_{q^2}^{1/2}
\left( {{1+q^{-4j-2\s}}\over{1+q^{-2\s}}}\right)^{1/2}
\Bigl( (-q^{2-2\s};q^2)_{l-j} (-q^{2+2\s};q^2)_{l+j}\Bigr)^{-1/2}.
$$
\endproclaim

To sketch the proof of this proposition we observe that the coefficients of
an eigenvector $\sum_{n=-l}^l c_ne^l_n$ of $t^l(X_\s A)$ satisfy a three-term
recurrence relation by \thetag{\vgldsefXsigma} and
\thetag{\vgldefreprU}, which can be identified with the three-term recurrence
relation for the dual $q$-Krawtchouk polynomials.

Since the vectors $v^{l,j}(\s)$ are orthonormal in a finite dimensional space,
we have the orthogonality relations
$$
\sum_{n=-l}^l v^{l,j}_n(\s) \overline{v^{l,i}_n(\s )} = \delta_{i,j},
\qquad
\sum_{i=-l}^l v^{l,i}_m(\s) \overline{v^{l,i}_n(\s )} = \delta_{n,m}.
\tag\eqname{\vglorthorelvector}
$$
The first part of \thetag{\vglorthorelvector} is equivalent to the
orthogonality relations \thetag{\vglorthorelqKrawtchouk}
for the $q$-Kraw\-tchouk polynomials, and
the second part of \thetag{\vglorthorelvector} is equivalent to the
orthogonality relations \thetag{\vglorthoreldualqKrawtchouk}
for the dual $q$-Krawtchouk polynomials.

The limit case $\s\to\infty$ of $t^l(X_\s A)v^{l,j}(\s)=\l_j(\s)
v^{l,j}(\s)$ reduces to $t^l(1-A^2)e^l_j = (1-q^{-2j})e^l_j$, cf.
\thetag{\vgldefXoneindig}, \thetag{\vgldefreprU}, because
$q^\s(q^{-1}-q)\l_j(\s) \to (1-q^{-2j})$ and
$v_n^{l,j}(\s)\to i^{l-j}\d_{j,n}$
as $\s\to\infty$. The last limit can be established using
transformation and summation formulas for basic hypergeometric series.

%&&&&&&&&&&&&&&&&&&&&&&&&&&&&&&&&&&&&
%%N E W   S E C T I O N%%%%%%%%%%%%**
%%%%%%%%%%%%%%%%%%%%%%%%%%%%%%%%%%%%%

\head\newsection . Generalised matrix elements\endhead

In this section we introduce the generalised matrix elements and we
derive various properties of these elements by employing the Hopf
$\ast$-algebra duality between the quantised polynomial algebra and the
quantised universal enveloping algebra. These properties, especially
proposition \the\sectionno.4, %forward reference
will be crucial in introducing orthogonal polynomials
in relation with these generalised matrix elements.

The vectors $v^{l,j}(\s)$, $j=-l,-l+1,\ldots,l$, of proposition
\propspectrumdualqK\ give an orthonormal basis for the representation space
of the $\ast$-representation $t^l$ of $\U$. We define linear functionals
on $\U$ by
$$
a^l_{i,j}(\t,\s) \bigl(X\bigr) = \langle t^l(X)v^{l,j}(\s),v^{l,i}(\t)\rangle,
\qquad \s,\t\in\R,\ i,j=-l,-l+1,\ldots,l,
\tag\eqname{\vgldefgeneralmatelt}
$$
where the inner product on the right hand side of
\thetag{\vgldefgeneralmatelt} is the inner product of the representation
space of $t^l$. It follows from proposition \propspectrumdualqK\ and
\ssectreprthey\  that
$$
a^l_{i,j}(\t,\s) = \sum_{n,m=-l}^l v^{l,j}_m(\s) \overline{v^{l,i}_n(\t)}
t^l_{n,m} \in\A,
\tag\eqname{\vglexprgenmateltinstanmatelt}
$$
so that we can write
$a^l_{i,j}(\t,\s) \bigl(X\bigr)=\langle X, a^l_{i,j}(\t,\s) \rangle$.
Note that $a^l_{n,m}(\infty,\infty)=i^{n-m}t^l_{n,m}$. So the
$a^l_{i,j}(\t,\s)$ are elements of the Hopf $\ast$-algebra $\A$ and we can
determine the action of the comultiplication, antipode and $\ast$-operator
on such an element.

\proclaim{Proposition \theoremname{\propHopfpropgenmatelt}}
The elements $a^l_{i,j}(\t,\s)$, $\s,\t\in\R$, $i,j=-l,-l+1,\ldots,l$,
$l\in\hZp$, defined in \thetag{\vgldefgeneralmatelt}, satisfy
$$
\gather
\Delta\bigl(a^l_{i,j}(\t,\s)\bigr) = \sum_{p=-l}^l a^l_{i,p}(\t,\mu)
\otimes a^l_{p,j}(\mu,\s), \quad \forall\,\mu\in\R,
\tag\eqname{\vglcomultongenmatelt} \\
\bigl(a^l_{i,j}(\t,\s)\bigr)^\ast = S\bigl(a^l_{j,i}(\s,\t)\bigr),\qquad
\varepsilon\bigl(a^l_{i,j}(\t,\s)\bigr) = \langle v^{l,j}(\s),
v^{l,i}(\t)\rangle  .
\tag\eqname{\vglSenepsiongenmatelt}
\endgather
$$
\endproclaim

\demo{Proof} These statements are proved by testing against appropriate
elements. Firstly,
$$
\align
\langle X\otimes Y, &\Delta\bigl(a^l_{i,j}(\t,\s)\bigr) \rangle =
\langle XY, a^l_{i,j}(\t,\s)\rangle =
\langle t^l(X)t^l(Y) v^{l,j}(\s), v^{l,i}(\t)\rangle \\
&= \sum_{p=-l}^l \langle t^l(X)v^{l,p}(\mu), v^{l,i}(\t)\rangle
\langle t^l(Y)v^{l,j}(\s),v^{l,p}(\mu)\rangle \\
&= \sum_{p=-l}^l \langle X, a^l_{i,p}(\t,\mu)\rangle
\langle Y, a^l_{p,j}(\mu,\s) \rangle , \qquad\forall\, X,Y\in\U,
\endalign
$$
by developing $t^l(Y)v^{l,j}(\s)$ in the basis $\{
v^{l,p}(\mu)\}_{p=-l,\ldots,l}$,
which proves \thetag{\vglcomultongenmatelt}. The first statement of
\thetag{\vglSenepsiongenmatelt} follows from
$$
\align
&\langle X,\bigl(a^l_{i,j}(\t,\s)\bigr)^\ast\rangle =
\overline{\langle S(X)^\ast, (a^l_{i,j}(\t,\s)\rangle} =
\overline{\langle t^l(S(X)^\ast)v^{l,j}(\s),v^{l,i}(\t)\rangle}\\
&=\langle t^l(S(X))v^{l,i}(\t),v^{l,j}(\s)\rangle
=\langle S(X), a^l_{j,i}(\s,\t)\rangle
=\langle X, S\bigl(a^l_{j,i}(\s,\t)\bigr)\rangle
\endalign
$$
for arbitrary $X\in\U$. Finally,
$$
\varepsilon\bigl( a^l_{i,j}(\t,\s)\bigr) = \langle 1, a^l_{i,j}(\t,\s)\rangle
= \langle v^{l,j}(\s), v^{l,i}(\t)\rangle
$$
proves the last statement. \qed\enddemo

Proposition \propHopfpropgenmatelt\ shows that by taking $\s=\t=\mu$ we obtain
an irreducible unitary matrix corepresentation
$a^l(\s)= \bigl(a^l_{i,j}(\s,\s)\bigl)_{i,j}$. This motivates to call
$a^l_{i,j}(\t,\s)$ generalised matrix elements. Of course,
$a^l(\s)$ is equivalent to $t^l$ by \ssectreprthey ,
so we can express the matrix elements in
terms of each other, cf. \thetag{\vglexprgenmateltinstanmatelt}.
But a more general expression holds, which is given in part (ii) of the
following corollary. Corollary \the\sectionno.2 %forward reference
also gives unitarity properties of the generalised matrix elements.

\proclaim{Corollary \theoremname{\corrHopfpropgenmatelt}}
(i) The matrix $a^l(\s)=\bigl(a^l_{i,j}(\s,\s)\bigl)_{i,j=-l,\ldots,l}$,
$l\in\hZp$, is an irredicible $(2l+1)$-dimensional unitary matrix
corepresentation of the Hopf $\ast$-algebra $\A$.

\noindent
(ii) The following connection between generalised matrix elements holds:
$$
a^l_{i,j}(\t,\s)= \sum_{n,m=-l}^l \langle v^{l,n}(\rho),v^{l,i}(\t)\rangle
\langle v^{l,j}(\s),v^{l,m}(\mu)\rangle a^l_{n,m}(\rho,\mu).
$$

\noindent
(iii) The generalised matrix elements satisfy the `unitarity property'
$$
\sum_{p=-l}^l a^l_{i,p}(\t,\mu)\bigl( a^l_{j,p}(\s,\mu)\bigr)^\ast =
\langle v^{l,j}(\s),v^{l,i}(\t)\rangle =
\sum_{p=-l}^l \bigl(a^l_{p,i}(\mu,\t)\bigr)^\ast a^l_{p,j}(\mu,\s).
$$
\endproclaim

\demo{Proof} Part (i) has been proved. To prove part (ii) we apply
$(\varepsilon\otimes id\otimes\varepsilon )\circ (\Delta\otimes id)$ to
\thetag{\vglcomultongenmatelt}. Because of the Hopf algebra axiom
$(\varepsilon \otimes id)\circ\Delta = id =
(id\otimes\varepsilon)\circ\Delta$, cf. \ssectHopfstaralgs ,
we see that the left hand side reduces to $a^l_{i,j}(\t,\s)$.
An application of \thetag{\vglcomultongenmatelt} shows that the right hand
side after this mapping yields for arbitrary $\rho\in\R$
$$
\sum_{p=-l}^l \sum_{r=-l}^l \varepsilon\bigl(a^l_{i,r}(\t,\rho)\bigr)
a^l_{r,p}(\rho,\mu) \varepsilon\bigl(a^l_{p,j}(\mu,\s)\bigr)
$$
and (ii) follows from \thetag{\vglSenepsiongenmatelt}.

To prove (iii) we apply $(id\otimes S)$ to
\thetag{\vglcomultongenmatelt}, and we use the Hopf algebra axiom
$m\circ(id\otimes S)=e\circ\varepsilon$, cf. \ssectHopfstaralgs .
The first equality follows from
\thetag{\vglSenepsiongenmatelt}. The second equality is proved similiarly
using the map $m\circ (S\otimes id)$. \qed\enddemo

For the generalised matrix elements we have the relative bi-invariance
property
$$
\align
&\langle (X_\t A)^\ast Y X_\s A, a^l_{i,j}(\t,\s)\rangle =
\langle t^l(Y) t^l(X_\s A)v^{l,j}(\s ), t^l(X_\t A)v^{l,i}(\t)\rangle \\
&= \lambda_j(\s)\lambda_i(\t) \langle t^l(Y) v^{l,j}(\s), v^{l,i}(\t)\rangle
=\lambda_j(\s)\lambda_i(\t)  \langle Y, a^l_{i,j}(\t,\s)\rangle ,
\quad \forall\,Y\in\U ,
\endalign
$$
by proposition \propspectrumdualqK .
Using \thetag{\vglweakactionUonA} we can reformulate this into the following
lemma, since $(X_\t A)^\ast= X_\t A$ by \thetag{\vglpropXsigma}.

\proclaim{Lemma \theoremname{\lemmarelbiinvgenmatelt}}
The generalised matrix elements $a^l_{i,j}(\t,\s)$ satisfy the
relative bi-invariance property
$$
(X_\s A).a^l_{i,j}(\t,\s)=\lambda_j(\s)\, a^l_{i,j}(\t,\s)
\quad{\hbox{\sl and}}\quad
a^l_{i,j}(\t,\s).(X_\t A) = \lambda_i(\t)\, a^l_{i,j}(\t,\s).
\tag\eqname{\vglbiinvariancegenmatelt}
$$
\endproclaim

Actually, lemma \lemmarelbiinvgenmatelt\ is a key observation in determining
the generalised matrix elements explicitly. But instead of determining
$a^l_{i,j}(\t,\s)$ we determine
$$
b^l_{i,j}(\t,\s)=A.a^l_{i,j}(\t,\s)
\tag\eqname{\vgldefdeformedgenmatelt}
$$
explicitly. Since $A.t^l_{n,m}=\sum_p t^l_{n,p}\langle t^l(A)e^l_m,
e^l_p\rangle = q^{-m}t^l_{n,m}$ we get from
\thetag{\vglexprgenmateltinstanmatelt} the explicit expression
$$
b^l_{i,j}(\t,\s) = \sum_{n,m=-l}^l v^{l,j}_m(\s) \overline{v^{l,i}_n(\t)}
q^{-m} t^l_{n,m} \in\Al \subset\A .
\tag\eqname{\vglexprdefgenmateltinstanmatelt}
$$

Note that $A.\colon\A\to\A$, $\xi\mapsto A.\xi$ is an algebra
homomorphism, because of \thetag{\vglactionUonprodA} and
\thetag{\vgldefDeltaonUq}. Moreover, this homomorphism is invertible, its
inverse being $D.\colon \A\to\A$, $\xi\mapsto D.\xi$. Lemma
\lemmarelbiinvgenmatelt\ can be rewritten as
$$
X_\s .b^l_{i,j}(\t,\s)=\lambda_j(\s)\, D.b^l_{i,j}(\t,\s)
\quad{\hbox{\rm and}}\quad
b^l_{i,j}(\t,\s).X_\t = \lambda_i(\t)\, b^l_{i,j}(\t,\s).D.
\tag\eqname{\vglbiinvariancedefgenmatelt}
$$
We have $\lambda_0(\s)=\lambda_0(\t)=0$, so that for $l\in\Z$
$b^l_{0,0}(\t,\s)$ is a $(\s,\t)$-spherical element
in the sense of \thetag{\vgldefspherical}.

\proclaim{Proposition \theoremname{\propeigdefmatelt}}
(i) Let $\xi\in\A$ be a $(\s,\t)$-spherical element, cf.
\thetag{\vgldefspherical}, and let $\eta\in\A$ satisfy
$$
X_\s .\eta=\lambda\, D.\eta
\quad{\hbox{\sl and}}\quad
\eta.X_\t = \mu\, \eta.D.
\tag\eqname{\vglgenbiinvariancet}
$$
for $\lambda,\mu\in\C$. Then $\eta\xi$ satisfies
\thetag{\vglgenbiinvariancet} for the same $\lambda$, $\mu$.
Moreover, if $\lambda,\mu\in\R$, then $\eta^\ast\eta$ is a $(\s,\t)$-spherical element.

\noindent
(ii) If $\eta\in\Al$ satisfies \thetag{\vglgenbiinvariancet} for
arbitrary $\lambda, \mu\in\C$ and $\eta$ is non-zero, then
$\lambda = \lambda_j(\s)$, $\mu=\lambda_i(\t)$ for some
$i,j\in\{-l,-l+1,\ldots,l\}$ and $\eta$ is a multiple of
$b^l_{i,j}(\t,\s)$.
\endproclaim

\demo{Proof} To prove (i) we first consider $\eta\xi$, then, by
\thetag{\vglactionUonprodA}, \thetag{\vglpropXsigma},
\thetag{\vglgenbiinvariancet}, \thetag{\vgldefspherical} and
\thetag{\vgldefDeltaonUq},
$$
X_\s .(\eta\xi) = (A.\eta)(X_\s .\xi)+(X_\s .\eta)(D.\xi) =
\lambda (D.\eta)(D.\xi) = \lambda\, D.(\eta\xi).
$$
Similarly we prove $(\eta\xi).X_\t = \mu\, (\eta\xi).D$.

To prove the other statement of (i) we proceed as before to obtain
$$
X_\s .(\eta^\ast\eta)= (A.\eta^\ast)(X_\s .\eta) + (X_\s .\eta^\ast)
(D.\eta) = (\lambda -\bar\lambda) (D.\eta)^\ast (D.\eta),
$$
by \thetag{\vglactionUonAandstar}, $S(A)^\ast=D$ and $S(X_\s)^\ast=-X_\s$.
This yields zero for $\lambda\in\R$. Similarly we prove
$(\eta^\ast\eta).X_\t=0$ for real $\mu$.

To prove (ii) we take $\eta$ of the form $\sum_{n,m=-l}^l \g_{n,m}t^l_{n,m}$,
then $\eta$ satisfies \thetag{\vglgenbiinvariancet} if and only if
$$
\align
t^l(AX_\s ) \sum_{m=-l}^l \g_{n,m}e^l_m &= \lambda \sum_{m=-l}^l \g_{n,m}e^l_m,
\quad \forall\, n\in\{ -l,-l+1,\ldots,l\} \\
t^l(X_\t A) \sum_{n=-l}^l \bar\g_{n,m}e^l_n &= \bar\mu \sum_{n=-l}^l
\bar\g_{n,m}e^l_n, \quad \forall\, m\in\{ -l,-l+1,\ldots,l\}
\endalign
$$
by the Peter-Weyl theorem for $\A$, cf. \ssectPWthmCGseries .
The operator $t^l(X_\t A)$ is
considered in proposition \propspectrumdualqK\ and the operator
$t^l(AX_\s )=t^l(A)t^l(X_\s A)t^l(D)$ is conjugated to such an operator.
Proposition \propspectrumdualqK\ then gives
$\g_{n,m}=c q^{-m}v^{l,j}_m(\s) \overline{v^{l,i}_n(\t)}$ for some
non-zero constant $c$ independent of $n$ and $m$, and
$\lambda=\lambda_j(\s)$, $\mu=\lambda_i(\t)$.\qed\enddemo

In the previous proposition we have shown that the product
of two elements can behave nicely with respect to the identity
\thetag{\vglgenbiinvariancet}.
Let us now consider $\eta\in\A$ satisfying \thetag{\vglgenbiinvariancet}
and multiply from the left by an arbitrary $\xi\in\A$. We get
similarly as in the proof of proposition \propeigdefmatelt
$$
\gather
X_\s.(\xi\eta) = (A.\xi)(X_\s.\eta) + (X_\s.\xi)(D.\eta) =
\bigl[ \l\, A.\xi + X_\s.\xi\bigr] (D.\eta), \\
(\xi\eta).X_\t  = (\xi.A)(\eta.X_\t) + (\xi.X_\t)(\eta.D) =
\bigl[ \mu\, \xi.A + \xi.X_\t\bigr] (\eta.D).
\endgather
$$
So, if $\xi$ satisfies
$$
\bigl[ \l\, A.\xi + X_\s.\xi\bigr] = \l_1\, D.\xi, \qquad
\bigl[ \mu\, \xi.A + \xi.X_\t\bigr] = \mu_1\, \xi.D
\tag\eqname{\vglrequironxi}
$$
for some $\l_1,\mu_1\in\C$, then we get
$$
X_\s.(\xi\eta)  = \l_1\, D.(\xi\eta), \qquad {\hbox{\rm and}} \qquad
(\xi\eta).X_\t = \mu_1\, (\xi\eta).D.
$$

Next we restrict $\xi\in{\Cal A}_q^{1/2}(SU(2))$, so that
$\xi = a\a+b\b+c\g+d\d$ for some complex constants $a$, $b$, $c$ and $d$.
From \ssectdualHopfalg\ we compute
$$
A.\pmatrix \a &\b\\ \g&\d \endpmatrix =
\pmatrix q^{1\over 2}\a &q^{-{1\over 2}}\b\\
q^{1\over 2}\g&q^{-{1\over 2}}\d \endpmatrix
, \quad B.\pmatrix \a &\b\\ \g&\d \endpmatrix =
\pmatrix 0&\a\\ 0&\g \endpmatrix
, \quad C.\pmatrix \a &\b\\ \g&\d \endpmatrix =
\pmatrix \b &0\\ \d&0 \endpmatrix.
$$
The first requirement of
\thetag{\vglrequironxi} is rewritten as
$$
\aligned
&\bigl[ q^{1/2} \l a - {{q^\s-q^{-\s}}\over{q-q^{-1}}} (q^{1/2}-q^{-1/2})a
+i q^{1/2} b - \l_1 q^{-1/2} a\bigr] \a + \\
&\bigl[ q^{-1/2} \l b - {{q^\s-q^{-\s}}\over{q-q^{-1}}} (q^{-1/2}-q^{1/2})b
-i q^{-1/2} a - \l_1 q^{1/2} b\bigr] \b + \\
&\bigl[ q^{1/2} \l c - {{q^\s-q^{-\s}}\over{q-q^{-1}}} (q^{1/2}-q^{-1/2})c
+i q^{1/2} d - \l_1 q^{-1/2} c\bigr] \g + \\
&\bigl[ q^{-1/2} \l d - {{q^\s-q^{-\s}}\over{q-q^{-1}}} (q^{-1/2}-q^{1/2})d
-i q^{-1/2} c - \l_1 q^{1/2} d\bigr] \d = 0.
\endaligned
\tag\eqname{\vglreqoneworkedout}
$$
Since the matrix elements are linearly independent in the algebra $\A$ we
see that each of the terms in square brackets in
\thetag{\vglreqoneworkedout} has to be zero, so we get two sets of two
equations in the unknowns $a$, $b$, $c$ and $d$. Or
$$
A \pmatrix a\\b\endpmatrix = \pmatrix 0\\0\endpmatrix, \qquad
A \pmatrix c\\d\endpmatrix = \pmatrix 0\\0\endpmatrix,
$$
where $A=A(\s,\l,\l_1)$ is the $2\times 2$-matrix given by
$$
A=A(\s,\l,\l_1) =
\left( {{ q^{1/2} \l - {{q^\s-q^{-\s}}\over{q-q^{-1}}} (q^{1/2}-q^{-1/2})
 - \l_1 q^{-1/2},\ iq^{1/2} }\atop{ -iq^{-1/2},\
 q^{-1/2} \l - {{q^\s-q^{-\s}}\over{q-q^{-1}}} (q^{-1/2}-q^{1/2})
- \l_1 q^{1/2} }} \right) .
\tag\eqname{\vgldefmatrixshifteen}
$$

There is only a non-trivial solution if the determinant of $A$ equals zero,
and this depends on the choice of $\l$ and $\l_1$.
From proposition \propeigdefmatelt (ii)
we see that without loss of generality we may assume that
$\l=\l_j(\s)$ and $\l_1=\l_n(\s)$ for some $j,n\in\hZp$. Now $\det A = 0$ if
and only if the product of the two diagonal elements equals $1$, which means
that the $\s$-dependence of this product must vanish. This gives
$n=j\pm {1\over 2}$. We find
$$
\gathered
A(\s,\l_j(\s),\l_{j+1/2}(\s)) =
\pmatrix q^{-2j-\s-1/2}, & iq^{1/2} \\
-iq^{-1/2}, & q^{\s+2j+1/2} \endpmatrix, \\
A(\s,\l_j(\s),\l_{j-1/2}(\s)) =
\pmatrix -q^{2j+\s-1/2}, & iq^{1/2} \\
-iq^{-1/2}, & -q^{1/2-\s-2j} \endpmatrix.
\endgathered
\tag\eqname{\vgldefmatrixshifttwee}
$$
Note that the matrices in \thetag{\vgldefmatrixshifttwee} only depend on
$\s+2j$, and that if $(x(\s),y(\s))^t$ is in the kernel of the first
matrix with $j=0$, then $(-\overline{x(-\s)},-\overline{y(-\s)})^t$ is
in the kernel of the second matrix of \thetag{\vgldefmatrixshifttwee}
with $j=0$.

The second requirement of \thetag{\vglrequironxi} leads analogously to
$$
B \pmatrix a\\c\endpmatrix = \pmatrix 0\\0\endpmatrix, \qquad
B \pmatrix b\\d\endpmatrix = \pmatrix 0\\0\endpmatrix,
$$
where $B^t=A(\t,\mu,\mu_1)$. So we only get non-trivial solutions for
$\mu=\l_i(\t)$, $\mu_1=\l_{i\pm 1/2}(\t)$. A similar remark on the relation
between the kernels of the two possible choices for $B$ holds here as well.
This proves the following proposition.

\proclaim{Proposition \theoremname{\propshifteigenval}}
Define elements in ${\Cal A}^{1/2}_q(SU(2))$ by
$$
\align
\a_{\t,\s} &= q^{1/2}\a-iq^{\s-1/2}\b +iq^{\t+1/2}\g + q^{\s+\t-1/2}\d, \\
\b_{\t,\s} &= -q^{\s+1/2}\a-iq^{-1/2}\b -iq^{\s+\t+1/2}\g + q^{\t-1/2}\d, \\
\g_{\t,\s} &= -q^{\t+1/2}\a+iq^{\t+\s-1/2}\b +iq^{1/2}\g + q^{\s-1/2}\d, \\
\d_{\t,\s} &= q^{\t+\s+1/2}\a+iq^{\t-1/2}\b -iq^{\s+1/2}\g + q^{-1/2}\d,
\endalign
$$
and let $\eta$ satisfy \thetag{\vglgenbiinvariancet} with $\l=\l_j(\s)$ and
$\mu=\l_i(\t)$, then

\noindent
(i) $\a_{\t+2i,\s+2j}\eta$ satisfies \thetag{\vglgenbiinvariancet} with
$\l=\l_{j-1/2}(\s)$ and $\mu=\l_{i-1/2}(\t)$,

\noindent
(ii) $\b_{\t+2i,\s+2j}\eta$ satisfies \thetag{\vglgenbiinvariancet} with
$\l=\l_{j+1/2}(\s)$ and $\mu=\l_{i-1/2}(\t)$,

\noindent
(iii) $\g_{\t+2i,\s+2j}\eta$ satisfies \thetag{\vglgenbiinvariancet} with
$\l=\l_{j-1/2}(\s)$ and $\mu=\l_{i+1/2}(\t)$,

\noindent
(iv) $\d_{\t+2i,\s+2j}\eta$ satisfies \thetag{\vglgenbiinvariancet} with
$\l=\l_{j+1/2}(\s)$ and $\mu=\l_{i+1/2}(\t)$.
\endproclaim

\demo{Remark} The normalisation for the elements $A(\t,\s)$, $B(\t,\s)$,
$C(\t,\s)$ and $D(\t,\s)$ has been chosen such that
$$
{1\over{\sqrt{(1+q^{2\s})(1+q^{2\t})} }}
\pmatrix \a_{\t,\s},& \b_{\t,\s}\\ \g_{\t,\s},&\d_{\t,\s} \endpmatrix =
b^{1/2}(\t,\s) =
\pmatrix b^{1/2}_{-1/2,-1/2}, & b^{1/2}_{-1/2,1/2}\\
b^{1/2}_{1/2,-1/2}, & b^{1/2}_{1/2,1/2}\endpmatrix.
$$
\enddemo

In particular if we take $\eta=b^l_{i,j}(\t,\s)$ we can write each of these
products as a sum of two similar elements by the Clebsch-Gordan series,
cf. \ssectPWthmCGseries .
In the extremal cases, i.e. $l=\max\lbrace\vert i\vert, |j|\rbrace$, this
sum reduces to one element and we get recursion relations. For instance,
if $l=i$, we obtain
$$
\gathered
\d_{\t+2l,\s+2m} b^l_{l,m} (\t,\s) = c\, b^{l+1/2}_{l+1/2,m+1/2} (\t,\s), \\
\g_{\t+2l,\s+2m} b^l_{l,m} (\t,\s) = c\,
b^{l+1/2}_{l+1/2,m-1/2} (\t,\s)
\endgathered
\tag\eqname{\vgliterposs}
$$
for some non-zero constants. These relations can be iterated in several ways.
Since the
elements are non-commuting we have to be careful about the ordering in
a product. We stick to the convention that
$\Pi_{j=0}^k \xi_k = \xi_0\xi_1\ldots \xi_k$.

Before giving the explicit
expression and the corresponding proof we derive symmetries for the
generalised matrix elements $a^l_{i,j}(\t,\s)$ enabling us to prove only
one case. First we introduce the
algebra map $\Psi\colon\A\to\A$ defined by interchanging $\beta$ and
$\gamma$. It follows directly from \thetag{4.1} that $\Psi$ does preserve
the commutation relations. Moreover, $\Psi(t^l_{n,m}) = t^l_{m,n}$ as follows
directly from \thetag{4.20}. This has already been observed by Koornwinder
\cite{\Koordrie} before knowing the explicit expressions for these
matrix elements. We also introduce the antilinear mapping
${}^-\colon\A\to\A$ defined by taking complex conjugates of all
coefficients.

\proclaim{Lemma \theoremname{\lemmsymmgenmateltas}} The generalised matrix
elements $a^l_{i,j}(\t,\s)$ satisfy the symmetry relations
$$
\aligned
a^l_{i,j}(\t,\s) & = \Psi\bigl( (a^l_{j,i}(\s,\t)^-\bigr) \\
                 & = (a^l_{-i,-j}(-\t,-\s))^- \\
                 & = \Psi\bigl( a^l_{-j,-i}(-\s,-\t)\bigr).
\endaligned
$$
\endproclaim

\demo{Proof} The last relation follows from the first two. The first relation
follows from $\Psi(t^l_{n,m})=t^l_{m,n}$ and \thetag{6.2}. The second relation
follows from \thetag{6.2} and
$$
v^{l,j}_n(\s) = \overline{v^{l,-j}_n(-\s)}.
\tag\eqname{\vglsymmvljsigmaenminsigma}
$$
To see that \thetag{\vglsymmvljsigmaenminsigma} is true we first note
that $C^{l,j}(\s)=C^{l,-j}(-\s)$. Next the transformation formula, cf.
\cite{\GaspR, (3.2.3)},
$$
{}_3\vp_2 \left( {{q^{-n},b,c}\atop{d,e}};q,q\right) =
{{(de/bc;q)_n}\over{(e;q)_n}} \left({{bc}\over d}\right)^n
{}_3\vp_2 \left( {{q^{-n},d/b,d/c}\atop{d,de/bc}};q,q\right)
$$
yields the identity
$$
\align
R_{l-n}(q^{2j-2l}-&q^{-2j-2l-2\s};q^{2\s},2l;q^2) = \\
&(-1)^{l-n} q^{-2\s(l-n)} R_{l-n}(q^{-2j-2l}-q^{2j-2l+2\s};q^{-2\s},2l;q^2),
\endalign
$$
for the dual $q$-Krawtchouk polynomial
$R_{l-n}$ in proposition \propspectrumdualqK\
and this proves \thetag{\vglsymmvljsigmaenminsigma}.
\qed\enddemo

In the following corollary we only work out one of the many possibilities for
one of the four cases. A similar expression for the other three cases can be
be obtained using lemma \lemmsymmgenmateltas\ and the observation
$\Psi(A.\xi)=\xi.A$.

\proclaim{Corollary \theoremname{\coroexplexprminelt}} With the notation of
proposition \propspectrumdualqK\ we have
$$
\align b^l_{l,m}(\t,\s)&= E^l_m(\t,\s)
\prod_{k=0}^{l+m-1} \d_{\t+2l-1-k,\s+2m-1-k} \\
&\qquad \times \prod_{j=0}^{l-m-1} \g_{\t+l-m-1-j,\s-l+m+1+j},
\endalign
$$
with
$$
E^l_m(\t,\s) = C^{l,m}(\s) C^{l,l}(\t) q^{\s (m-l)}
q^{{1\over 2}(l-m)(l-m-1)}.
$$
\endproclaim

\demo{Remark} This corollary has been stated by Noumi and Mimachi in an
unpublished announcement that extends the announcement \cite{\NoumM}.
\enddemo

\demo{Proof} Iteration of \thetag{\vgliterposs} proves the existence of the
product as in the corollary. It remains to determine the constant. For
this we apply the one-dimensional $\ast$-representation $\hp$ to it and
we compare the coefficients of $e^{-il\theta}$ on both sides. The
left hand side gives $C^{l,m}(\s) C^{l,l}(\t) q^{-l}$ and the right
hand side follows from proposition \propshifteigenval\ and we find
$E^l_m(\t,\s) q^{-l} q^{\s(l-m)}
q^{-{1\over 2}(l-m)(l-m-1)}$, from which the corollary follows.
\qed\enddemo

%&&&&&&&&&&&&&&&&&&&&&&&&&&&&&&&&&&&&
%%N E W   S E C T I O N%%%%%%%%%%%%**
%%%%%%%%%%%%%%%%%%%%%%%%%%%%%%%%%%%%%

\head\newsection . Generalised matrix elements and Askey-Wilson polynomials
\endhead

We have now developed all the necessary ingredients for the
interpretation of Askey-Wilson polynomials on the quantum $SU(2)$ group.
In this section we first show that orthogonal polynomials are of
importance in describing generalised matrix elements. Next these
polynomials are explicitly calculated in terms of Askey-Wilson
polynomials, giving a full proof of the theorem announced by Noumi and
Mimachi \cite{\NoumM , thm.~3}, see also \cite{\NoumMtwee , thm.~4}.

\proclaim{Theorem \theoremname{\thmabstrorthpolgenmatelt}}
For fixed $i,j\in\hZp$ such that $i-j\in\Z$, there exists a system of
orthogonal polynomials $(p_k)_{k\in\Zp}$ of degree $k$ such that for
$l\geq m =\max (\vert i\vert,\vert j\vert)$, $l-m\in\Zp$,
$$
b^l_{i,j}(\t,\s) = b^m_{i,j}(\t,\s) p_{l-m}(\r ),
\tag\eqname{\vglabstrorthpolgenmatelt}
$$
where $\r \in\A$ is the self-adjoint element given in
\thetag{\vgldefrhosigmatau}.
\endproclaim

\demo{Proof} We first prove that an expression of the form
\thetag{\vglabstrorthpolgenmatelt} exists. Consider for any polynomial
$s_{l-m}$ of degree $l-m$ the expression
$b^m_{i,j}(\t,\s) s_{l-m}(\r )$. If we decompose this product with respect to
the decomposition of $\A$ in \thetag{\vglPeterWeylforA} we get
$$
b^m_{i,j}(\t,\s) s_{l-m}(\r ) = \sum_{k=\vert 2m-l\vert}^l b^k,
\qquad b^k \in {\Cal A}_q^k(SU(2)),
\tag\eqname{\vgleenproofthmabstrorthpolgenmatelt}
$$
by the Clebsch-Gordan series \thetag{\vglCGseriesonA}.
Now \thetag{\vgldefunitarycorep} implies
that fo any $X\in\U$ the mappings $X.$ and $.X$ preserve
${\Cal A}_q^k(SU(2))$. Proposition \propeigdefmatelt (i) shows that the left
hand side \thetag{\vgleenproofthmabstrorthpolgenmatelt} satisfies
\thetag{\vglgenbiinvariancet} with $\lambda=\lambda_j(\s)$ and
$\mu=\lambda_i(\t)$. Consequently, each $b^k$ has to satisfy
\thetag{\vglgenbiinvariancet} with $\lambda=\lambda_j(\s)$ and
$\mu=\lambda_i(\t)$. Proposition \propeigdefmatelt (ii) implies that
$b^k=0$ for $k<m$ and $b^k=c_kb^k_{i,j}(\t,\s)$ for $k\geq m$ and some
constants $c_k$. Hence, \thetag{\vgleenproofthmabstrorthpolgenmatelt}
reduces to
$$
b^m_{i,j}(\t,\s) s_{l-m}(\r ) = \sum_{k=m}^l c^k b^k_{i,j}(\t,\s).
$$
Since both sides contain the same degree of freedom,
\thetag{\vglabstrorthpolgenmatelt} follows once we know that the mapping
$s_{l-m}\mapsto b^m_{i,j}(\t,\s)s_{l-m}(\r )$ is injective. This can be seen
by applying the one-dimensional $\ast$-representation $\p$ of $\A$, cf.
\thetag{\vgldefonedimreprA},
and use of the explicit expression of $b^m_{i,j}(\t,\s)$, cf.
\thetag{\vglexprdefgenmateltinstanmatelt}. This proves that an expression
as in \thetag{\vglabstrorthpolgenmatelt} exists.

For $l,k\geq m$, $l-m,k-m\in\Zp$ we have $b^l_{i,j}(\t,\s)\in\Al$,
$b^k_{i,j}(\t,\s)\in{\Cal A}_q^k(SU(2))$, so that the Schur orthogonality
relations for the Haar functional $h$, cf.
\thetag{\vglstanmateltSchurorthrel}, imply
$$
h\Bigl( \bigl(b^l_{i,j}(\t,\s)\bigr)^\ast b^k_{i,j}(\t,\s)\Bigr) =
\delta_{k,l}h_l, \qquad h_l>0 .
\tag\eqname{\vgltweeproofthmabstrorthpolgenmatelt}
$$
Now $\bigl(b^m_{i,j}(\t,\s)\bigr)^\ast b^m_{i,j}(\t,\s)=w_m(\r)$ for some
polynomial $w_m$ of degree $2m$,
by proposition \propeigdefmatelt (i), theorem \thmzonaleen\ and the
Clebsch-Go\-rdan series, cf. \thetag{\vglCGseriesonA}.
Hence, \thetag{\vgltweeproofthmabstrorthpolgenmatelt} equals
$$
h\bigm( \bar p_{l-m}(\r)p_{k-m}(\r) w_m(\r)\bigr) = \delta_{l,k} h_l,
\qquad h_l >0,
$$
since we have already established \thetag{\vglabstrorthpolgenmatelt} for
some polynomial and since $\r^\ast=\r$. Consequently,
by theorem \thmzonaleen\ the polynomials
$p_{l-m}$, $l-m\in\Zp$, form a system of orthogonal polynomials with
respect to the moment functional ${\Cal L}$ with moments given by
$$
{\Cal L}[x^n] = \int_\R x^n w_m(x) \, dm(x;a,b,c,d\mid q^2) <\infty,
\tag\eqname{\vgldefmomentfunct}
$$
where $dm(x;a,b,c,d\mid q^2)$ is the measure described in theorem
\thmzonaleen . \qed\enddemo

\demo{Remark \theoremname{\remquadnormothproof}}
(i) The constant $h_l$ in \thetag{\vgltweeproofthmabstrorthpolgenmatelt}
can be given explicitly using the Schur orthogonality relations
\thetag{\vglstanmateltSchurorthrel} and
\thetag{\vglexprdefgenmateltinstanmatelt}. This yields
$$
\aligned
h_l &= {{(1-q^2)q^{2l}}\over{1-q^{4l+2}}} \sum_{m=-l}^l
\vert v^{l,j}_m(\s)\vert^2 q^{-2m} \sum_{n=-l}^l
\vert v^{l,i}_n(\t)\vert^2 q^{-2n} \\
&= {{(1-q^2)q^{2l}}\over{1-q^{4l+2}}}
\t_{\sqrt{q}}\bigr( b^l_{j,j}(\s,\s)\bigl)
\t_{\sqrt{q}}\bigr( b^l_{i,i}(\t,\t)\bigl) \\
&= {{(1-q^2)q^{2l}}\over{1-q^{4l+2}}}
\t_{\sqrt{q}}\bigr( b^{\vert j\vert}_{j,j}(\s,\s)\bigl)
r_{l-\vert j\vert}\bigl( (q+q^{-1})/2\bigr)
\t_{\sqrt{q}}\bigr( b^{\vert i\vert}_{i,i}(\t,\t)\bigl)
s_{l-\vert i\vert}\bigl( (q+q^{-1})/2\bigr),
\endaligned
\tag\eqname{\vglsquarenormformeen}
$$
for some polynomials $r_{l-\vert j\vert}$ and $s_{l-\vert i\vert}$
given by theorem \thmabstrorthpolgenmatelt .
Here $\t_{\sqrt{q}}$ is the one-dimensional representation of $\A$
defined in \ssectreprsA .

\noindent
(ii) Another way of seeing that the polynomials in
\thetag{\vglabstrorthpolgenmatelt} form a set of orthogonal polynomials
once the the existence of a polynomial satisfying
\thetag{\vglabstrorthpolgenmatelt} is established, is the following.
Multiply both sides of \thetag{\vglabstrorthpolgenmatelt} from the right
by $\r$ to get
$$
\gather
 b^m_{i,j}(\t,\s) p_{l-m}(\r )\r = b^l_{i,j}(\t,\s)\r = \\
 A_{l-m}b^{l+1}_{i,j}(\t,\s) + B_{l-m}b^l_{i,j}(\t,\s) +
C_{l-m}b^{l-1}_{i,j}(\t,\s) = \\
 A_{l-m}b^m_{i,j}(\t,\s) p_{l-m+1}(\r ) + B_{l-m}
b^m_{i,j}(\t,\s) p_{l-m}(\r ) + C_{l-m} b^m_{i,j}(\t,\s) p_{l-m-1}(\r ),
\endgather
$$
by the Clebsch-Gordan series, \ssectPWthmCGseries , and
\thetag{\vglabstrorthpolgenmatelt}.
Using the injectivity of the map
$s_{l-m}\mapsto b^m_{i,j}(\t,\s)s_{l-m}(\r )$ proves that the
polynomials satisfy a three-term recurrence relation implying the
orthogonality. See also remark \the\sectionno .6(i). %forward reference
\enddemo

In order to be able to identify the the system of orthogonal polynomials
described in theorem \thmabstrorthpolgenmatelt , we have to calculate
the polynomial $w_m$. Apply the one-dimensional $\ast$-representation
$\hp$ to the identity
$w_m(\r)= \bigl(b^m_{i,j}(\t,\s)\bigr)^\ast b^m_{i,j}(\t,\s)$ to obtain
$$
w_m(\cos\theta)= \bigl\vert \hp\bigl( b^m_{i,j}(\t,\s)\bigr)\bigr\vert^2.
\tag\eqname{\vgldefwm}
$$
Hence, it is sufficient to calculate $\hp\bigl(b^m_{i,j}(\t,\s)\bigr)$.
These values follow from the following proposition by taking
$\l=e^{i\theta/2}$, cf. \ssectreprsA .

\proclaim{Proposition \theoremname{\propexplvalueonedimonmingenmatelt}}
For $i,j\in\hZp$, $i-j\in\Z$ and $m=\max (\vert i\vert,\vert j\vert)$
we have

\noindent
(i) In case $m=i$ or $-i\leq j\leq i$
$$
\tl \bigl( b^i_{i,j}(\t,\s)\bigr) = C^{i,j}(\s)C^{i,i}(\t) q^{-i}
\l^{-2i} (\l^2q^{1+\t-\s};q^2)_{i-j}
(-\l^2q^{1+\t+\s};q^2)_{i+j}.
$$

\noindent
(ii) In case $m=j$ or $-j\leq i\leq j$
$$
\tl \bigl( b^j_{i,j}(\t,\s)\bigr) = C^{j,i}(\t)C^{j,j}(\s) q^{-j}
\l^{-2j} (\l^2q^{1+\s-\t};q^2)_{j-i}
(-\l^2q^{1+\t+\s};q^2)_{i+j}.
$$

\noindent
(iii) In case $m=-i$ or $i\leq j\leq -i$
$$
\tl \bigl( b^{-i}_{i,j}(\t,\s)\bigr) = C^{-i,-j}(-\s)C^{-i,-i}(-\t) q^i
\l^{2i} (\l^2q^{1-\t+\s};q^2)_{j-i}
(-\l^2q^{1-\t-\s};q^2)_{-i-j}.
$$

\noindent
(iv) In case $m=-j$ or $j\leq i\leq -j$
$$
\tl \bigl( b^{-j}_{i,j}(\t,\s)\bigr) = C^{-j,-i}(-\t)C^{-j,-j}(-\s) q^j
\l^{2j} (\l^2q^{1-\s+\t};q^2)_{i-j}
(-\l^2q^{1-\t-\s};q^2)_{-i-j}.
$$
\endproclaim

\demo{Proof} First we observe that
the function $\tl \bigl( b^l_{i,j}(\t,\s)\bigr)$ of $\l$ satisfies
the symmetry relations
$$
\aligned
\tl \bigl( b^l_{i,j}(\t,\s)\bigr) &=
\overline{ \t_{\bar \l} \bigl( b^l_{j,i}(\s,\t)\bigr)} \\
&= \tl \bigl( b^l_{-j,-i}(-\s,-\t)\bigr)\\
&= \overline{ \t_{\bar \l} \bigl( b^l_{-i,-j}(-\t,-\s)\bigr)} .
\endaligned
\tag\eqname{\vglsymmetriesgenmatelt}
$$
This a straightforward consequence of lemma \lemmsymmgenmateltas\ and
and $\tl \circ \Psi=\tl$.

By \thetag{\vglsymmetriesgenmatelt}
it suffices to prove the first statement.
Application of $\tl$ to corollary \coroexplexprminelt\ gives
$$
\align
\tl \bigl( b^l_{l,m}(\t,\s)\bigr) & = E^l_m(\t,\s)
\prod_{k=0}^{l+m-1} (q^{\t+\s+2l+2m-2-2k+1/2}\l + q^{-1/2}\l^{-1}) \\
&\qquad\qquad\times
\prod_{j=0}^{l-m-1} (-q^{\t+l-m-1-j+1/2}\l + q^{\s-l+m+1+j-1/2}\l^{-1}) \\
& = C^{l,m}(\s)C^{l,l}(\t) q^{-l} \l^{-2l} (\l^2q^{1+\t-\s};q^2)_{l-m}
(-\l^2q^{1+\t+\s};q^2)_{l+m}
\endalign
$$
and this proves (i). \qed\enddemo

\proclaim{Corollary \theoremname{\corogendualqKraw}}
A generating function for the dual $q$-Krawtchouk polynomials is
$$
\aligned
\sum_{n=0}^N &t^n q^{n(N+\s)/2} {{(q^{-N};q)_n}\over{(q;q)_n}}
R_n(q^{-x}-q^{x-N-\s};q^\s,N;q) \\
&= (-tq^{-(N+\s)/2};q)_x (tq^{(\s-N)/2};q)_{N-x}
\endaligned
\tag\eqname{\vglgenfunctiedualqKrawtchouk}
$$
for $N\in\N$, $x\in\{ 0,\ldots, N\}$.
\endproclaim

\demo{Remark}
This generating function can be obtained from the generating function for the
Askey-Wilson polynomials given by Ismail and Wilson \cite{\IsmaW, (1.9)}, cf.
\cite{\Koel, (3.3)} for this derivation. In \cite{\Koel} the generating
function \thetag{\vglgenfunctiedualqKrawtchouk} is used to prove
proposition \propexplvalueonedimonmingenmatelt\ in a special case.
\enddemo

\demo{Proof} Apply the one-dimensional representation $\tl$ to the expression
for $b^i_{i,j}(\t,\s)$ from \thetag{\vglexprdefgenmateltinstanmatelt}. Since
this also equals proposition \propexplvalueonedimonmingenmatelt (i) we obtain
\thetag{\vglgenfunctiedualqKrawtchouk} with $q$, $N$, $n$, $x$, and $t$ replaced
by $q^2$, $2i$, $i-n$, $i-j$ and $-\l^2q^{\t+2i+1}$.
\qed\enddemo

We can now prove the main result of this section in which we relate
Askey-Wilson polynomials to generalised matrix elements. The Askey-Wilson
polynomials involve four continuous parameters. In the following theorem we
establish an interpretation of the Askey-Wilson polynomials with two
continuous and two discrete parameters.

\proclaim{Theorem \theoremname{\thmdefgenmateltAskeyWpol}}
For $i,j\in\hZp$, $i-j\in\Z$ and $l-m\in\Zp$, $m=\max (\vert i\vert,\vert
j\vert)$ we have

\noindent
(i) In case $m=i$ or $-i\leq j\leq i$:
$$
b^l_{i,j}(\t,\s) = d^l_{i,j}(\t,\s) b^i_{i,j}(\t,\s)
p^{(i-j,i+j)}_{l-i}(\r;q^\t, q^\s \mid q^2).
$$

\noindent
(ii) In case $m=j$ or $-j\leq i\leq j$:
$$
b^l_{i,j}(\t,\s) = d^l_{j,i}(\s,\t) b^j_{i,j}(\t,\s)
p^{(j-i,i+j)}_{l-j}(\r;q^\s, q^\t \mid q^2).
$$

\noindent
(iii) In case $m=-i$ or $i\leq j\leq -i$:
$$
b^l_{i,j}(\t,\s) = d^l_{-i,-j}(-\t,-\s) b^{-i}_{i,j}(\t,\s)
p^{(j-i,-i-j)}_{l+i}(\r;q^{-\t}, q^{-\s} \mid q^2).
$$

\noindent
(iv) In case $m=-j$ or $j\leq i\leq -j$:
$$
b^l_{i,j}(\t,\s) = d^l_{-j,-i}(-\s,-\t) b^{-j}_{i,j}(\t,\s)
p^{(i-j,-i-j)}_{l+j}(\r;q^{-\s}, q^{-\t} \mid q^2).
$$

\noindent
Here the constant is given by
$$
d^l_{i,j}(\t,\s) = {{C^{l,j}(\s)C^{l,i}(\t)}\over{C^{i,j}(\s)C^{i,i}(\t)}}
{{q^{i-l}}\over{(q^{4l};q^{-2})_{l-i}}} .
$$
\endproclaim

\demo{Proof} The explicit form \thetag{\vglabstrorthpolgenmatelt} given
in theorem \thmabstrorthpolgenmatelt , the symmetry relations of
\thetag{\vglsymmetriesgenmatelt} and $\hp (\r) =\cos\theta$ being
independent of $\s$, $\t$ show that it suffices to prove the first statement.

From the explicit form of the Askey-Wilson weight measure, cf.
proposition \propAWorthogonality ,
we immediately get
$$
(az,a/z;q)_r dm(x;a,b,c,d\mid q) = {{ (ab,ac,ad;q)_r}\over{(abcd;q)_r}}
\, dm(aq^r,b,c,d\mid q)
$$
for $r\in\Zp$, $x=(z+z^{-1})/2$. A double application shows that for
$r,s\in\Zp$, $x=(z+z^{-1})/2$,
$$
\align
&(az,a/z;q)_r (dz,d/z;q)_s dm(x;a,b,c,d\mid q) = \\
&(ab,ac;q)_r(bd,cd;q)_s {{(ad;q)_{r+s}}\over{(abcd;q)_{r+s}}} \,
dm(aq^r,b,c,dq^s\mid q).
\endalign
$$
Use this in conjuction with \thetag{\vgldefhaaropzonal},
\thetag{\vgldefwm}, proposition \propexplvalueonedimonmingenmatelt (i)
and \thetag{\vgldefmomentfunct} to find that the looked-for polynomials
are multiples of the Askey-Wilson polynomials
$$
p_{l-i}(\r;-q^{\s+\t+1+2i+2j}, -q^{-\s-\t+1}, q^{\s-\t+1}, q^{\t-\s+1+2i-2j}
\mid q^2),
$$
which we rewrite using \thetag{\vgldefAWalsqJacobi}.

It remains to calculate the constant. We apply the one-dimensional
$\ast$-representation $\hp$ to both sides of (i), and next we compare the
coefficient of $e^{-il\theta}$ on both sides. The coefficient of
$e^{-il\theta}$ on the left hand side is
$v^{l,i}_l(\s)\overline{v^{l,j}_l(\t)}q^{-l} = C^{l,j}(\s)C^{l,i}(\t)q^{-l}$.
The coefficient of $e^{-i(l-i)\theta}$ of $p_{l-i}$ is
$(q^{2l+2i+2};q^2)_{l-i}= (q^{4l};q^{-2})_{l-i}$, cf. \cite{\AskeW, p.5}, so
that the coefficient of $e^{-il\theta}$ on the right hand side
equals
$C^{i,j}(\s)C^{i,i}(\t)q^{-i}(q^{4l};q^{-2})_{l-i}$, from which we obtain the
value for $d^l_{i,j}(\t,\s)$. \qed\enddemo

\demo{Remark \theoremname{\remotherproofsortpols}} (i) The orthogonal
polynomials are determined in theorem \thmdefgenmateltAskeyWpol\ by
explicitly determining the orthogonality measure for these polynomials
as is done in \cite{\Koel} to determine associated spherical elements.
Since there are other characteristics of orthogonal polynomials, theorem
\thmdefgenmateltAskeyWpol\ can also be proved in other ways. We will
sketch two other possible lines of proof.

The first alternative uses the three-term recurrence relation for
orthogonal polynomials. Since $\r\in{\Cal A}^1_q(SU(2))$ we obtain from
the Clebsch-Gordan series, cf. \ssectPWthmCGseries , and proposition
\propeigdefmatelt\ the relation
$$
b^l_{i,j}(\t,\s)\r = A_lb^{l+1}_{i,j}(\t,\s) + B_l b^l_{i,j}(\t,\s)
+C_l b^{l-1}_{i,j}(\t,\s)
\tag\eqname{\vgleenproofviarecurrence}
$$
for certain constants $A_l$, $B_l$, $C_l$, cf. the proof of theorem
\thmabstrorthpolgenmatelt\ and remark \remquadnormothproof (ii).
Applying the one-dimensional
$\ast$-representation $\hp$ to \thetag{\vgleenproofviarecurrence} and
using \thetag{\vglexprdefgenmateltinstanmatelt} yields an identity for
trigonometric polynomials from which $A_l,B_l,C_l\in\C$ can be
determined explicitly by comparing the coefficients of
$e^{i(l+1)\theta}$, $e^{il\theta}$$e^{i(l-1)\theta}$. It follows that
the polynomials of theorem \thmabstrorthpolgenmatelt\ satisfy the
recurrence relation
$$
\cos\theta p_{l-m+1}(\cos\theta) = A_l p_{l-m+1}(\cos\theta)
+ B_l p_{l-m}(\cos\theta) + C_l p_{l-m-1}(\cos\theta)
\tag\eqname{\vgltweeproofviarecurrence}
$$
with initial conditions $p_{-1}(\cos\theta)=0$, $p_0(\cos\theta)=1$.
From this we can obtain the desired expression of $p_{l-m}$ in terms of
the Askey-Wilson polynomials by identifying
\thetag{\vgltweeproofviarecurrence} with the three-term recurrence
relation for the Askey-Wilson polynomials, cf. \cite{\AskeW ,
(1.24)}.

The second alternative proof uses the second order $q$-difference
equation for the Askey-Wilson polynomials as the identifying
characteristic. This approach is also used by
Koornwinder \cite{\Koortwee} for the $(\s,\t)$-spherical elements and
the proof was sketched in \S 5. This time we have to use
$$
\langle A^k\Omega , b^l_{i,j}(\t,\s)\rangle =
q^{1-2l}
{{(1-q^{2l+1})^2}\over{(1-q^2)^2}} \langle A^k, b^m_{i,j}(\t,\s)\rangle
p_{l-m}\bigl( (q^k+q^{-k})/2\bigr),
$$
where $\Omega$ is the Casimir element of $\U$, cf.
\thetag{\vgldefCasimir}. But $\langle A^k,b^m_{i,j}(\t,\s)\rangle$ is
known by proposition \propexplvalueonedimonmingenmatelt , since $\langle
A^k,\xi\rangle=\t_{q^{k/2}}(\xi)$ for all $\xi\in\A$. By rewriting
$A^k\Omega$ as a linear combination of $A^{k+2}$, $A^k$, $A^{k-2}$
modulo $\U (AX_s-\l_j(\s)) + (X_\t A-\l_i(\t))\U)$ and using the
relative bi-invariance property
\thetag{\vglbiinvariancedefgenmatelt} we get a relation of the form
$$
\gather
q^{1-2l}
{{(1-q^{2l+1})^2}\over{(1-q^2)^2}} \langle A^k, b^m_{i,j}(\t,\s)\rangle
p_{l-m}\bigl( (q^k+q^{-k})/2\bigr) =\\
\sum_{p=-1}^1 C_p \langle A^{k+2p}, b^m_{i,j}(\t,\s)\rangle
p_{l-m}\bigl( (q^{k+2p}+q^{-k-2p})/2\bigr),
\endgather
$$
where the value of $C_p$ follows from the rewriting of $A^k\Omega$. The
resulting identity can be related to the second order $q$-difference
equation for the Askey-Wilson polynomials, cf. \cite{\AskeW , (5.7)},
from which the identification of the $p_{l-m}$'s with the Askey-Wilson
polynomials follows up to a scalar.

It must be noted that these two alternatives do not use the a priori
knowledge of the explicit form of the Haar functional on
$(\s,\t)$-spherical elements as given in theorem \thmzonaleen .
The second alternative is an extension of the method sketched in \S 5 to
arbitrary matrix elements. It would be instructive to have a proof of the
expression of the Haar functional as in theorem \thmzonaleen\ without using
the explicit expression of the zonal spherical elements in terms of the
Askey-Wilson polynomials.

\noindent
(ii) The polynomials in theorem \thmdefgenmateltAskeyWpol\ are
identified using the orthogonality measure and the constant involved is
calculated by comparing certain coefficients. So we have not used
\thetag{\vglsquarenormformeen}, but the expression in
\thetag{\vglsquarenormformeen} is explicitly known by theorem
\thmdefgenmateltAskeyWpol\  and proposition
\propexplvalueonedimonmingenmatelt . On the other hand, now that we have
identified the polynomials as Askey-Wilson polynomials, the value of
$h_l$ can also be read off from proposition \propAWorthogonality . It is
straighforward to check that these values agree.
\enddemo

Now that we have a full description of the matrix elements $b^l_{i,j}(\t,\s)$
we can state the corresponding results on the generalised matrix elements
$a^l_{i,j}(\t,\s)$ defined in \thetag{\vgldefgeneralmatelt}. The first
part of each of the four statements of the following corollary is
equivalent to the theorem announced by Noumi and Mimachi \cite{\NoumM ,
thm.~3}, \cite{\NoumMtwee , thm.~4}.

\proclaim{Corollary \theoremname{\corexplgenmateltAWpols}}
For $i,j\in\hZp$, $i-j\in\Z$ and $l-m\in\Zp$, $m=\max (\vert i\vert,\vert
j\vert)$ and the same constant $d^l_{i,j}(\t,\s)$ as in theorem
\thmdefgenmateltAskeyWpol\ we have

\noindent
(i) In case $m=i$ or $-i\leq j\leq i$:
$$
\align
a^l_{i,j}(\t,\s) &= d^l_{i,j}(\t,\s) a^i_{i,j}(\t,\s)
p^{(i-j,i+j)}_{l-i}(D.\r;q^\t, q^\s \mid q^2), \\
\tl \bigl( a^i_{i,j}(\t,\s)\bigr) &= C^{i,j}(\s)C^{i,i}(\t)
\l^{-2i} (\l^2q^{\t-\s};q^2)_{i-j}
(-\l^2q^{\t+\s};q^2)_{i+j}.
\endalign
$$

\noindent
(ii) In case $m=j$ or $-j\leq i\leq j$:
$$
\align
a^l_{i,j}(\t,\s) &= d^l_{j,i}(\s,\t) a^j_{i,j}(\t,\s)
p^{(j-i,i+j)}_{l-j}(D.\r;q^\s, q^\t \mid q^2), \\
\tl \bigl( a^j_{i,j}(\t,\s)\bigr) &= C^{j,i}(\t)C^{j,j}(\s)
\l^{-2j} (\l^2q^{\s-\t};q^2)_{j-i}
(-\l^2q^{\t+\s};q^2)_{i+j}.
\endalign
$$

\noindent
(iii) In case $m=-i$ or $i\leq j\leq -i$:
$$
\align
a^l_{i,j}(\t,\s) &= d^l_{-i,-j}(-\t,-\s) a^{-i}_{i,j}(\t,\s)
p^{(j-i,-i-j)}_{l+i}(D.\r;q^{-\t}, q^{-\s} \mid q^2), \\
\tl \bigl( a^{-i}_{i,j}(\t,\s)\bigr) &= C^{-i,-j}(-\s)C^{-i,-i}(-\t)
\l^{2i} (\l^2q^{-\t+\s};q^2)_{j-i}
(-\l^2q^{-\t-\s};q^2)_{-i-j}.
\endalign
$$

\noindent
(iv) In case $m=-j$ or $j\leq i\leq -j$:
$$
\align
a^l_{i,j}(\t,\s) &= d^l_{-j,-i}(-\s,-\t) a^{-j}_{i,j}(\t,\s)
p^{(i-j,-i-j)}_{l+j}(D.\r;q^{-\s}, q^{-\t} \mid q^2), \\
\tl \bigl( a^{-j}_{i,j}(\t,\s)\bigr) &= C^{-j,-i}(-\t)C^{-j,-j}(-\s)
\l^{2j} (\l^2q^{-\s+\t};q^2)_{i-j}
(-\l^2q^{-\t-\s};q^2)_{-i-j}.
\endalign
$$
\endproclaim

\demo{Proof} As already remarked, $D.b^l_{i,j}(\t,\s)=a^l_{i,j}(\t,\s)$,
cf. \thetag{\vgldefdeformedgenmatelt}.
Moreover, $D.$ is an algebra
homomorphism, so the first statements in (i)--(iv) are nothing but
theorem \thmdefgenmateltAskeyWpol .

To prove the other statements we first observe that from the explicit
duality \thetag{\vglexpldualityongenAU} between $\A$ and $\U$ we get
$$
D.\a = q^{-1/2}\a,\quad D.\b=q^{-1/2}\b,\quad D.\g=q^{1/2}\g,
\quad D.\d=q^{1/2}\d.
$$
Hence, $\tl(D.\xi)=\t_{\l q^{-1/2}}(\xi)$ for any $\xi\in\A$ since $D.$
is an algebra homomorphism. This observation and
proposition \propexplvalueonedimonmingenmatelt\ prove the other
statements. \qed\enddemo

The special case $\s,\t\to\infty$ of corollary \corexplgenmateltAWpols\
corresponds precisely to \thetag{\vglmatrixeltalslittleqJacobi} if we
use the transition of Askey-Wilson polynomials to the little $q$-Jacobi
polynomials, cf. \cite{\Koortwee , \S 6}. Note that the only property
of the matrix elements $t^l_{n,m}\in\A$ that we have used in the proofs
until now is
$\tl(t^l_{n,m})=\d_{n,m}\l^{-2n}$, which can also be replaced by
$\langle A^k,t^l_{n,m}\rangle = \d_{n,m} q^{-nk}$, cf.
\thetag{\vgldefreprU}. So we can obtain
\thetag{\vglmatrixeltalslittleqJacobi} as a special case of
corollary \corexplgenmateltAWpols\ if we can determine the minimal
elements $a^m_{i,j}(\infty,\infty)$, $m=\max(\vert i\vert , \vert
j\vert)$, explicitly. This is possible using \thetag{\vgldefreprU} and
the explicit duality given for suitable bases of the underlying linear
spaces of $\A$ and $\U$, cf. Masuda et al.
\cite{\MasuMNNSUeen , lemma~5}.

%&&&&&&&&&&&&&&&&&&&&&&&&&&&&&&&&&&&&
%%N E W   S E C T I O N%%%%%%%%%%%%**
%%%%%%%%%%%%%%%%%%%%%%%%%%%%%%%%%%%%%

\head\newsection . Some applications\endhead

Before giving some new applications of the interpretation of the
Askey-Wilson polynomials on the quantum $SU(2)$ group as established in
the previous sections we give three recent applications. These
applications concern addition formulas for $q$-analogues of the Legendre
polynomial, the action of the quantum $SU(2)$ group on quantum spheres
and non-negative linearisation coefficients for products of
$q$-Legendre polynomials. The proof in \S \the\sectionno .3 is somewhat
different from the proof given by Koornwinder \cite{\KoorHG , \S 7}.
For detailed results and proofs the reader may consult the references
cited in the following subsections.

%&&&&&  NEW SUBSECTION &&&&&&&&&&&&&&&&&&&&&&&&&&&&&&&&&&&&&&&&&&&&&&&&&
\subhead\newsubsection{\ssectadditionformqLegendre}
Addition formulas for $q$-Legendre polynomials\endsubhead
The group theoretic proof of the addition formula for Legendre
polynomials is a consequence of the classical group counterpart of
\thetag{\vglcomultongenmatelt} for $i=j=0$, cf. Vilenkin \cite{\Vile ,
Ch.~3, \S 4.2}, Vilenkin and Klimyk \cite{\VileK , Vol.~1, \S 6.6.2}.
The addition formula for Legendre polynomials can also be found in Askey
\cite{\Aske , Lecture~4} and Erd\'elyi et al. \cite{\ErdeHTF , Vol.~2,
10.11(47)}.
However in the quantum group case \thetag{\vglcomultongenmatelt} for
$i=j=0$ does not readily give addition formulas for $q$-analogues of the
Legendre polynomials, since \thetag{\vglcomultongenmatelt} is an
identity involving non-commuting variables. So we have to use
representations of $\A$, cf. \ssectreprsA , to represent
\thetag{\vglcomultongenmatelt} for $i=j=0$ as an identity involving
$q$-special functions with ordinary arguments.

Application of the one-dimensional representation
$\t_{\sqrt{\l}}\otimes \t_{\sqrt{\nu}}$ to \thetag{\vglcomultongenmatelt} for
$i=j=0$ yields the following addition formula for a two-parameter family
of $q$-Legendre polynomials, cf. Noumi and Mimachi \cite{\NoumM ,
(4.3)}, see also Koelink \cite{\Koel , (3.15)};
$$
\aligned
&(q^2;q^2)_l q^{-l} p_l^{(0,0)}(\xi(\l\nu);q^\t,q^\s\mid q^2) =
{{p_l^{(0,0)}(\xi(\l);q^\mu,q^\t \mid q^2)
p_l^{(0,0)}(\xi(\nu);q^\mu,q^\s \mid q^2)}\over
{(-q^{2-2\mu},-q^{2+2\mu};q^2)_l}} \\
& + \sum_{p=1}^l {{(1+q^{4p+2\mu})(q^2;q^2)_{l+p} (\l\nu)^{-p}
(\l q^{\mu-\t},-\l q^{\t+\mu},\nu q^{\mu-\s}, -\nu
q^{\mu+\s};q^2)_p}\over{(1+q^{2\mu}) (q^2;q^2)_{l-p}
(-q^{2-2\mu};q^2)_{l-p} (-q^{2+2\mu};q^2)_{l+p}}} \\
&\qquad\qquad \times
p_{l-p}^{(p,p)}(\xi(\l);q^\mu,q^\t \mid q^2)
p_{l-p}^{(p,p)}(\xi(\nu);q^\mu,q^\s \mid q^2) \\
& + \sum_{p=1}^l {{(1+q^{4p-2\mu})(q^2;q^2)_{l+p} (\l\nu)^{-p}
(\l q^{\t-\mu},-\l q^{-\t-\mu},\nu q^{\s-\mu}, -\nu
q^{-\mu-\s};q^2)_p}\over{(1+q^{-2\mu}) (q^2;q^2)_{l-p}
(-q^{2+2\mu};q^2)_{l-p} (-q^{2-2\mu};q^2)_{l+p}}}\\
&\qquad\qquad \times
p_{l-p}^{(p,p)}(\xi(\l);q^{-\mu},q^{-\t} \mid q^2)
p_{l-p}^{(p,p)}(\xi(\nu);q^{-\mu},q^{-\s} \mid q^2),
\endaligned
\tag\eqname{\vgladdformqLegtwoparam}
$$
with $\xi(\l) = {1\over 2}(q^{-1}\l + q \l^{-1})$ by corollary
\corexplgenmateltAWpols\ and straightforward calculation.

This addition formula is probably not the most general addition formula
which can be derived for the $q$-Legendre polynomial
$p_l^{(0,0)}(x;q^\t,q^\s\mid q^2)$ for two  reasons. The first reason is
that the Rahman-Verma addition formula for the continuous $q$-Legendre
polynomial, cf. \cite{\RahmV , (1.24) with $a=q^{1/4}$}
corresponding to the case $\s=\t=0$ is not contained in
\thetag{\vgladdformqLegtwoparam}. The second reason is that the mapping
$\t_{\sqrt{\l}}\otimes \t_{\sqrt{\nu}}$ has a very large kernel, so that
this map kills a lot of information contained in
\thetag{\vglcomultongenmatelt}. However, it is possible to obtain the
addition formula for Legendre polynomials in full generality from
\thetag{\vgladdformqLegtwoparam} by replacing $q^\t$ by $t$ before
letting $q\uparrow 1$.

In case we take $\s=\t=\mu=0$, $i=j=0$ in \thetag{\vglcomultongenmatelt}
we obtain an identity for continuous $q$-ultraspherical polynomials from
which it is possible to derive the Rahman-Verma addition formula for the
continuous $q$-Legendre polynomial, cf. Koelink \cite{\Koel , \S 4}. The
method uses one-dimensional representations of the Hopf $\ast$-algebra
$\A$, but it also makes fundamental use of the relative bi-invariance
properties, cf. lemma \lemmarelbiinvgenmatelt , to introduce an extra
independent variable.

The first example of an addition formula for $q$-Legendre polynomials
derived using quantum group techniques is Koornwinder's \cite{\KoorAF}
addition formula for the little $q$-Legendre polynomials. Now we start
from \thetag{\vglcomultongenmatelt}  with $i=j=0$ and $\s=\t=\mu=\infty$
with the explicit form of the matrix elements given by
\thetag{\vglmatrixeltalslittleqJacobi}. Application of the infinite
dimensional $\ast$-representation $\pi^\infty_0\otimes \pi^\infty_0$,
cf. \ssectreprsA , yields an identity for bounded operators on the
Hilbert space $\ell^2(\Zp) \hat\otimes \ell^2(\Zp)$. Letting act the
operators on suitable vectors and taking inner products leads to the
addition formula for the little $q$-Legendre polynomials.

Another way of deriving addition formulas for $q$-Legendre polynomials
starts from corollary \corrHopfpropgenmatelt (ii). If we specialise
$i=j=0$, $\rho=\mu=\infty$ and we apply the $\ast$-representation
$\hp$ to the resulting identity, then we get, after a suitable
rescaling,
the Fourier series
$$
p_l^{(0,0)}(\cos\theta;q^\s,q^\t\mid q^2) = C \sum_{n=-l}^l v^{l,0}_n(\s)
\overline{v^{l,0}_n} q^{-n} e^{-in\theta},
\tag\eqname{\vglaltaddformqLeg}
$$
for some explicitly known constant $C$,
cf. Koornwinder \cite{\Koortwee , (5.9)}. Koornwinder \cite{\Koortwee ,
remark~5.4} also announces a generalisation of
\thetag{\vglaltaddformqLeg} to the corresponding $q$-ultraspherical
polynomials. To see that \thetag{\vglaltaddformqLeg} can be viewed as a
$q$-analogue of the addition formula for the Legendre polynomials we
recall that the $v^{l,0}_n(\s)$, cf. proposition \propspectrumdualqK , is
given by a dual $q$-Krawtchouk polynomial which tends to the Krawtchouk
polynomial \thetag{\vgldefKrawtchouk} as $q\uparrow 1$.
Now the Krawtchouk polynomial can be
rewritten in terms of a Jacobi polynomial, cf. e.g. Nikiforov and Uvarov
\cite{\NikiU , \S 12, \S 22}, and then we obtain the addition formula for
the Legendre polynomials in the limit $q\uparrow 1$.

Another special case of corollary \corrHopfpropgenmatelt (ii), namely
$\rho=\mu=\t=\infty$, is turned into an addition formula for the big
$q$-Legendre polynomials by use of an infinite dimensional
$\ast$-representation of $\A$, cf. Koelink \cite{\KoelAF}.

%&&&&&  NEW SUBSECTION &&&&&&&&&&&&&&&&&&&&&&&&&&&&&&&&&&&&&&&&&&&&&&&&&
\subhead\newsubsection{\ssectactionquantumspheres}
Action on quantum spheres\endsubhead
The Lie group $SU(2)$ acts on the two-dimensional sphere $S^2$. This can
be seen by noting that $SU(2)$ is a double cover of $SO(3)$, the group
of $3\times 3$ orthogonal matrices with determinant $1$, i.e.
$SU(2)/\{\pm I\} \cong SO(3)$, and by letting $SO(3)$ act on $S^2$ in
the natural way. We can identify $S^2 \cong SU(2)/K$ with
$K=S(U(1)\times U(1))$ as in \S 2. So $SU(2)$ also acts on the algebra
of functions on $S^2$.
Here we consider the analogue of these actions by letting act  $\A$ on
certain deformed function algebras of quantum spheres originally
introduced by Podle\'s \cite{\Podl}. As an alternative for this
subsection there is an action of the quantised universal enveloping
algebra $\U$ on the real two-dimensional sphere $S^2$, cf. Rideau and
Winternitz \cite{\RideW} for details.

Define the subspace $\AS = \{\xi\in\A \mid X_\s.\xi=0\}$ of
$\A$ of right invariant elements with respect to $X_\s\in\U$, cf.
\thetag{\vgldsefXsigma}. Applying similar techniques as in the proof of
proposition \propeigdefmatelt (i) we see that $\AS$ is
a $\ast$-subalgebra of $\A$. From the coassociativity axiom, cf.
\ssectHopfstaralgs , we get $(id\otimes X.)\Delta(\xi) = \Delta(X.\xi)$
for any $X\in\U$, $\xi\in\A$, cf. \thetag{\vgldefactionUonA}. This
implies $\Delta \colon \AS \to \A \otimes \AS$.

A left corepresentation
of the Hopf $\ast$-algebra $\A$ in a linear space $V$ is a map $L\colon
V\to \A\otimes V$ such that $(id\otimes L)\circ L = (\Delta\otimes
id)\circ L$ and $(\varepsilon \otimes id)\circ L=id$. If we choose a basis
in a left corepresentation space $V$, then the matrix elements in $\A$
form a matrix corepresentation of $\A$, cf. \ssectreprthey . The
Hopf algebra axioms, cf. \ssectHopfstaralgs ,
imply that we have a left corepresentation
in $\AS$ of the Hopf $\ast$-algebra $\A$ with the left coaction $L$ given
by the comultiplication $\Delta$. This is the analogue of the left
regular representation of $SU(2)$ on $L^2(S^2) \cong L^2(SU(2)/K)$.

It is a consequence of proposition \propeigdefmatelt (ii) that
${\Cal A}^l_{q,\s}(S^2) = \AS\cap\Al$ equals $\{ 0\}$ for $l\in{1\over
2}+\Zp$ and that ${\Cal A}^l_{q,\s}(S^2)$ equals the
$(2l+1)$-dimensional subspace spanned by $b^l_{i,0}(\t,\s)$,
$i=-l,-l+1,\ldots,l$, for $l\in\Zp$ and for any $\t\in\R\cup\{\infty\}$.
So, by \ssectreprthey , the corepresentation of $\A$ in $\AS$ splits
multiplicity free as the direct sum of the finite dimensional matrix
corepresentations $t^l$, $l\in\Zp$.

Using $(id\otimes X.)\Delta(\xi) = \Delta(X.\xi)$ with $X=A$ and
\thetag{\vglcomultongenmatelt} we get
$$
\Delta\bigl( b^l_{i,0}(\t,\s)\bigr) = \sum_{p=-l}^l
a^l_{i,p}(\t,\mu)\otimes b^l_{p,0}(\mu,\s),
\tag\eqname{\vglactiononStwo}
$$
so that we recover the generalised matrix elements in this realisation
of the irreducible representation of spin $l\in\Zp$ where the associated
spherical elements play the role of basis in $\AS$.

Two seperate cases of \thetag{\vglactiononStwo} have been worked out by
Noumi and Mimachi before this interpretation of the Askey-Wilson
polynomials on the quantum $SU(2)$ group, namely for $\t=\mu=\infty$,
yielding the associated spherical elements as a basis of $\AS$ in terms of
big $q$-Jacobi polynomials with $\a=\b\in\Zp$, cf.
\cite{\NoumMCMP}, and for $\t=\s=\mu=0$, yielding the associated spherical
elements as a basis of ${\Cal A}_{q,0}(S^2)$ in terms of the continuous
$q$-ultraspherical polynomials, cf. \cite{\NoumMDMJ}.

So we can view $\AS$ as the deformed algebra of polynomials on the sphere
$S^2$. It is possible to identify $\AS$ with one of the quantum
$2$-spheres of Podle\'s \cite{\Podl}, cf. Noumi and Mimachi
\cite{\NoumMtwee , prop.~2} (without proof) and Dijkhuizen and
Koornwinder \cite{\DijkK , \S 2}. Moreover, each of Podle\'s's quantum
$2$-spheres can be obtained in this way.

%&&&&&  NEW SUBSECTION &&&&&&&&&&&&&&&&&&&&&&&&&&&&&&&&&&&&&&&&&&&&&&&&&
\subhead\newsubsection{\ssectpositivityresults}
Non-negative linearisation coefficients for $q$-Legendre polynomials
\endsubhead
In this subsection we consider the linearisation coefficients of products
of $q$-Legendre polynomials corresponding to the matrix elements
$a^l_{0,0}(\s,\s)$. The linearisation coefficients are shown to be
non-negative as already proved by Koornwinder \cite{\KoorHG , \S 7} using
the framework of compact quantum Gelfand pairs. See also Floris
\cite{\Flor} and Vainerman \cite{\Vain}. Here we give a somewhat
different proof of this fact using Clebsch-Gordan coefficients.

Let $V^l$ be the $(2l+1)$-dimensional vector space in which the
representation $t^l$ of $\U$ acts, then we have for any $X\in\U$, cf.
\ssectPWthmCGseries ,
$$
(t^{l_1}\otimes t^{l_2})\circ \Delta(X) = C^\ast\Bigl(
{\sum_{l=\vert l_1-l_2\vert}^{l_1+l_2}} \hskip-0.3truecm {}^\oplus\ \,
t^l(X) \Bigr) C
$$
where $C\colon V^{l_1}\otimes V^{l_2}\to \bigoplus_{
l=\vert l_1-l_2\vert}^{l_1+l_2} V^l$ is a unitary intertwining
operator. Its matrix coefficients with respect to the standard bases
are known in terms of $q$-Hahn polynomials, cf. \cite{\KoelK}. Here we
consider $C$ with respect to the orthonormal basis
$v^{l_1,j_1}(\s)\otimes v^{l_2,j_2}(\s)$, $j_1=-l_1,\ldots,l_1$,
$j_2=-l_2,\ldots, l_2$ of $V^{l_1}\otimes V^{l_2}$ and the orthonormal
basis $v^{l,j}(\s)$, $j=-l,\ldots,l$, $l=\vert l_1-l_2\vert,\ldots,
l_1+l_2$ of $\bigoplus_{l=\vert l_1-l_2\vert}^{l_1+l_2} V^l$. Then $C$
becomes a unitary matrix with matrix coefficients defined by
$$
C \colon v^{l_1,j_1}(\s)\otimes v^{l_2,j_2}(\s) \mapsto
\sum_{l=\vert l_1-l_2\vert}^{l_1+l_2} \sum_{j=-l}^l
C^{l_1,l_2,l}_{j_1,j_2,j}(\s) v^{l,j}(\s).
$$

If we now take $l_1,l_2\in\Zp$ we see that $v^{l_1,0}(\s)\otimes
v^{l_2,0}(\s)$ belongs to the kernel of $(t^{l_1}\otimes t^{l_2})
\Delta(X_\s A)$. Since $C$ is an intertwiner we obtain from proposition
\propspectrumdualqK
$$
0= \sum_{l=\vert l_1-l_2\vert}^{l_1+l_2} \sum_{j=-l}^l
C^{l_1,l_2,l}_{0,0,j}(\s) \l_j(\s) v^{l,j}(\s).
$$
Since $v^{l,j}(\s)$ forms a basis and $\l_j(\s)\not= 0$ for $j\not= 0$
we get $ C^{l_1,l_2,l}_{0,0,j}(\s) =0$ for $j\not= 0$.

Using the intertwining property we obtain
$$
\aligned
& \sum_{l=\vert l_1-l_2\vert }^{l_1+l_2} \Bigl\vert
C^{l_1,l_2,l}_{0,0,0}(\s)\Bigr\vert^2 \langle t^l(X) v^{l,0}(\s),
v^{l,0}(\s)\rangle = \\
& \langle \bigl( (t^{l_1}\otimes t^{l_2})\circ\Delta(X) \bigr)
v^{l_1,0}(\s)\otimes v^{l_2,0}(\s),
v^{l_1,0}(\s)\otimes v^{l_2,0}(\s) \rangle = \\
& \sum_{(X)} \langle t^{l_1}(X_{(1)}v^{l_1,0}(\s), v^{l_1,0}(\s)\rangle
\langle t^{l_2}(X_{(2)})v^{l_2,0}(\s), v^{l_1,0}(\s) \rangle,
\endaligned
$$
where $\Delta(X) = \sum_{(X)}X_{(1)}\otimes X_{(2)}$. Using the duality
between $\A$ and $\U$ and the notation of \S 6 we have
$$
\sum_{l=\vert l_1-l_2\vert }^{l_1+l_2} \Bigl\vert
C^{l_1,l_2,l}_{0,0,0}(\s)\Bigr\vert^2 \langle a^l_{0,0}(\s,\s), X\rangle
= \langle a^{l_1}_{0,0}(\s,\s)a^{l_2}_{0,0}(\s,\s), X\rangle.
$$
Apply the one-dimensional representation $\t_{e^{i\theta/2}q^{1/2}}$
to the resulting identity in $\A$ and use
corollary \corexplgenmateltAWpols\ to find
$$
p_{l_1}^{(0,0)} (\cos\theta;q^\s,q^\s\mid q^2)
p_{l_2}^{(0,0)} (\cos\theta;q^\s,q^\s\mid q^2) =
\sum_{l=\vert l_1-l_2\vert}^{l_1+l_2} c_l(l_1,l_2)
p_l^{(0,0)} (\cos\theta;q^\s,q^\s\mid q^2)
\tag\eqname{\vgllinearisatieqLegendre}
$$
with non-negative coefficients $c_l(l_1,l_2) \geq 0$. This also remains
valid for the little $q$-Legendre polynomials corresponding to
$\s=\infty$, but then we also have an explicit formula for the
linearisation coefficient in terms of the square of a certain $q$-Hahn
polynomial, since the Clebsch-Gordan coefficient is explicitly known in
that case, cf. \cite{\KoelK}.

The linearisation coefficients in \thetag{\vgllinearisatieqLegendre}
for $\s=0$, corresponding to
the continuous $q$-ultraspherical polynomials, were explicitly given by
Rogers in 1894, cf. \cite{\AskeI , (4.18)}, \cite{\AskeW , (4.8)},
\cite{\GaspR , \S 8.5}, and non-negativity can be read off from that
expression.

Finally, we note that $\A$ and the $\ast$-invariant coideal $J$ of $\U$
generated by $X_\s A$ is a compact quantum Gelfand pair in the sense of
Floris \cite{\Flor} and Koornwinder \cite{\KoorHG , \S 5}.

%&&&&&&&&&&&&&&&&&&&&&&&&&&&&&&&&&&&&
%%N E W   S E C T I O N%%%%%%%%%%%%**
%%%%%%%%%%%%%%%%%%%%%%%%%%%%%%%%%%%%%

\head\newsection . Discrete orthogonality relations\endhead

In this section we translate some of the unitarity properties of the
generalised matrix elements $a^l_{i,j}(\t,\s)$, cf. corollary
\corrHopfpropgenmatelt (iii), into identities for $q$-special functions.
We obtain a formula for dual $q$-Krawtchouk polynomials which contains
as special cases the generating function
\thetag{\vglgenfunctiedualqKrawtchouk} and the Poisson kernel for the dual
$q$-Krawtchouk polynomials and the orthogonality relations for
$q$-Krawtchouk polynomials. Discrete orthogonality relations for certain
functions extending the orthogonality relations for the dual
$q$-Krawtchouk polynomials are also obtained.

\proclaim{Proposition \theoremname{\propAWnaarqKraw}}
For $l\in\hZp$, $i,j=-l,-l+1,\ldots,l$ we have
$$
\aligned
\tl \bigl( a^l_{i,j}(\t,\s)\bigr) = &(-1)^{l-i} C^{l,j}(\s)C^{l,i}(\t)
q^{(l-i)(l-i-1)}q^{(\s-\t)(l-i)} \l^{-2i}
(\l^{-2}q^{\t-\s+2-2l+2i};q^2)_{l-i}\\
&\times (\l^2q^{\t-\s+2i-2l};q^2)_{l-j}
{{(-\l^2q^{\t+\s};q^2)_\infty}\over{(-\l^2q^{\t+\s+2i+2j};q^2)_\infty}}\\
&\times \,{}_4\vp_3 \left( {{q^{-2(l-i)},q^{-2(l-j)},-q^{2\t-2l+2i},
-q^{-2\s-2l-2j}}\atop{q^{-4l},q^{2-\s+\t-2l+2i}\l^{-2},
q^{\t-\s-2l+2i}\l^2}};q^2,q^2\right).
\endaligned
$$
\endproclaim

\demo{Proof} We have to show that each of the four expressions in
corollary \corexplgenmateltAWpols\ can be rewritten in this form. In
the cases (i) and (ii) of corollary \corexplgenmateltAWpols\ we
rewrite
the terminating $_4\vp_3$-series of the Askey-Wilson polynomial by reversing
the series summation. So we use, cf. \cite{\GaspR, ex.~1.4.(ii)}
$$
\align
&{}_4\vp_3 \left( {{q^{-n},abcdq^{n-1},az,a/z}\atop{ab,ac,ad}};q,q\right)
= \tag\eqname{\vglrevseriesAWpols} \\
{{(abcdq^{n-1},az,a/z;q)_n}\over{(ab,ac,ad;q)_n}}& (-1)^nq^{-{1\over 2}n(n-1)}
{}_4\vp_3 \left( {{q^{-n}, q^{1-n}/(ab), q^{1-n}/(ac), q^{1-n}/(ad)}\atop
{q^{2-2n}/(abcd),q^{1-n}/(az),q^{1-n}z/a}};q,q\right),
\endalign
$$
with $q$ replaced by $q^2$.
In case (i) we use \thetag{\vglrevseriesAWpols} with $a=q^{1+\s-\t}$ and
in case (ii) we use \thetag{\vglrevseriesAWpols}
with $a=q^{1+\s-\t+2j-2i}$. The values of the other
parameters in the Askey-Wilson polynomials follow from corollary
\corexplgenmateltAWpols\ and \thetag{\vgldefAWalsqJacobi}. Note that we
use the symmetry of the Askey-Wilson polynomials in its four parameters.
The complete expression of the proposition follows by a straightforward
calculation.

In case (iii) and (iv) we use \thetag{\vglrevseriesAWpols} as well. In
case (iii) we take $a=-q^{1+\s+\t}$ and in case (iv) we take
$a=-q^{1-\s-\t-2i-2j}$. The resulting $_4\vp_3$-series are not of the
form as in the proposition, but they are transformed to it by use of
Sears's transformation, cf. \cite{\GaspR, (2.10.4)},
$$
{}_4\vp_3\left( {{q^{-n},a,b,c}\atop{d,e,f}};q,q\right) =
a^n {{(e/a,f/a;q)_n}\over{(e,f;q)_n}}
{}_4\vp_3\left( {{q^{-n},a,d/b,d/c}\atop{d,aq^{1-n}/e,aq^{1-n}/f}};
q,q\right),
\tag\eqname{\vglSearsstransform}
$$
where $abc=defq^{n-1}$.
The observation $C^{l,-j}(-\s)= C^{l,j}(\s)$
and a straightforward calculation prove the
correctness of the factor in front of the $_4\vp_3$-series in the
proposition. \qed\enddemo

\demo {Remark} Proposition \propAWnaarqKraw\ can also be proved from
rewriting case (i) of corollary \corexplgenmateltAWpols\ as in the proof
given above and then using the symmetry relations of lemma
\lemmsymmgenmateltas , or \thetag{\vglsymmetriesgenmatelt}, for the other
cases.
However, we still have to use transformations for ${}_4\vp_3$-series to obtain
proposition \propAWnaarqKraw\ in full generality.\enddemo

Two applications of proposition \propAWnaarqKraw\ are given in the rest
of this section. First we derive a formula for dual $q$-Krawtchouk
polynomials, from which we can obtain some special cases.

\proclaim{Proposition \theoremname{\propPoissonkerndualqKraw}}
For $N\in\Zp$, $i,j\in\{ 0,1,\ldots,N\}$, $\s,\t\in\R$ we have
$$
\aligned
&\sum_{n=0}^N q^{{1\over 2}n(\s+\t)}q^{{1\over 2}n(n-1)} {{(q^N;q^{-1})_n}
\over {(q;q)_n}} t^n R_n(q^{-j}-q^{j-N-\s};q^\s,N;q)
R_n(q^{-i}-q^{i-N-\t};q^\t,N;q) \\
&=(-1)^i q^{{1\over 2}i(i-1)} q^{{1\over 2}i(\s-\t)} t^i (t^{-1}q^{{1\over
2}(\s-\t)+1-i};q)_i (tq^{{1\over 2}(\t-\s)-i};q)_j
{{(-tq^{{1\over 2}(\t+\s)};q)_\infty}\over{(-tq^{{1\over
2}(\t+\s)+N-i-j};q)_\infty}} \\
&\qquad\times {}_4\vp_3 \left( {{q^{-i},q^{-j},-q^{\t-i},-q^{j-N-\s}}\atop{
q^{-N},t^{-1}q^{{1\over 2}(\t-\s)+1-i}, tq^{{1\over 2}(\t-\s)-i}}};q,q\right).
\endaligned
$$
\endproclaim

\demo{Proof} Apply the one-dimensional representation $\tl$ to
\thetag{\vglexprgenmateltinstanmatelt} and use proposition \propAWnaarqKraw\
to rewrite the left hand side as a ${}_4\vp_3$-series. Then use the
explicit
value for $v^{l,j}_n(\s)$ of proposition \propspectrumdualqK\ and replace
$q^2$, $l-j$, $l-i$, $l-n$, $2l$, $\l^2$ by $q$, $j$, $i$, $n$, $N$, $t$ to
obtain the proposition. \qed\enddemo

\demo{Remark} Proposition \propPoissonkerndualqKraw\ can also be proved
analytically as follows. Use for the dual $q$-Krawtchouk polynomials the
${}_2\vp_1$-series representation which can be obtained from
\cite{\GaspR , (III.7)}. Interchange the summations and note that the
summation over $n$ from $0$ to $N$ can be summed by the $q$-binomial formula
\cite{\GaspR , (II.4)}. Write the resulting double sum as a sum with a
terminating balanced ${}_3\vp_2$-series in the summand. Using the
$q$-Saalsch\"utz summation formula \cite{\GaspR , (II.12)} we obtain a
${}_4\vp_3$-series, which can be transformed using Sears's transformation
formula \thetag{\vglSearsstransform} to get a ${}_4\vp_3$-series as in the
right hand side of proposition \propPoissonkerndualqKraw .
I thank Mizan Rahman for providing me with this analytic proof of
proposition \propPoissonkerndualqKraw .
\enddemo

Let us note some special cases of proposition \propPoissonkerndualqKraw .
Firstly, take $i=0$, so that one of the dual $q$-Krawtchouk polynomials and
the ${}_4\vp_3$-series reduce to $1$. We obtain the generating function
\thetag{\vglgenfunctiedualqKrawtchouk} for the dual $q$-Krawtchouk
polynomials as a special case of proposition \propPoissonkerndualqKraw .
Secondly, take $\s=\t$, $t=1$ and note that
$$
\aligned
&(q^{1-i};q)_i (q^{-i};q)_j
\,{}_4\vp_3 \left( {{q^{-i},q^{-j},-q^{\s-i},-q^{j-N-\s}}\atop{
q^{-N},q^{1-i}, q^{-i}}};q,q\right) \\
= & \sum_{p=0}^{\min(i,j)}
{{(q^{-i},q^{-j},-q^{\s-i},-q^{j-N-\s};q)_p}\over{
(q^{-N},q;q)_p}} q^p (q^{1-i+p};q)_{i-p} (q^{p-i};q)_{j-p} \\
= &\ \d_{i,j}
{{(q^{-i},q^{-i},-q^{\s-i},-q^{i-N-\s};q)_i}\over{
(q^{-N},q;q)_i}} q^i ,
\endaligned
$$
since $(q^{1-i+p};q)_{i-p}=\d_{i,p}$ and $(q^{p-i};q)_{j-p}=0$ for $j>i$. We
obtain the orthogonality relations for the $q$-Krawtchouk polynomials, cf.
\thetag{\vglorthorelqKrawtchouk}, as a special case of proposition
\propPoissonkerndualqKraw . Finally, note that if we only
specialise $\s=\t$ in proposition \propPoissonkerndualqKraw , we obtain a
Poisson kernel for the dual $q$-Krawtchouk polynomials.
The Poisson kernel is a special case of the Gasper-Rahman
formula for the Poisson kernel for the $q$-Racah polynomials, cf.
\cite{\GaspRPK , (3.11)}, \cite{\GaspR , (8.7.13),
but introduce a factor $(q;q)_s$ in the denominator}.
This calculation is far from trivial and I thank George Gasper and especially
Mizan Rahman for making this identification precise.
Note that for $\s\not=\t$ proposition \propPoissonkerndualqKraw\ does
not fit into the Gasper-Rahman Poisson kernel.

As another application of proposition \propAWnaarqKraw\ we rewrite the
orthogonality relations of corollary \corrHopfpropgenmatelt (iii) for
$\s=\t$ in terms of the ${}_4\vp_3$-series derived in proposition
\propAWnaarqKraw . Unfortunately, the orthogonality relations we find
are not the orthogonality relations for the $q$-Racah polynomials
\cite{\AskeWtwee}, but orthogonality relations for the function
$$
\aligned
r_n(q^p;\s,\t,e^{i\theta},N;q) = &(e^{i\theta}q^{{1\over 2}(\t-\s)-n};q)_p
(-e^{i\theta}q^{{1\over 2}(\t+\s)-n};q)_{N-p} \\
&\times \, {}_4\vp_3 \left( {{q^{-n},q^{-p},-q^{\t-n},-q^{p-N-\s}}\atop{
q^{-N},q^{1+{1\over 2}(\t-\s)-n}e^{-i\theta}, q^{{1\over
2}(\t-\s)-n}e^{i\theta}}};q,q\right)
\endaligned
\tag\eqname{\vgldefstrangeqKrawtchoukextension}
$$
for $N\in\Zp$, $n,p\in\{ 0,1,\ldots,N\}$, $\s,\t\in\R$. Use $(a;q)_p =
(a;q)_\infty/(aq^p;q)_\infty$ to see that $r_n$ is indeed a function of
$q^p$. The dual $q$-Krawtchouk polynomial \thetag{\vgldefdualqKrawtchoukpol}
can be obtained as a special case by letting $\t\to\infty$ in
\thetag{\vgldefstrangeqKrawtchoukextension}.

As already stated these orthogonality relations do not correspond to the
orthogonality relations for the $q$-Racah polynomials. It is known that in
the special case $\s=\t=\infty$ the discrete orthogonality relations for the
so-called quantum $q$-Krawtchouk polynomials can be obtained in this way,
cf. Koornwinder \cite{\Koordrie , \S 6}. This result is extended to the
special case $\s=\infty$ by Noumi and Mimachi, cf.
\cite{\NoumMCompM , \S 4}, from which the discrete orthogonality
relations of the (dual) $q$-Hahn polynomials,
cf. \ssectqHahn , can be obtained.
These derivations take place on the non-commutative level,
in the algebra $\A$ in case $\s=\t=\infty$ and in an extension of the
algebra $\A$ in case $\s=\infty$.
It might be possible that the orthogonality relations for the
$q$-Racah polynomials can be obtained in a still larger extension of $\A$
in a way similar as in Noumi and Mimachi \cite{\NoumMCompM}. I thank Masatoshi
Noumi for explaining the ideas and backgrounds of \cite{\NoumMCompM} to me.

\proclaim{Proposition \theoremname{\proporthorelforextqKrawtchouk}}
The function $r_n(q^p)= r_n(q^p;\s,\t,e^{i\theta},N;q)$ satisfies the
following discrete orthogonality relations;
$$
\sum_{p=0}^N r_n(q^p) \overline{r_m(q^p)}
{{(1+q^{2p-N-\s})}\over{(-q^{p-2N-\s})^p}}
{{(q^{-N},-q^{-N-\s};q)_p}\over{(q,-q^{1-\s};q)_p}} = \d_{n,m} h_n^{-1}
\tag\eqname{\vglorthorelforextqKrawtchouk}
$$
where
$$
\aligned
h_n = &(-q^{2N-\s+1})^n {{ q^{2N+\s+\t}(1+q^{2n-N-\t})
(q^{-N},-q^{-N-\t};q)_n}\over{ (-q^\t,-q^\s;q)_{N+1}
(q,-q^{1-\t};q)_n }} \\
&\times {{ (q^{{1\over 2}(\s-\t)}e^{i\theta},
q^{{1\over 2}(\s-\t)}e^{-i\theta};q)_n}\over{
(-q^{1-{1\over 2}(\t+\s)}e^{i\theta},
-q^{1-{1\over 2}(\t+\s)}e^{-i\theta};q)_n }}.
\endaligned
$$
\endproclaim

\demo{Remark} Note that the weights in
\thetag{\vglorthorelforextqKrawtchouk} are equal to weights in the
orthogonality relations for the dual $q$-Krawtchouk polynomials, cf.
\thetag{\vglorthoreldualqKrawtchouk}.  Moreover,
\thetag{\vglorthorelforextqKrawtchouk} reduces to the orthogonality
relations for the dual $q$-Krawtchouk polynomials if $\t\to\infty$.
\enddemo

\demo{Proof} Using the notation
\thetag{\vgldefstrangeqKrawtchoukextension} we can rewrite proposition
\propAWnaarqKraw\ as
$$
\align
\hp\bigl( a^l_{i,p}(\t,\s)\bigr) =&\ (-1)^{l-i} C^{l,i}(\t) q^{(l-i)(l-i-1)}
q^{(l-i)(\s-\t)} e^{-ii\theta}
{{(e^{-i\theta}q^{\t-\s+2-2l+2i};q^2)_{l-i}}\over
{(-e^{i\theta}q^{\t+\s-2l+2i};q^2)_{l-i}}} \\
&\times C^{l,p}(\s)\ r_{l-i}(q^{2(l-p)};\s,\t,e^{i\theta},2l;q^2),
\tag\eqname{\vglhpalipasgeneralisedqKrawtcouk}
\endalign
$$
where we used $\hp=\t_{e^{i\theta/2}}$, cf. \ssectreprsA .
From corollary \corrHopfpropgenmatelt (iii)
we obtain the orthogonality relations
$$
\sum_{p=-l}^l \hp\bigl( a^l_{i,p}(\t,\s)\bigr) \overline{
\hp\bigl( a^l_{j,p}(\t,\s)\bigr)} = \d_{i,j},
$$
since $\hp$ is a $\ast$-representation of $\A$.
In this expression we use
\thetag{\vglhpalipasgeneralisedqKrawtcouk} and
$$
\vert C^{l,p}(\s)\vert^2 = {{ q^{4l+2\s} (1+q^{-4p-2\s})}\over
{(-q^{2(l-p)-8l-2\s})^{l-p} (-q^{2\s};q^2)_{2l+1} }}\, ,
$$
cf. proposition \propspectrumdualqK .
Now replace $l-i$, $l-j$, $l-p$ $2l$, $q^2$ by $n$, $m$, $p$, $N$, $q$ to
obtain the statement of the proposition. \qed \enddemo

The orthogonality relations dual to
\thetag{\vglorthorelforextqKrawtchouk}, which are obtainable from
$$
\sum_{p=-l}^l \hp\bigl( a^l_{p,j}(\t,\s)\bigr) \overline{
\hp\bigl( a^l_{p,i}(\t,\s)\bigr)} = \d_{i,j},
$$
cf. corollary \corrHopfpropgenmatelt (iii), and proposition
\propAWnaarqKraw , can be rewritten as the orthogonality relations
\thetag{\vglorthorelforextqKrawtchouk} with $\s$ and $\t$ interchanged
by lemma \lemmsymmgenmateltas , cf. also \thetag{\vglsymmetriesgenmatelt}.

%&&&&&&&&&&&&&&&&&&&&&&&&&&&&&&&&&&&&
%%N E W   S E C T I O N%%%%%%%%%%%%**
%%%%%%%%%%%%%%%%%%%%%%%%%%%%%%%%%%%%%

\head\newsection . Spherical and associated spherical elements
\endhead

In this section we show how the associated spherical elements can be
obtained from the spherical elements. In case of the continuous
$q$-ultraspherical polynomials we explicitly give an operator mapping
the continuous $q$-Legendre polynomial to the continuous
$q$-ultraspherical polynomials. This operator is a $q$-Hahn polynomial
in a very simple shift operator. This section is motivated by the paper
\cite{\BadeK} by Badertscher and Koornwinder.

Consider the irreducible unitary matrix corepresentation
$a^l(\s)=\bigl( a^l_{i,j}(\s,\s)\bigr)_{i,j=-l,\ldots,l}$, cf.
corollary \corrHopfpropgenmatelt (i). Then we have
$$
\langle t^l(D^2)v^{l,j}(\s),v^{l,i}(\s)\rangle =
\langle D^2, a^l_{i,j}(\s,\s)\rangle = \t_{q^{-1}}\bigl(
a^l_{i,j}(\s,\s)\bigr),
$$
since the algebra homomorphisms $\t_{q^{k/2}}\colon\A\to\C$ given in
\ssectreprsA , and $\xi\mapsto\langle A^k,\xi\rangle$, $k\in\Z$, coincide
on the generators of $\A$. From corollary \corexplgenmateltAWpols\ we get
$$
\langle t^l(D^2)v^{l,j}(\s),v^{l,i}(\s)\rangle = C
(q^{-2};q^2)_{\vert i-j\vert}
$$
for all $i,j=-l,-l+1,\ldots,l$. So $t^l(D^2)$ is a tridiagonal matrix with
respect to the orthonormal basis $v^{l,i}(\s)$, $i=-l,-l+1,\ldots,l$. The
subdiagonal and superdiagonal are non-zero for $\s\not= \pm\infty$.

In the remainder of this section we assume $l\in\Zp$ fixed. Let
$A(\s)=\bigl( A_{i,j}(\s)\bigr)$ be the matrix of $t^l(D^2)$ with respect to
the basis $\{ v^{l,i}(\s)\mid i=-l,-l+1,\ldots, l\}$.

\proclaim{Proposition \theoremname{\proptridiagonalmatrix}}
$A(\s)$ is a tridiagonal, real, symmetric matrix with positive diagonal
entries. Moreover, $A_{-i,-j}(-\s)=A_{i,j}(\s)$.
\endproclaim

\demo{Proof} The operator $t^l(D^2)$ is self-adjoint, and
$$
A_{i,j}(\s) = \langle t^l(D^2)v^{l,j}(\s),v^{l,i}(\s)\rangle =
\sum_{m=-l}^l v^{l,j}_m(\s) \overline{v^{l,i}_m(\s)}q^{2m}.
$$
The explicit value of the coefficients $v^{l,j}_n(\s)$ given in proposition
\propspectrumdualqK\ shows that this yields a real value, which is obviously
positive for $i=j$. The last statement follows from
\thetag{\vglsymmvljsigmaenminsigma}. \qed\enddemo

If we denote the $k$-th power of $A(\s)$ by $A^k(\s)$, then we get
$A^k_{i,j}(\s)=A^k_{-i,-j}(-\s)$ by induction with respect to $k$.
Indeed, the cases $k=0,1$ being known, we get
$$
A^{k+1}_{i,j}(\s) =\sum_{p=-l}^l A_{i,p}(\s)A^k_{p,j}(\s) =
\sum_{p=-l}^l A_{-i,-p}(-\s)A^k_{-p,-j}(-\s) =
A^{k+1}_{-i,-j}(-\s).
$$
Consequently, for $k=0,1,\ldots,l$ and with $A_k(\s)=A^k_{k,0}(\s)$
$$
t^l(D^{2k})v^{l,0}(\s) = A_k(\s)v^{l,k}(\s)+A_k(-\s)v^{l,-k}(\s)
+\sum_{p=-k+1}^{k-1} c_p(\s) v^{l,p}(\s)
\tag\eqname{\vglcaltlDtweekopsphervect}
$$
for certain constants $c_p(\s)=A^k_{p,0}(\s)$
satisfying $c_p(\s)=c_{-p}(-\s)$.

So the (ordered) set of vectors
$$
v^{l,0}(\s),\quad t^l(D^2)v^{l,0}(\s),\quad t^l(D^4)v^{l,0}(\s),\quad
\ldots,\quad t^l(D^{2l})v^{l,0}(\s)
$$
is linearly independent. Applying the Gram-Schmidt orthogonalisation process
to these vectors yields
$$
v^{l,0}(\s), \quad A_1(\s)v^{l,1}(\s)+A_1(-\s)v^{l,-1}(\s),\quad
\ldots,\quad A_l(\s)v^{l,l}(\s)+A_l(-\s)v^{l,-l}(\s),
$$
as a set of orthogonal vectors because of the specific form of
\thetag{\vglcaltlDtweekopsphervect}.

\proclaim{Proposition \theoremname{\propdiscreteorthpols}}
There exists a set of monic orthogonal polynomials $\{ p_k\}_{k=0}^l$ of
degree $k$ such that
$$
t^l\bigl( p_k(D^2)\bigr) v^{l,0}(\s) =
A_k(\s)v^{l,k}(\s)+A_k(-\s)v^{l,-k}(\s)
\tag\eqname{\vgldefdiscreteortpols}
$$
for $k=0,1,\ldots,l$.
\endproclaim

\demo{Proof} The Gram-Schmidt orthogonalisation process yields the $k$-th
orthogonal vector as a linear combination of the zeroth until the $k$-th
vector of the linearly independent set to be orthogonalised. This shows the
existence of the polynomial $p_k$ in \thetag{\vgldefdiscreteortpols}, which
has to have degree $k$ and leading coefficient $1$
because of \thetag{\vglcaltlDtweekopsphervect}.

Consider the moment functional $\L$ defined by
$$
\L (x^k) = \langle t^l(D^{2k})v^{l,0}(\s),v^{l,0}(\s)\rangle,
\tag\eqname{\vgldefmomentsfunctional}
$$
then, for $0\leq n,m\leq l$,
$$
\aligned
\L \bigl(\bar p_n(x)p_m(x)\bigr) &=
\langle t^l\bigl(p_m(D^2)\bigr) v^{l,0}(\s),
t^l\bigl(p_n(D^2)\bigr) v^{l,0}(\s)\rangle \\
&= \delta_{n,m}
\cases A_n^2(\s) + A_n^2(-\s), &\text{$n>0$,}\\
1, &\text{$n=0$,}\endcases
\endaligned
$$
since $D^\ast=D$. \qed\enddemo

More explicit formulas for the moments of the moment functional $\L$ are given
by
$$
\L (x^k)  = \sum_{n=-l}^l \vert v^{l,0}_n(\s)\vert^2 q^{2kn}
 = {}_4\vp_3 \left( {{q^{-2l}, q^{2+2l}, q^{-2k}, q^{2+2k}}\atop
{q^2,-q^{2-2\s},-q^{2+2\s}}};q^2,q^2\right).
$$
The first equality follows form proposition \propspectrumdualqK ,
\thetag{\vgldefreprU} and
\thetag{\vgldefmomentsfunctional}, and the second equality follows from the
observation that \thetag{\vgldefmomentsfunctional} equals
$\t_{q^{-k}}\bigl( a^l_{0,0}(\s,\s)\bigr)$ and corollary
\corexplgenmateltAWpols .
So the orthogonality measure corresponding to the moment functional $\L$ is
the finite discrete measure with positive weights $\vert v^{l,0}_n(\s)\vert^2$
at the points $x_n=q^{2n}$, $n=-l,-l+1,\ldots,l$.

Unfortunately, I have not (yet) been able to obtain an explicit formula for
the orthogonal polynomials $\{ p_k\}_{k=0}^l$ unless $\s=0$. In this case we
can obtain special $q$-Hahn polynomials by identifying the orthogonality
measures. We start with the lemma in which we do the necessary computation,
cf. \cite{\Koortwee , p.~805}.

\proclaim{Lemma \theoremname{\lemmacompincasesigmaisnul}}
In case $\s=0$ the weights for $l-n$ even are given by
$$
\vert v^{l,0}_n(0)\vert^2 = q^{l+n} {{ (q^2;q^4)_{(l-n)/2}(q^2;q^4)_{(l+n)/2}}
\over{(q^4;q^4)_{(l-n)/2}(q^4;q^4)_{(l+n)/2}}}
$$
and $\vert v^{l,0}_n(0)\vert^2=0$ for $l-n$ odd.
\endproclaim

\demo{Proof} The proof is based on Andrews's summation formula, cf.
\cite{\GaspR, (II.17)},
$$
{}_4\vp_3 \left( {{q^{-n}, aq^n,c,-c}\atop{(aq)^{1/2},-(aq)^{1/2},c^2}};
q,q\right) = \cases 0, &\text{if $n$ is odd,} \\
{\displaystyle{ {{c^n(q,aq/c^2;q^2)_{n/2}}\over{(aq,c^2q;q^2)_{n/2}}}, }}
&\text{if $n$ is even.}\endcases
$$
An application of this formula with $a$, $c$, $n$ and $q$ replaced by
$0$, $q^{-2l}$, $l-n$ and $q^2$ shows that, cf.
\thetag{\vgldefdualqKrawtchoukpol},
$$
R_{l-n}(q^{-2l}-q^{-2l};1,2l;q^2)
= \cases 0, &\text{if $l-n$ is odd,} \\
{\displaystyle{ {{q^{-2l(l-n)}(q^2;q^4)_{(l-n)/2}}\over
{(q^{2-4l};q^4)_{(l-n)/2}}}, }} &\text{if $l-n$ is even.}\endcases
$$
This proves that $\vert v^{l,0}_n(0)\vert^2=0$ for $l-n$ odd, so it remains to
calculate the value for $l-n$ even.

It follows from proposition \propspectrumdualqK\ that we can rewrite
$\vert C^{l,0}(0)\vert^2 = q^{2l} (q^2;q^4)_l/(q^4;q^4)_l$.
Hence, proposition \propspectrumdualqK\ and the evaluation of the dual
$q$-Krawtchouk polynomial in this case imply that
$$
\aligned
\vert v^{l,0}_n(0)\vert^2= & q^{2l} {{(q^2;q^4)_l}\over{(q^4;q^4)_l}}
q^{(l-n)(l-n-1)}{{(q^{4l};q^{-2})_{l-n}}\over{(q^2;q^2)_{l-n}}} \\
& \times q^{-4l(l-n)} q^{(l-n)(4l-2)} q^{-(l-n)(l-n-2)}
\Bigl( {{(q^2;q^4)_{(l-n)/2}}\over{(q^{4l-2};q^{-4})_{(l-n)/2}}}
\Bigr)^2
\endaligned
$$
by inverting one of the $q$-shifted factorials. Now we use
$$
\gather
{{(q^2;q^4)_{(l-n)/2}}\over{(q^2;q^2)_{l-n}}} = {1\over{(q^4;q^4)_{(l-n)/2}}},
\quad
{{(q^{4l};q^{-2})_{l-n}}\over{(q^{4l-2};q^{-4})_{(l-n)/2}}} =
(q^{4l};q^{-4})_{(l-n)/2}, \\
{{(q^{4l};q^{-4})_{(l-n)/2}}\over{(q^4;q^4)_l}} =
{1\over{(q^4;q^4)_{(l+n)/2}}}, \quad
{{(q^2;q^4)_l}\over{(q^{4l-2};q^{-4})_{(l-n)/2}}} =
(q^2;q^4)_{(l+n)/2}.
\endgather
$$
Collecting the powers of $q$ proves the lemma. \qed\enddemo

Using the moment functional ${\Cal L}^l_0$ introduced in the proof of
lemma \lemmacompincasesigmaisnul\ we see that in case $\s=0$ we have to
find orthogonal polynomials $\{ p_n\}_{n=0}^l$ satisfying
$$
\gather
\sum_{n=-l, \ l-n\in 2\Z}^l
q^{l+n} {{ (q^2;q^4)_{(l-n)/2}(q^2;q^4)_{(l+n)/2}}
\over{(q^4;q^4)_{(l-n)/2}(q^4;q^4)_{(l+n)/2}}} \bar p_k(q^{2n})
p_m(q^{2n}) = 0, \quad k\not= m, \\
\Longleftrightarrow
\sum_{n=0}^l q^{2l-2n} {{ (q^2;q^4)_n(q^2;q^4)_{l-n}}
\over{(q^4;q^4)_n(q^4;q^4)_{l-n}}} \bar p_k(q^{2l-4n}) p_m(q^{2l-4n}) = 0,
\quad k\not= m.
\tag\eqname{\vglorthdiscrexpl}
\endgather
$$
The weights in \thetag{\vglorthdiscrexpl}
almost coincide with the Fourier coefficients of the continuous
$q$-Legendre polynomial, cf. \ssectcontqultraspherical ,
\cite{\AskeI, (3.1)}, \cite{\GaspR , (7.4.2)} with $\b=q^{1/2}$,
which is analogous to the
general classical case described in \cite{\BadeK, prop. 5.2}.
Comparing \thetag{\vglorthdiscrexpl} with the orthogonality relations for the
$q$-Hahn polynomials, cf. \thetag{\vglorthorelqHahn}, leads to
$$
p_m(x) = q^{2lm} {{(q^2;q^4)_m(q^{-4l};q^4)_m}\over{(q^{4m};q^4)_m}}
Q_m(q^{-2l}x;q^{-2}, q^{-2},l;q^4), \qquad m=0,1,\ldots
l,
\tag\eqname{\vglintpqHahnpols}
$$
where the constant follows from the condition that $p_m$ is monic. So we
have for the polynomial $p_m$ defined in \thetag{\vglintpqHahnpols}
$$
t^l\bigl( p_m(D^2)\bigr) v^{l,0}(0) =A_m(0) \bigl(
v^{l,m}(0) + v^{l,-m}(0)\bigr)
$$
and consequently, using the notation of \S\S 4, 6
$$
\aligned
\langle X,&p_m(D^2).a^l_{i,0}(\t,0)\rangle =
\langle t^l(X)t^l\bigl(p_m(D^2)\bigr)v^{l,0}(0), v^{l,i}(\t)\rangle \\
&=A_m(0)\bigl(\langle t^l(X)v^{l,m}(0), v^{l,i}(\t)\rangle +
\langle t^l(X)v^{l,-m}(0), v^{l,i}(\t)\rangle\bigr) \\
&=A_m(0)\bigl(\langle X, a^l_{i,m}(\t,0)\rangle+\langle X,
a^l_{i,-m}(\t,0)\rangle\bigr),
\endaligned
$$
so that as an identity in $\A$ we have
$$
p_m(D^2).a^l_{i,0}(\t,0) = A_m(0)\bigl(
a^l_{i,m}(\t,0)+a^l_{i,-m}(\t,0)\bigr),
\tag\eqname{\vglstartintaflultraq}
$$
with $p_m$ defined in terms of a $q$-Hahn polynomial as in
\thetag{\vglintpqHahnpols}.

\proclaim{Proposition \theoremname{\propquantumqultra}}
For $m,l\in\Zp$, $m\leq l$ we have
$$
\aligned
&Q_m(q^{-l}E_-;q^{-1},q^{-1},l;q^2)
p_l^{(0,0)}(\cos\theta;1,q^{{1\over 2}\t}\mid q) =\\
&{1\over 2} (q^{l+1};q)_m e^{-im\theta} \Bigl\lbrace
(e^{i\theta}q^{{1\over 2}(1-\t)},-e^{i\theta}q^{{1\over 2}(1+\t)};q)_m
p_{l-m}^{(m,m)}(\cos\theta;1,q^{{1\over 2}\t}\mid q) \\
&\qquad\qquad +
(e^{i\theta}q^{{1\over 2}(1+\t)},-e^{i\theta}q^{{1\over 2}(1-\t)};q)_m
p_{l-m}^{(m,m)}(\cos\theta;1,q^{-{1\over 2}\t}\mid q) \Bigr\rbrace ,
\endaligned
\tag\eqname{\vglqHahndifopongenqLegendre}
$$
where $Q_m$ denotes a $q$-Hahn polynomial
\thetag{\vgldefqHahnpol}, $p_{l-m}^{(m,m)}$ an Askey-Wilson polynomial
\thetag{\vgldefAWalsqJacobi} and $E_-$ is the operator acting on
trigonometric polynomials $f(e^{i\theta})$ by $\bigl( E_-f\bigr)
(e^{i\theta}) = f(q^{-1}e^{i\theta})$.

In particular, for $\t=0$ we have
$$
\aligned
&Q_m(q^{-l}E_-;q^{-1},q^{-1},l;q^2) C_l(\cos\theta;q\mid q^2)\\
= & {{(q;q^2)_m}\over{(q^{2l-2m+2};q^2)_m}} e^{-im\theta}
(qe^{2i\theta};q^2)_m C_{l-m}(\cos\theta;q^{1+2m}\mid q^2),
\endaligned
\tag\eqname{\vglidqHahncontqultra}
$$
where $C_l(\cos\theta;\b\mid q)$
denotes a continuous $q$-ultraspherical polynomial
\thetag{\vgldefcontqultraspherpols}.
\endproclaim

\demo{Remark} A generalisation of \thetag{\vglidqHahncontqultra} to
arbitrary continuous $q$-ultraspherical polynomials using a generating
function for the dual $q$-Hahn polynomials
is given in \cite{\KoelUP , prop.~1.1}. \enddemo

\demo{Proof} The proof consists of applying the one-dimensional
representation $\tl$ to \thetag{\vglstartintaflultraq} for $i=0$.
In the proof of corollary \corexplgenmateltAWpols\ we already noted
that $\tl(D.\xi)=\t_{\l q^{-1/2}}(\xi)$, $\xi\in\A$. More generally we
have $\tl(D^{2j}.\xi)=\t_{\l q^{-j}}(\xi)$, $\xi\in\A$. Introduce the
operator $F_-$ acting on Laurent polynomials in $\l$ by
$\bigl( F_-f\bigr)(\l) = f(q^{-1}\l)$, then we get
$\tl\bigl(p(D^2).\xi\bigr) = p(F_-) \tl(\xi)$ for arbitrary
polynomials $p$.

From this observation, corollary \corexplgenmateltAWpols\
and \thetag{\vglintpqHahnpols}
we get for some constant $C$ the expression
$$
\align
&Q_m(q^{-2l}F_-;q^{-2},q^{-2},l;q^4)
p_l^{(0,0)}((q^{-1}\l^2+q\l^{-2})/2;1,q^\t\mid q^2) = \\
&C \l^{-2m} \Bigl\lbrace
(\l^2q^{-\t},-\l^2q^\t;q^2)_m
p_{l-m}^{(m,m)}((q^{-1}\l^2+q\l^{-2})/2;1,q^\t\mid q^2) \\
&\qquad\qquad +
(\l^2q^\t,-\l^2q^{-\t};q^2)_m
p_{l-m}^{(m,m)}((q^{-1}\l^2+q\l^{-2})/2;1,q^{-\t}\mid q^2) \Bigr\rbrace
\endalign
$$
after applying $\t_\l$ to \thetag{\vglstartintaflultraq} for $i=0$.
In this expression we replace $q^{-1}\l^2$, $q^2$ by $e^{i\theta}$, $q$,
so that the operator $F_-$ goes over into $E_-$. We obtain the
expression of the proposition apart from the explicit expression for the
constant involved. The constant is calculated by comparing the
coefficients of $e^{-il\theta}$ on both sides of
\thetag{\vglqHahndifopongenqLegendre}. Since the coefficient of
$e^{-i(l-m)\theta}$ in $p_{l-m}^{(m,m)}(\cos\theta;1,q^{\pm{1\over
2}\t}\mid q)$ is $(q^{2l};q^{-1})_{l-m}$, cf. \cite{\AskeW , p.5}, and
the $q$-Hahn polynomial of argument $1$ equals $1$ we have to solve
$(q^{2l};q^{-1})_l=2C(q^{2l};q^{-1})_{l-m}$. This proves
\thetag{\vglqHahndifopongenqLegendre}.

To prove the special case $\t=0$ we note that both terms on the right
hand side of \thetag{\vglqHahndifopongenqLegendre} are equal and that we
can combine the $q$-shifted factorials involving $e^{i\theta}$.
Because of \thetag{\vglquadraticcqusAWp}
the last statement of the proposition
follows after a short manipulation of $q$-shifted factorials.
\qed\enddemo

%&&&&&&&&&&&&&&&&&&&&&&&&&&&&&&&&&&&&
%%N E W   S E C T I O N%%%%%%%%%%%%**
%%%%%%%%%%%%%%%%%%%%%%%%%%%%%%%%%%%%%

\head\newsection . Characters \endhead

The characters of the irreducible unitary representations $t^l$ of the
quantum $SU(2)$ group are investigated in this section. We obtain
expansions for the Chebyshev polynomial of the second kind in terms of
Askey-Wilson polynomials.

For any finite dimensional
representation $t$ of the algebra $\U$ we define the character $\chi$
of the representation as a linear functional on $\U$ by
$\langle X, \chi\rangle = tr(t(X))$. Here $tr$ denotes the trace
operator. In particular for the irreducible $\ast$-representations $t^l$
of $\U$, cf. \ssectreprthey , of spin $l\in\hZp$, we define the
character $\chi_l$ by
$\langle X,\chi_l\rangle = tr\bigl( t^l(X)\bigr)$. Since the matrix
elements of the irreducible $\ast$-representation $t^l$ are elements of
$\A$ we see that $\chi_l\in\A$ and
$$
\chi_l = \sum_{n=-l}^l t^l_{n,n} = \sum_{n=-l}^l a^l_{n,n}(\s,\s),
\tag\eqname{\vgldefcharacterchil}
$$
where the last equality follows from corollary \corrHopfpropgenmatelt
(i) and \ssectreprthey , or from the fact that the trace operator is
independent of the basis.

Since $t^l$ is a representation we find for $X,Y\in\U$
$$
\langle XY, \chi_l\rangle = tr\bigl( t^l(X)t^l(Y)\bigr)
= tr\bigl( t^l(Y)t^l(X)\bigr) = \langle YX, \chi_l\rangle
$$
and from \thetag{\vgldualcomult} we then find that the character
$\chi_l$, and, for that matter, any character of a finite dimensional
representation, is cocentral.
This means $\Delta(\chi_l) = \omega \circ \Delta(\chi_l)$ where
$\omega \colon \eta\otimes \xi \mapsto \xi\otimes\eta$ denotes the flip
automorphism of $\A \otimes \A$ and $\Delta$ is the comultiplication of
$\A$, cf. \ssectquanfunalg .

The characters of the irreducible unitary representations and the
algebra of cocentral elements have been
investigated by Woronowicz \cite{\Woro , \S 5, App.~A.1}. We collect
some of Woronowicz's results in the following proposition.

\proclaim{Proposition \theoremname{\propWoroonChars}}
(i) The algebra of cocentral elements is generated by $(\a+\d)/2$. \par
\noindent
(ii) The Haar functional on the algebra of
cocentral elements is given by
$$
h\biggl( p\bigl( {1\over 2}(\a+\d)\bigr)\biggr) = {2\over \pi}
\int_{-1}^1 p(x) \sqrt{1-x^2}\ dx
$$
for any polynomial $p$. \par
\noindent
(iii) For $l\in\hZp$ we have $\chi_l= U_{2l}\bigl( (\a+\d)/2\bigr)$,
where $U_n(\cos\theta) =
\sin (n+1)\theta/\sin\theta$ is the Chebyshev polynomial of the second
kind.
\endproclaim

By \thetag{\vgldefcharacterchil} and proposition \propWoroonChars (iii)
we get
$$
U_{2l}\bigl( (\a+\d)/2\bigr) = \sum_{n=-l}^l a^l_{n,n}(\s,\s)
$$
for arbitrary $\s\in\R$. Next we apply the one-dimensional
representation $\t_{q^{1/2}e^{i\theta/2}}$, cf. \ssectreprsA , to
this identity. We obtain from corollary \corexplgenmateltAWpols\
the identity
$$
\aligned
&\sum_{m=0}^{\lceil l-1\rceil} {{\vert C^{l,l-m}(\s)\vert^2
q^{-l}}\over{(q^{4l};q^{-2})_m}} e^{-i(l-m)\theta}
(-e^{i\theta}q^{1+2\s};q^2)_{2l-2m}
p_m^{(0,2l-2m)} (\cos\theta;q^\s,q^\s\mid q^2)      \\
+&\sum_{m=0}^{\lceil l-1\rceil} {{\vert C^{l,l-m}(-\s)\vert^2
q^{-l}}\over{(q^{4l};q^{-2})_m}} e^{-i(l-m)\theta}
(-e^{i\theta}q^{1-2\s};q^2)_{2l-2m}
p_m^{(0,2l-2m)} (\cos\theta;q^{-\s},q^{-\s}\mid q^2)      \\
+'&\, {{\vert C^{l,0}(\s)\vert^2
q^{-l}}\over{(q^{2l};q^{-2})_l}}
p_l^{(0,0)}(\cos\theta;q^\s,q^\s\mid q^2) =
U_{2l}((q^{1/2}e^{i\theta/2}+q^{-1/2}e^{-i\theta/2})/2),
\endaligned
\tag\eqname{\vglexpanforUl}
$$
where we have replaced $l-n$ by $m$.
The $+'$ indicates that the last term in \thetag{\vglexpanforUl}
is only occurring for $l\in\Zp$, $\lceil a\rceil$ denotes the smallest
integer greater than or equal to $a\in\R$ and $C^{l,m}(\s)$ is given in
proposition \propspectrumdualqK .
Note that the right hand side of
\thetag{\vglexpanforUl} is independent of $\s$. Askey and
Wilson \cite{\AskeW , (2.18)} already noted that the Chebyshev
polynomials of the second kind can be obtained as a special case of the
Askey-Wilson polynomials. Equation \thetag{\vglexpanforUl} is
a $q$-analogue of Vilenkin \cite{\Vile , Ch.~3, \S 7.1(10)}, Vilenkin
and Klimyk \cite{\VileK , Vol.~1, \S 6.9.2(6)}.

Now that we have represented \thetag{\vgldefcharacterchil} for $\s\in\R$
as \thetag{\vglexpanforUl} we can proceed to find an identity for
special functions corresponding to the case $\s\to\infty$ of
\thetag{\vgldefcharacterchil}, i.e. the first equality of
\thetag{\vgldefcharacterchil}. Applying the one-dimensional
$\ast$-representation $\p$ gives the known equality
$U_{2l}(\cos\theta) = \sin (2l+1)\theta/\sin\theta$.

To obtain another identity we apply the
infinite dimensional $\ast$-representation $\pi^\infty=\pi^\infty_0$,
cf. \ssectreprsA , to the first equality in
\thetag{\vgldefcharacterchil}. So we obtain an identity for bounded
linear operators acting on $\ell^2(\Zp)$. To rewrite this as an identity
for functions we try to find eigenvectors for $\pi^\infty({1\over
2}(\a+\d))$. Note that the matrix elements $t^l_{n,m}$ already act
nicely on the standard orthonormal basis $\{ f_n\}_{n\in\Zp}$ of
$\ell^2(\Zp)$ by \thetag{\vglmatrixeltalslittleqJacobi}, cf. \ssectreprsA .

Now $\sum_{n=0}^\infty c_nf_n$ is an eigenvector for the
eigenvalue $\l$ of the self-adjoint operator $\pi^\infty({1\over
2}(\a+\d))$ if and only if
$$
c_{n+1} \sqrt{1-q^{2n+2}} + c_{n-1} \sqrt{1-q^{2n}} = 2\l c_n, \qquad
n\in\Zp
$$
with the convention $c_{-1}=0$, $c_0=1$. So we can consider $c_n$ as a
polynomial of degree $n$ in $\l$. A suitable normalisation,
$c_n=(q^2;q^2)_n^{-1/2}p_n(\l)$, leads to the three-term recurrence
relation
$$
p_{n+1}(\l) + (1-q^{2n})p_{n-1}(\l) = 2\l p_n(\l), \qquad n\in\Zp,
$$
with initial conditions $p_{-1}(\l)=0$, $p_0(\l)=1$.
The solution $c_n= (q^2;q^2)_n^{-1/2} H_n(\l\mid q^2)$ follows
by comparing this with
the three-term recurrence relation for the continuous $q$-Hermite
polynomials, cf. \thetag{\vgldrietermcontqHermite}.

However, the corresponding eigenvector $\sum_{n=0}^\infty (q^2;q^2)_n^{-1/2}
H_n(\l\mid q^2)f_n$ is not an element of the Hilbert space $\ell^2(\Zp)$,
as can be seen from the asymptotic expansion for the Askey-Wilson
polynomials, cf. Ismail and Wilson \cite{\IsmaW , (1.11)-(1.13)}, or
from the $q$-analogue of Mehler's formula for the continuous $q$-Hermite
polynomials, cf. Askey and Ismail \cite{\AskeI , (6.5)}.
To circumvent this problem we introduce the
vector $v_\l(N) = \sum_{n=0}^N (q^2;q^2)_n^{-1/2} H_n(\l\mid q^2)f_n$
for which we have
$$
\pi^\infty \bigl({1\over 2}(\a+\d)\bigr) v_\l(N) = \l v_\l(N) +
{{H_N(\l\mid q^2)}\over{2\sqrt{(q^2;q^2)_N} }}
\sqrt{1-q^{2N+2}} \, f_{N+1}
-{{H_{N+1}(\l\mid q^2)}\over{2\sqrt{(q^2;q^2)_N} }} \, f_N
$$
and
$$
\pi^\infty\biggl( U_{2l} \bigl({1\over 2}(\a+\d)\bigr) \biggr)
v_\l(N) = U_{2l}(\l) v_\l(N) +
\sum_{p=N-2l+1}^{N+2l} c_p\, f_p
\tag\eqname{\vglchoppedoffeigenvector}
$$
for certain constants $c_p$, $N>2l$ and $\l\in\hZp$.

\proclaim{Proposition \theoremname{\propUlandHpexpansion}}
For $l\in\hZp$, $p\in\Zp$ we have
$$
\aligned
&U_{2l}(\cos\theta)H_p(\cos\theta\mid q^2) = \sum_{m=0}^{\lceil l-1
\rceil} p_m(q^{2p};1,q^{4l-4m};q^2) H_{p+2l-2m}(\cos\theta\mid q^2) \\
& + \sum_{m=0}^{\lceil l-1\rceil} (q^{2p};q^{-2})_{2l-2m}
p_m(q^{2p-4l+4m};1,q^{4l-4m};q^2) H_{p-2l+2m}(\cos\theta\mid q^2) \\
&+'\, p_l(q^{2p};1,1;q^2) H_p(\cos\theta\mid q^2)
\endaligned
$$
where $+'$ means that the last term is absent for $l\in{1\over 2}+\Zp$
and $\lceil a\rceil$ denotes the smallest
integer greater than or equal to $a\in\R$.
Here $U_{2l}$
denotes a Chebyshev polynomial of the second kind, $H_p$ a continuous
$q$-Hermite polynomial, cf. \ssectcontqultraspherical ,
and $p_m$ a little $q$-Jacobi polynomial
\thetag{\vgldeflittleqJacobi}.
\endproclaim

\demo{Proof} From \thetag{\vglmatrixeltalslittleqJacobi}
and \ssectreprsA\ we get for $l\in\hZp$,
$l-n\in\Zp$, $n\geq 0$, the operators, which act on a basis vector $f_p$
by
$$
\aligned
\pi^\infty(t^l_{n,n})f_p & = p_{l-n}(q^{2p};1,q^{4n};q^2)
\sqrt{(q^{2p+2};q^2)_{2n}} f_{p+2n} \\
\pi^\infty(t^l_{-n,-n})f_p & = p_{l-n}(q^{2p-4n};1,q^{4n};q^2)
\sqrt{(q^{2p};q^{-2})_{2n}} f_{p-2n},
\endaligned
\tag\eqname{\vglactionlittleqJaconfp}
$$
with the convention that $f_{-n}=0$ for $n\in\N$. Now we have the
identity
$$
\langle \sum_{n=-l}^l \pi^\infty(t^l_{n,n})f_p, v_\l (N) \rangle =
\langle U_{2l}\bigl( (\a+\d)/2\bigr) f_p, v_\l (N) \rangle =
\langle f_p , U_{2l}\bigl( (\a+\d)/2\bigr) v_\l (N) \rangle ,
\tag\eqname{\vglrightinnerproduct}
$$
since $\a+\d$ is a self-adjoint element of $\A$. Now we use
\thetag{\vglactionlittleqJaconfp} and \thetag{\vglchoppedoffeigenvector}
in \thetag{\vglrightinnerproduct}. For $N\geq p+2l$ the
additional terms on the right hand side of
\thetag{\vglchoppedoffeigenvector} are of no importance and we can use
the explicit value $\langle f_{p\pm 2n}, v_\l(N)\rangle =
(q^2;q^2)_{p\pm 2n}^{-1/2} H_{p\pm 2n}(\l\mid q^2)$.
We obtain the proposition after substitution $m=l-n$.
\qed\enddemo

Instead of a chopped-off eigenvector $v_\l(N)$ of
$\pi^\infty \bigl({1\over 2}(\a+\d)\bigr)$ we can also use the
spectral theorem for the self-adjoint operator
$\pi^\infty \bigl({1\over 2}(\a+\d)\bigr)$, which is a Jacobi matrix
with respect to the basis $f_n$ of $\ell^2(\Zp)$, cf. e.g.
Dombrowski \cite{\Domb}.

Since $U_{2l}(\cos\theta)$ and $H_n(\cos\theta\mid q)$ are
continuous $q$-ultraspherical polynomials for $\b=q$ and $\b=0$,
the result in proposition \propUlandHpexpansion\ can also be
derived from Rogers's linearisation and connection formulas for the
continuous $q$-ultraspherical polynomials, cf. \cite{\AskeI , (4.15),
(4.18)}, \cite{\AskeW , (4.7), (4.8)}.

The limit $q\uparrow 1$ of proposition \propUlandHpexpansion\ to the
relation can be handled using the devices developed by Van Assche and
Koornwinder \cite{\VAsscK , thm.~1} to show that Koornwinder's
addition formula for the little $q$-Legendre polynomials \cite{\KoorAF}
tends to the addition formula for the Lgendre polynomials as $q\uparrow
1$. See also Koelink \cite{\KoelAF , \S 5} for another application of
that theorem.

A straightforward limit sending $q\uparrow 1$ of proposition
\propUlandHpexpansion\ can also be given. Since $\lim_{q\uparrow 1}
H_p(\cos\theta\mid q^2) = 2^p \cos^p\theta$, by either letting $q\uparrow
1$ in the three-term recurrence relation
\thetag{\vgldrietermcontqHermite} for the continuous
$q$-Hermite polynomials or by letting $q\uparrow 1$ in the explicit
expression, cf. \thetag{\vgldefcontqultraspherpols},
for the continuous $q$-Hermite polynomials and using the
binomial theorem, we obtain, after replacing $2l$ by $l$,
$$
U_l(\cos\theta) 2^p\cos^p\theta = \sum_{m=0}^{\lfloor l/2\rfloor}
{}_2F_1 \left( {{-m,l-m+1}\atop{1}};1\right) 2^{p+l-2m}
\cos^{p+l-2m}\theta ,
$$
where $\lfloor a\rfloor$ is the largest integer smaller than or equal to
$a\in\R$.
Applying Vandermonde's formula ${}_2F_1(-n,b;c;1) = (c-b)_n/(c)_n$, cf.
\cite{\ErdeHTF , Vol.~1, \S 2.8(18)}, yields the classical formula for
the Chebyshev polynomial of the second kind, cf. \cite{\ErdeHTF ,
Vol.~2, \S 10.11(23)},
$$
U_l(\cos\theta) = \sum_{m=0}^{\lfloor l/2\rfloor} {{(m-l)_m}\over{m!}}
2^{l-2m} \cos^{l-2m}\theta .
$$

The orthogonality relations for the continuous $q$-Hermite polynomials,
cf. \thetag{\vglorthrelcontqHermite},
lead to the following integral representation of the little
$q$-Legendre polynomials, i.e. little $q$-Jacobi polynomials
\thetag{\vgldeflittleqJacobi} with $a=b=1$.
For $l\in\Zp$ multiply the result of proposition \propUlandHpexpansion\
by $H_p(\cos\theta\mid q^2)$ and integrate over $[0,\pi]$ with respect to
the orthogonality measure $(e^{2i\theta}, e^{-2i\theta};q)_\infty$.

\proclaim{Corollary \theoremname{\corrintegralreplittleqleg}}
The little $q$-Legendre polynomial has the integral representation
$$
p_l(q^p;1,1;q) = {{(q^{p+1};q)_\infty}\over{2\pi}} \int_0^\pi
U_{2l}(\cos\theta) \bigl( H_p(\cos\theta\mid q)\bigr)^2
(e^{2i\theta}, e^{-2i\theta};q)_\infty\, d\theta .
$$
\endproclaim

This integral representation is an alternative for the $q$-integral
representation for the little $q$-Legendre polynomial derived by
Koornwinder \cite{\KoorAF , thm.~5.1}.

%%%%%%%%%%%%%%%%%References%%%%%%%%%%%%%%%%%%%%%%%%%%%%%%%%%%%%%%%%%%%%%%

\enddocument